\documentclass[onefignum,onetabnum]{siamart171218}

\usepackage{ulem}

\usepackage{lipsum}
\usepackage{amsfonts}
\usepackage{graphicx}
\usepackage{epstopdf}

\usepackage{multirow}

\usepackage{tikz}
\usetikzlibrary{matrix}
\usetikzlibrary{shapes,patterns,snakes,arrows}
\usetikzlibrary{positioning,fit,calc}
\usetikzlibrary{graphs}

\usepackage{algorithmic}
\ifpdf
  \DeclareGraphicsExtensions{.eps,.pdf,.png,.jpg}
\else
  \DeclareGraphicsExtensions{.eps}
\fi

\newsiamremark{remark}{Remark}
\newsiamremark{hypothesis}{Hypothesis}
\crefname{hypothesis}{Hypothesis}{Hypotheses}
\newsiamthm{claim}{Claim}

\def\bb{\begin{bmatrix}}
\def\eb{\end{bmatrix}}

\newcommand{\TheTitle}{Learning Adaptive Coarse Spaces Using 
Transferable Neural Network Models for Linear and Nonlinear Overlapping Domain Decomposition Methods}

\headers{Learning adaptive coarse spaces for Schwarz methods}{A. Klawonn, M. Lanser, and J. Weber-Hamacher}

\title{{\TheTitle}}

\author{Axel Klawonn\footnotemark[2]\ \footnotemark[3] \and Martin Lanser\footnotemark[2]\ \footnotemark[3] \and Janine Weber-Hamacher\footnotemark[2] \footnotemark[3]}

\usepackage{amsopn}

\usepackage{ifthen}%
\newcommand{\mycomment}[1]{%
\ifthenelse{\isodd{\value{page}}}{%
\normalmarginpar%
\marginpar{\tiny {#1}}%
}{%
\reversemarginpar%
\marginpar{\tiny {#1}}%
}}%

\def\edge{\mathcal{E}}

\usepackage{cite}
\usepackage[english]{babel}
\usepackage{hyperref} 
\usepackage{amsmath}
\usepackage{amssymb}
\usepackage{amsfonts}
\usepackage[T1]{fontenc}
\usepackage[utf8]{inputenc}

\usepackage{subcaption}

\usepackage{makeidx}  
\usepackage{enumerate}
\usepackage{tikz}
\usetikzlibrary{matrix}
\usetikzlibrary{shapes,patterns,snakes,arrows,decorations.markings,arrows.meta}
\usetikzlibrary{positioning,fit,calc}
\usetikzlibrary{graphs}
\usepackage{ifthen}
\usepackage{nicefrac}
\usepackage{relsize}
\usepackage{bbm}  

\usepackage{listings}  
\usepackage{booktabs} 

\definecolor{darkspringgreen}{rgb}{0.09, 0.45, 0.27}
\definecolor{bostonuniversityred}{rgb}{0.8, 0.0, 0.0}
\definecolor{ferrarired}{rgb}{1.0, 0.11, 0.0}
\definecolor{coralred}{rgb}{1.0, 0.25, 0.25}
\definecolor{pastelred}{rgb}{1.0, 0.41, 0.38}
\definecolor{darkblue}{RGB}{0,0,139}
\definecolor{darkred}{RGB}{105,0,0}

\def\Real{{\rm I\!R}}

\newcommand{\VGDSW}[1]{V_{\rm{GDSW}}}

\newcommand{\takeout}[1]{ }  

\tikzset{%
  every neuron/.style={
    circle,
    draw,
    minimum size=0.6cm
  },
  every input neuron/.style={
    circle,
    draw,
    minimum size=0.6cm,
    fill=green!50
  },
  every output neuron/.style={
    circle,
    draw,
    minimum size=0.6cm,
    fill=orange!30
  },
  every hidden neuron/.style={
    circle,
    draw,
    minimum size=0.6cm,
    fill=blue!40
  },
  neuron missing/.style={
    draw=none, 
    scale=1.5,
    text height=0.3cm,
    execute at begin node=\color{black}$\vdots$
  },
}

\begin{document}

\maketitle

\renewcommand{\thefootnote}{\fnsymbol{footnote}}

\footnotetext[2]{Department of Mathematics and Computer Science, University of Cologne, Weyertal 86-90, 50931 K\"oln, Germany, \email{\{axel.klawonn, martin.lanser, janine.weber\}@uni-koeln.de}, url: \url{http://www.numerik.uni-koeln.de}} 
\footnotetext[3]{Center for Data and Simulation Science, University of Cologne, Germany, url: \url{http://www.cds.uni-koeln.de}} 

\begin{center}
	\today
\end{center}

\begin{abstract}
Domain decomposition methods have been established as efficient and parallel scalable iterative solvers and preconditioners for the  solution of large-scale systems arising from the discretization of partial differential equations. In particular, overlapping Schwarz methods have been successfully applied to a wide range of linear and nonlinear problems. However, for problems with highly heterogeneous coefficients, standard domain decomposition methods typically suffer from deteriorating convergence rates. Robustness with respect to the coefficient contrast can be achieved by enriching the coarse space with adaptively selected constraints obtained from local generalized eigenvalue problems. The construction of these adaptive coarse spaces, however, can account for a significant part of the overall computing time.

In the present work, machine learning techniques are employed to reduce this part of the computing time in the context of the adaptive Generalized Dryja-Smith-Widlund (AGDSW) coarse space. A two-stage approach is proposed in which regression neural networks are used to predict the adaptive coarse basis functions, while a classification neural network is employed to predict the number of basis functions required to ensure robustness. As a consequence, adaptive coarse spaces can be set up in the online phase without solving any eigenvalue problem. Particular attention is paid to problem-specific aspects, including sign-invariant loss functions and post-processing strategies to significantly improve the predicted constraints. The proposed approach is first investigated for scalar diffusion problems with high coefficient contrasts and is subsequently transferred, without retraining, to problems of linear elasticity and to nonlinear $p$-Laplace problems, also within a nonlinear Schwarz framework.  

\end{abstract}



\section{Introduction}
\label{sec:intro}

Discretizing second-order elliptic partial differential equations (PDEs), for example by using finite elements, leads to linear or nonlinear systems of equations which are often very large and thus have to be solved on parallel computers. Parallel scalable, preconditioned iterative methods have been proven to be a good choice for the solution  of these kind of problems. Prominent examples belong to the family of domain decomposition methods (DDMs). In this work, we will focus on a specific problem that occurs when solving such discretized elliptic PDEs iteratively. If the contrast of the maximum and the minimum of certain coefficients of the PDEs becomes too large, usually the convergence rate of standard DDMs will deteriorate. There exist extensions of standard two-level DDMs where the second-level coarse space is enriched with additional information on the coefficient contrast such that the convergence rate becomes robust and the iterative methods converge independently of it. In the present work, we will focus on one member of the family of these robust DDMs, namely the adaptive Generalized  Dryja-Smith-Widlund (AGDSW) overlapping Schwarz method. In the following, we will describe in a broad picture the main ingredient which makes this method robust with regard to a high coefficient contrast, discuss briefly the computational cost, and how AGDSW can be made computationally much more efficient and faster by using machine learning regression and classification approaches. 
To facilitate the flow of reading and present the ideas of this work as concisely as possible, we have initially omitted source references in the following overview. Interested readers can subsequently inform themselves about the current state of research on the individual topics; however, those paragraphs may also be skipped during the first reading.

For the time being, to be more specific, we restrict ourselves to two-dimensional domains, the stationary scalar elliptic and self-adjoint diffusion equation without zero-order term, the equations of linear elasticity, and the nonlinear $p$-Laplace equation. 
As mentioned above, from the meanwhile large family of adaptive DDMs, we choose AGDSW in which certain generalized eigenvalue problems are solved on edges of neighboring subdomains and a selected subset of the resulting eigenvectors is integrated into the second-level coarse space. 
 
We will now give a conceptual representation of a typical AGDSW eigenvalue problem without being too detailed. Let us consider an edge shared by two neighboring subdomains and the finite element stiffness matrix assembled on the union of the two subdomains including the edge. Then, we consider the Schur complement matrix obtained by eliminating all degrees of freedom but those on the edge. Next, we consider the principal minor of the stiffness matrix that belongs to the degrees of freedom on the edge. These two subdomain matrices form a generalized eigenvalue problem. Given a user-defined tolerance, we take all eigenvectors belonging to eigenvalues not larger than this tolerance and include them in the coarse problem defining the second level of the AGDSW Schwarz domain decomposition method. We often denote these eigenvectors as adaptive coarse basis functions or constraints.

The overall computational effort involved in assembling the finite element matrix on the union of the two subdomains, the Schur complement matrix, and for solving the generalized eigenvalue problem, is high. Thus, we consider a machine learning approach based on neural networks to reduce it. Our approach is two-fold. First, we develop regression neural networks to predict the adaptive coarse basis functions. Let us note that we do not have to train these neural networks for each edge in the domain decomposition and for each new decomposition. Under certain assumption explained later in the paper, we can restrict ourselves to the case of two typical adjacent subdomains sharing an edge. Numbering the adaptive coarse basis functions, we have to train a regression neural network for each of them separately. Since usually only a low number of adaptive coarse basis functions is needed to obtain a robust method, we choose this number $k$ a priori; this is problem dependent and the choice has to be made for a certain class of problems but we will see that it can be applied to different linear and nonlinear problems and also to the system of linear elasticity. Most of the times you do not need all $k$ adaptive coarse basis functions on every edge. Thus, in a second step, we train a single classification neural network which predicts the number of adaptive coarse basis functions necessary for robustness for each edge.
Combining these two steps results in a strategy to efficiently predict the adaptive coarse basis functions without solving any eigenvalue problems in the online phase. 
Numerical experiments proving this strategy to work well are first given for the scalar diffusion equation. It is then shown experimentally that the same regression and classification neural networks can also be used for the system of linear elasticity with almost no loss in efficiency and for the nonlinear $p$-Laplace operator using a nonlinear version of the two-level Schwarz DDM using an AGDSW coarse space. Let us note that the neural networks were not newly trained from scratch or retrained when applied to linear elasticity and the $p$-Laplace equation. Finally, we would like to point out that the training of the regression neural networks needed certain problem-aware considerations related to specific aspects coming from the theoretical aspects of DDMs, such as, for example, a post-processing of the adaptive coarse basis functions obtained from the regression neural network, a sign-invariant loss function, or the choice of a specific AGDSW method operating just on a subset of the subdomains.

Learning adaptive coarse constraints and the application of the learned models in linear and nonlinear Schwarz methods brings together three highly relevant research topics: the combination of DDMs and machine learning (ML), adaptive DDMs, and nonlinear DDMs. We provide a brief overview over selected relevant references for all three topics.

{\bf Adaptive DDMs.}
As already said, DDMs are essentially preconditioned iterative methods for solving discretized  PDEs~\cite{toselli_widlund}. They are all based on a divide-and-conquer principle and are robust and highly scalable parallel solvers~\cite{kl_p1,kl_p3,kl_p4,ZampiniTuLi2016,frosch1,SousedikSistekMandel2013}. In addition to the class of nonoverlapping methods, such as FETI-DP (Finite Element Tearing and Interconnecting - Dual Primal)~\cite{fetidp1,fetidp2,fetidp3,fetidp4,Klawonn:2002:DPF_face_constraints} or BDDC (Balancing Domain Decomposition by Constraints)~\cite{bddc1,bddc2}, there is a variety of overlapping variants. Here, we focus on additive two-level Schwarz methods~\cite{toselli_widlund} and two different coarse spaces: The first one is GDSW~\cite{gdsw1,frosch1}, an easy-to-compute coarse space that, however, ignores complex heterogeneous structures in the coefficients of the differential equation. In such scenarios, adaptive coarse spaces such as adaptive GDSW (AGDSW)~\cite{agdsw1,agdsw2,Knepper2022} are still robust. 
For the sake of completeness, it should be noted that there are also alternative adaptive coarse spaces for overlapping methods delivering a comparable robustness; see, for example,~\cite{geneo,shem}. Similar developments exist for nonoverlapping methods such as FETI-DP or BDDC~\cite{SousedikSistekMandel2013,MandelSousedik2007,PechsteinDohrmann2017,KlawonnRadtkeRheinbach2016}.

{\bf Nonlinear DDMs.}
We also consider the solution of nonlinear PDEs. In addition to the classical approach of linearization using Newton's method and applying the DDM to the linearized system, we also examine modern nonlinear DDMs. These typically involve a nonlinear preconditioner based on the principles of domain decomposition, with a wide variety of left-~\cite{nl_left1,nl_left2,nl_left3} and right-preconditioners~\cite{nl_right2,nl_right1}. In particular, we focus here on nonlinear two-level Schwarz methods~\cite{nl_left3,nl_adaptive} and, above all, the combination with AGDSW and our learned coarse spaces. In~\cite{nl_adaptive}, already tests with AGDSW and nonlinear two-level Schwarz have been made. Let us note that there is a much longer list of nonlinear DDMs, but we are mentioning only a small selection here, as the focus of this work is entirely different, namely, on the training and application of learned coarse spaces.

{\bf DDMs and machine learning.}
With the emerging success of scientific machine learning (SciML)~\cite{USDpt:SciML} as an independent research field, in particular, the combination of DDMs and machine learning methods has become of increasing importance. In~\cite{klawonn2024survey}, the authors give a structured overview of publications in this interdisciplinary field up to mid-2024 and suggest a first taxonomy.
Within the scope of the present paper, we especially focus on approaches where different machine learning methods are used to systematically learn coarse spaces in adaptive DDMs. In~\cite{KLW:2024:JCP,KLW:2023:DD_nonlin,HKLW:2018:ML_acc,chung2021learning,Weber:2022:Diss}, regression and classfication neural networks have been trained to predict adaptive constraints in FETI-DP, BDDC, and overlapping Schwarz methods. In~\cite{Dolean:2026:NO}, different neural operators are trained and compared to learn a coarse space in a two-level DDM. For a more detailed review of methods combining machine learning and DDMs, we refer again to~\cite{klawonn2024survey} since a complete list would be beyond the scope of this paper.

\section{GDSW and adaptive GDSW coarse spaces for additive Schwarz methods}
\label{sec:feti-dp_mp}
In the numerical results section, we consider different linear and nonlinear PDEs to challenge our learning-based coarse spaces as well as the classical GDSW~\cite{gdsw1,frosch1} and AGDSW~\cite{agdsw1,agdsw2} coarse spaces. To introduce the construction of all these coarse spaces, we first only consider a linear stationary diffusion problem with coefficient jumps. Let us note that only this linear and scalar problem is used to generate training data for all our learning-based approaches and the same machine learning models are transferred to be used for all other considerd PDEs without retraining.

\subsection{Stationary diffusion}
\label{sec:mp}

 The problem of   \textit{stationary diffusion}, is
\begin{equation}
\begin{aligned}
- \nabla \cdot (\rho \nabla u) & = f \ \text{in}\ \Omega \\
u &= 0\ \text{on}\ \partial \Omega_D \\
\rho \nabla u \cdot n &= g\ \text{on}\ \partial \Omega_N. 
\end{aligned}
\end{equation}  
Here, $\rho: \Omega \subset \mathbb{R}^2 \rightarrow \mathbb{R}$ is a sufficiently smooth coefficient function, $f: \Omega \rightarrow \mathbb{R}$ and $g: \partial \Omega_N \rightarrow \mathbb{R}$ are appropriate right-hand sides, and $n$ denotes the outer unit normal on $\partial \Omega_N$. For the remainder of this paper, we exclusively consider $\partial \Omega_D := \partial \Omega$ for the stationary diffusion case.
We define the Sobolev space $V := H_0^1(\Omega, \partial \Omega_D) := \{ v \in H^1(\Omega ) \, : \, v_{|\partial \Omega_D} = 0 \}$. Thus, for a piecewise constant parameter distribution $\rho \in L^\infty (\Omega)$ with $\rho \geq \rho_{\min} > 0$ and $f \in L^2(\Omega)$, we obtain the weak formulation: 
Find $u \in V$ such that
\begin{equation}
\label{eq:diff_weakForm}
a(u,v) = F(v) \quad \forall v \in V,
\end{equation}
where
\begin{equation}
a(u,v) := 
\int_{\Omega} \rho \nabla u \cdot \nabla v \ dx \quad \text{and} \quad
F(v) := \int_{\Omega} f v\ dx.
\end{equation}
To compute a numerical solution of the given stationary diffusion problem, we discretize~\eqref{eq:diff_weakForm} with finite elements and denote the respective finite element space by $V^h$. We thus obtain a linear system of equations
\begin{equation}
A u = f. 
\label{eq:lin_system}	
\end{equation}
with $u,\, f \in V^h$. 

\subsection{Additive two-level overlapping Schwarz method}
\label{sec:schwarz}

For large problems with many degrees of freedom, a direct solve of (\ref{eq:lin_system}) is infeasible. Instead, an iterative solution of the preconditioned system 
\begin{equation}
M^{-1}A u = M^{-1}f 
\label{eq:lin_system_pre}	
\end{equation}
with the preconditioned conjugate gradient (CG) or generalized minimal residual (GMRES) method is beneficial. Throughout this work, we consider the two-level additive Schwarz preconditioner for $M^{-1}$. Let us briefly describe the Schwarz preconditioner in more detail. We first decompose $\Omega$ into nonoverlapping subdomains $\Omega_i,\;i=1,...,N$, and construct overlapping ones $\Omega_i'$ by adding layers of finite elements around $\Omega_i$, creating an overlap of width $\delta$. Let $R_i : V^h \rightarrow V^h(\Omega_i'),\;i=1,...,N$, be restrictions to the local finite element spaces on the overlapping subdomains.  Additionally, let $V_0(\Omega)$ be a coarse space, spanned by a number of coarse basis functions, which usually have a local support the size of a few subdomains. Discretizing the coarse basis functions in $V_h$ and collecting them in a tall and skinny matrix, we obtain a prolongation operator $\Phi: V_0 \rightarrow V^h$. We will specify how to build $\Phi$ in Section~\ref{sec:gdsw}. 

Finally, the preconditioner writes
\begin{equation}
M^{-1} = \sum\limits_{i=1}^N R_i^T \left(R_i A R_i^T\right)^{-1} R_i + \Phi \left(\Phi^T A \Phi\right)^{-1} \Phi^T,
\label{eq:pre}	
\end{equation}
where $\left(R_i A R_i^T\right)^{-1}$  is each implemented with sparse direct solvers on the overlapping subdomains. This part can be parallelized easily by assigning the subdomains to different parallel processors. The inverse $\left(\Phi^T A \Phi\right)^{-1}$ represents the global coarse problem, which is indispensable for robustness. It is also implemented using a sparse direct solver, but due to its global nature, in the many-subdomain case, a parallel bottleneck. This also implies that one target is always a coarse space that is as small as possible but still robust. This will become important in the numerical results and analysis of different learned coarse spaces.

\subsection{Different coarse spaces: GDSW, AGDSW, and AGDSW-slab}
\label{sec:gdsw}

In general, all considered coarse spaces are built using the nonoverlapping partition with subdomains $\Omega_i,\;i=1,...,N$.  The interface is defined by
\begin{equation} \label{eq:interface}
\Gamma = \bigcup_{i \neq j} \left( \partial \Omega_i \cap \partial \Omega_j \right) \setminus \partial \Omega_D.
\end{equation}
 In all three coarse spaces, the coarse basis functions {$\Phi$} are constructed in two steps: first, the interface values $\Phi_\Gamma$ are defined, and then, the values in the interior of the nonoverlapping subdomains are computed using energy minimizing extensions
\begin{equation} \label{eq:ext}
	\Phi_I = - A_{II}^{-1} A_{I\Gamma} \Phi_\Gamma{.}
\end{equation}
Here, $A_{II}$ and $A_{I\Gamma}$ are submatrices of the global matrix 
\begin{equation*}
	A 
	=
	\begin{bmatrix}
		A_{II} & A_{I\Gamma} \\
		A_{\Gamma I} & A_{\Gamma\Gamma}
	\end{bmatrix},
\end{equation*}
with the indices corresponding to the interior degrees of freedom (DOFs) $I$ and the indices corresponding to the interface DOFs $\Gamma$.  Let us note that the extensions can be computed locally on the subdomains and efficiently in parallel, but we do not discuss the details here for brevity. The different coarse spaces exclusively differ in the definition of $\Phi_\Gamma$. 
In GDSW, $\Phi_\Gamma$ only depends on the decomposition and is therefore coefficient-agnostic. This will be different in the adaptive cases. Let us first partition the interface into edges shared by two subdomains and vertices, which are single nodes shared by more than two subdomains. In all three coarse spaces, one coarse basis function is defined for each vertex $v$. The corresponding column in $\Phi_\Gamma$ is set to $1$ in the vertex and $0$ on the remaining interface. In the case of GDSW, we do a similar thing for edges. For each edge $\edge$, we define a single column in $\Phi_\Gamma$, which is set to $1$ in all DOFs belonging to $\edge$ and $0$ on the remaining interface. By construction, the sum over all columns of $\Phi_\Gamma$ is  a partition of unity on the interface and, as long as $A$ is symmetric positive definite, also the discrete harmonic extension $\Phi_I$ of $\Phi_\Gamma$ builds a partition of unity. 

In the adaptive variants, the coarse basis functions corresponding to edges are built differently, based on generalized eigenvalue problems. We first describe the AGDSW case. Let $\edge$ again be an edge and $\Omega_i$ and $\Omega_j$ the two adjacent nonoverlapping subdomains. Then, we consider the two matrices
\begin{align*}
	S^{(i,j)}_{\edge} 
	:= 
	A^{(i,j)}_{\edge\,\edge} - 
	A^{(i,j)}_{\edge\,\mathcal{R}}
	\big(A^{(i,j)}_{\mathcal{R}\,\mathcal{R}}\big)^{-1}
	A^{(i,j)}_{\mathcal{R}\,\edge}
\end{align*}
and $A^{(i,j)}_{\edge \edge}$, where $A^{(i,j)}$ is the Neumann matrix on $\Omega^{(i,j)} = \overline{\Omega_i \cup \Omega_j}$, $\edge$ corresponds to the DOFs on the (interior) edge $\edge$, and $\mathcal{R}$ corresponds to all remaining DOFs of $\Omega^{(i,j)}$.
Using these matrices, we solve the generalized eigenvalue problem: find $(\tau_{\edge}, \mu_{\edge}) \in V_0^h (\edge) \times \Real$ such that
\begin{align}
\label{eq:agdsw:evp}
	\theta^T \, S^{(i,j)}_{\edge} \, \tau_{\edge} 
	= 
	\lambda_\edge^{-1} \ 	
	\theta^T \, A^{(i,j)}_{\edge\, \edge} \, \tau_{\edge} 
	\quad \forall \theta \in V_{0}^{h}\left(\edge\right).
\end{align}
Here, $V_{0}^{h}\left(\edge\right)$ is the finite element {space} on the {interior nodes of the} edge. 
Let the eigenvalues $\lambda_\edge$ of the eigenvalue problem~\eqref{eq:agdsw:evp} be sorted in ascending order. Then, we select all eigenpairs {$(\lambda_\edge,\tau_\edge)$} with $\lambda_\edge$ below a user-chosen tolerance $tol$.  Each of the $\tau_\edge$  then defines one column in $\Phi_\Gamma$. Let us note that the smallest eigenvalue is always zero and the corresponding eigenvector is exactly the edge constraint of GDSW, that is, the constant function on the edge. 

In the literature, efficient variants of adaptive coarse spaces avoiding the expensive computation of $S^{(i,j)}_{\edge}$ are called economic or slab variants. We denote the variant used here by AGDSW-slab; see~\cite{Knepper2022,agdsw1} for further details. Instead of setting up $A^{(i,j)}$ on the union $\Omega^{(i,j)} = \overline{\Omega_i \cup \Omega_j}$, we consider a smaller area $\Omega_s^{(i,j)} = \overline{(\Omega_i \cup \Omega_j) \cap \Omega_\edge}$. Here, the slab $\Omega_\edge$ is an area in the neighborhood of the edge $\edge$, usually a few layers of finite elements. Consequently, the matrix $A_s^{(i,j)}$ with Neumann-type boundary condition on $\partial \Omega_s^{(i,j)}$ is smaller then $A^{(i,j)}$ and the computation of the Schur complement on the interior of the edge is much cheaper. The eigenvalue problem is nearly identical, only the Schur complement on the left hand side is built from the smaller matrix defined on the slab. The only negative effect of using AGDSW-slab instead of AGDSW is the risk of slightly larger coarse spaces, since a higher number of small eigenvalues might occur. Nevertheless, we will see that AGDSW-slab is better suited for learning approaches which replace the computation of adaptive constraints later on.   

\subsection{Condition number bounds}
\label{sec:cond}

For GDSW, the classical two-level estimate reads
\[
\kappa(M^{-1}A)\le C\Bigl(1+\tfrac{H}{\delta}\Bigr)\Bigl(1+\log\!\Bigl(\tfrac{H}{h}\Bigr)\Bigr)^2,
\]
where $H$ and $h$ are the coarse and fine mesh sizes, respectively, $\delta$ the overlap width, and $C$ depends on the coefficient contrast between the highest and lowest diffusion coefficient $\rho$ occuring but not on the mesh resolution itself.
For AGDSW, the local generalized eigenvalue problems are solved on each interface component with tolerance $tol$, and the preconditioned operator satisfies a bound of the form
\[
\kappa(M^{-1}A)\le C\Bigl(1+\tfrac{1}{\mathrm{tol}}\Bigr)^2,
\]
where $C$ is independent of $H$, $h$, and the coefficient contrast; the logarithmic mesh terms are replaced by dependence on the user-prescribed EVP tolerance. 

These bounds hold for stationary diffusion as well as linear elasticity problems and the AGDSW estimate also covers AGDSW-slab, where only the adaptive coarse space might be slightly larger to obtain a similar condition number. For details and proofs, we refer to~\cite{agdsw1,agdsw2,Knepper2022}.

\section{Learning adaptive constraints using supervised regression neural network models}
\label{sec:ml_nn}

 As already stated in the introduction, the setup and solution of the local eigenvalue problems in the setup of AGDSW-type coarse spaces is quite costly and can take up more than 50\% of the runtime in a parallel implementation. In the following, we describe how to accelerate our computations replacing this computationally expensive step by SciML approaches. 
To obtain an efficient, robust, and trustworthy adaptive Schwarz method, our aim is to directly learn and predict the adaptive coarse spaces in Schwarz methods completely by supervised regression neural network models, that is, to predict $\Phi_\Gamma$ on crucial edges. 
This extends our work from~\cite{KLW:2024:JCP,HKLW:2018:ML_acc,KLW:2023:DD_nonlin, Weber:2022:Diss} on adaptive FETI-DP as well as preliminary work in~\cite{HKLW:2023:DD_schwarz} on linear and nonlinear Schwarz. 
In particular, in this work, we significantly modify our machine learning approach in comparison to previous work by making the method more problem-aware, that is, we use a modified loss function, a slab-variant to compute the training data, and we combine a classifier and a regression neural network model to completely automate the creation of the adaptive coarse space.

\subsection{Defining a mesh-independent neural network regression model}
\label{sec:nn_model}

\begin{figure}
\centering
\includegraphics[width=0.95\textwidth]{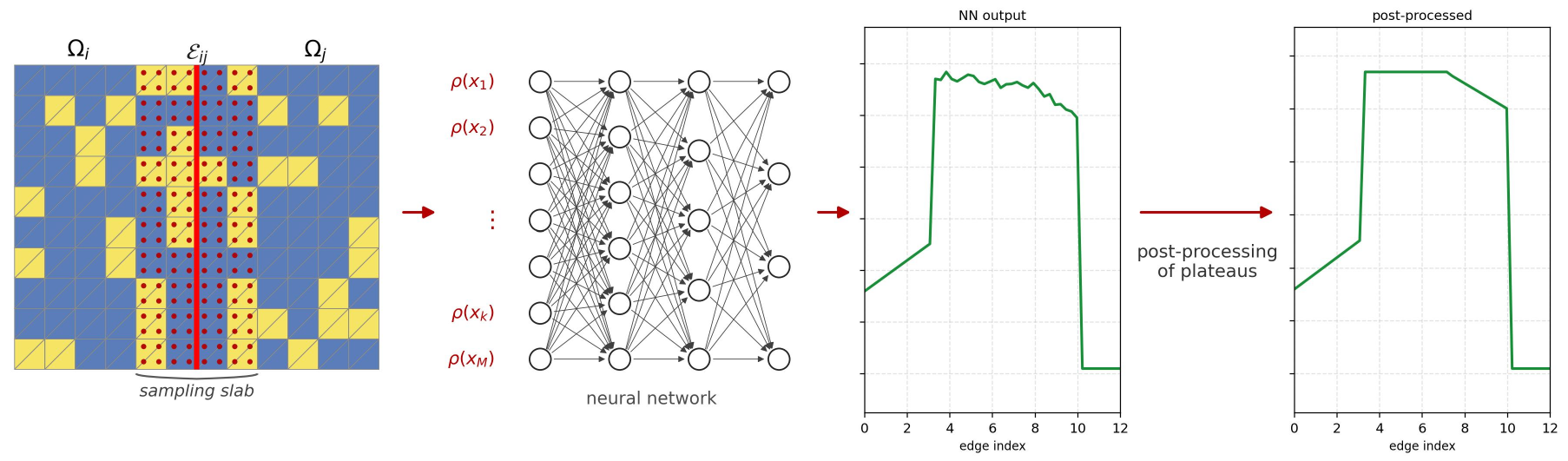}
\caption{Visualization of our network models $N_l$ to predict discrete approximations of the AGDSW edge constraint. Dark blue corresponds to a low coefficient and yellow corresponds to a high coefficient. The red samples within slabs are used as input data for a feedforward neural network (left and middle). A plateau post-processing is applied to the approximated edge constraints obtained by the neural network (middle and right). }
\label{fig:fnn}
\end{figure}

For the prediction of discretized adaptive edge constraints for an edge $\edge$, we use $k \in \mathbb{N}$ separate regression network models to predict $k$ edge constraints resulting from the edge eigenvalue problem~\eqref{eq:agdsw:evp} of AGDSW; see also Fig~\ref{fig:fnn} for a schematic overview. In the following, we denote each of these regression neural network models by $N_l,\ l\leq k$. Note again that in AGDSW, the first adaptive constraint is always a constant function, and, hence, has not to be predicted. Thus, our first regression neural network, that is, $N_1$ always predicts the second constraint, $N_2$ predicts the third constraint etc. 

As input data for the neural network models $N_l,\ l\leq k, $ we aim at a mesh-independent image representation of the coefficient function $\rho$ of the two subdomains $\Omega_i$ and $\Omega_j$ sharing an edge $\edge$. In principle, we use the same sampling procedure as described in detail in~\cite{KLW:2024:JCP} for adaptive FETI-DP. That is, we compute function evaluations of the coefficient function with a consistent ordering and a fixed length as input data for the neural networks. 
This is achieved by computing a geometric grid of points which covers the coefficient function in the two neighboring subdomains of an edge $\edge$. 
However, in contrast to our previous work on predicting the coarse basis functions in FETI-DP~\cite{KLW:2024:JCP}, as an alternative, we choose an economic variant here for the sampling, which closely aligns with the AGDSW-slab variant. Hence, instead of using the sampling, that is, the coefficient values within the entire subdomains $\Omega_i$ and $\Omega_j$ as input data, we exclusively use the coefficient values in a smaller area $\Omega_s^{(i,j)} = \overline{(\Omega_i \cup \Omega_j) \cap \Omega_\edge}$, where the slab $\Omega_\edge$ is the area around the edge $\edge$, just as in AGDSW-slab. 
The intuition behind this is two-fold. First, we know from our experience with adaptive coarse spaces, that usually the area close to an edge is the most important for the decision whether adaptive constraints are necessary for robustness, and, additionally, the coefficient jumps in the close neighborhood of the edge are the most important to determine the actual shape of the coarse basis function. Second, our hope is that the regression neural networks can focus more on the relevant, significant information during the learning process, given less input data only on smaller slabs, to predict the approximate shape of the coarse basis functions, rather than being misguided by potentially unnecessary information on the coefficient function far away from the edge $\edge$. 

As output for the different network models $N_l, l \leq k,$ we use discrete values of the adaptive edge constraints obtained by the solution of the local edge eigenvalue problem~\eqref{eq:agdsw:evp}. More precisely, for the training of the $l$-th network $N_l$, we use the respective constraint as output data which results from the eigenvector $\tau_{\edge}$ belonging to the eigenvalue $\mu_{\edge}$.

The architecture of the regression neural network models have been optimized by a grid search and using cross validation on the training and validation data. The obtained optimized parameters are summarized in Table~\ref{tab:hyperparams}.
Let us note that, in contrast to our previous work on FETI-DP, here, we always use a sign-invariant MSE (sign-inv-MSE) loss function for the regression neural networks which is defined as follows
\begin{equation}
\min \left(\text{MSE}(y_{\text{true}},y_{\text{pred}}), \text{MSE}(y_{\text{true}},-y_{\text{pred}}) \right),
\end{equation}
where MSE denotes the classical mean squared error, $y_{\text{true}}$ the ground truth, and $y_{\text{pred}}$ the prediction by the regression neural network. 
We have deliberately chosen the sign-inv-MSE in this study, since we have experienced in AGDSW and also adaptive FETI-DP that mirrored adaptive constraints are also robust and do not change the iteration count or condition number of the adaptive DDM. 
Therefore, we have strategically also included this information directly in the network learning problem by choosing a problem-aware loss function. In preliminary tests on simple coefficient distributions with channels of different widths, we also saw a massive improvement compared to using a simple MSE loss; we do not include these results here for brevity and all our models  use the sign-inv-MSE loss.
For the implementation of the sign-inv-MSE and the training of the different network models we use TensorFlow-GPU 2.5~\cite{tensorflow2015-whitepaper}.

\begin{table}[t]
\centering
\begin{tabular}{l|l}
{\bf Hyperparameter} 	& {\bf Optimal choice} \\
\hline
\# Hidden layers 							& 3 \\
\# Neurons per hidden layer					& \{150, 100, 80\} \\
Dropout per hidden layer & 5\% \\
Activation function & GeLU \\
Optimizer			& Adam \\
Initial learning rate &  0.001 \\
Loss function & sign-inv-MSE \\
\end{tabular}
\caption{\label{tab:hyperparams} Hyperparameters for the regression neural networks obtained by a grid search.}
\end{table}

\subsection{Results on training and validation data}
\label{sec:nn_traindata}

As mentioned before, the training of the regression neural network models is completely localized, such that, for the prediction of the adaptive constraint for an edge $\edge$ we exclusively consider the coefficient distributions in the two neighboring subdomains $\Omega_i$ and $\Omega_j$ sharing the edge $\edge$. 
For the training and validation of the network models $N_l, l=1,\ldots, k,$ we have generated $36\,000$ datapoints consisting of pairs of neighboring subdomains sharing an edge with varying coefficient distributions for a mesh size of $H/h=20$; Here, $H$ is the width of a square subdomain and $h$ the width of a finite element.   
To generate an appropriate amout of training data, where also a sufficient amount of edges is present, with a higher number of adaptive constraints necessary for robustness, we use a set of randomized coefficient distributions where we also impose a slight structure in the coefficients. 
 For the generation of this training data set, we create piecewise-constant coefficient fields on a structured reference grid with $H/h=10$, and then interpolate them to the fine finite element mesh with $H/h=20$. We choose $\rho=1$ as the low coefficient and $\rho=1e6$ in the high regions. In general, the patterns consist of horizontal and/or vertical strips.

For each random instance, we draw two integer strip lengths $\ell_h,\ell_v \in \{0,1,\ldots,10\}$ and two strip probabilities $p_h,p_v \in [0.05,0.25]$ uniformly at random. If one length is smaller than 4, it is set to zero. This enforces the occurance of distributions with vertical or horizontal orientations only. If both lengths are smaller than $4$, we enforce that at least one orientation is active by setting one length to $0$ and the other to $4$. Thus, each field is either horizontally striped, vertically striped, or combines both orientations.

\begin{figure}[ht]
\centering
\includegraphics[width=0.45\textwidth]{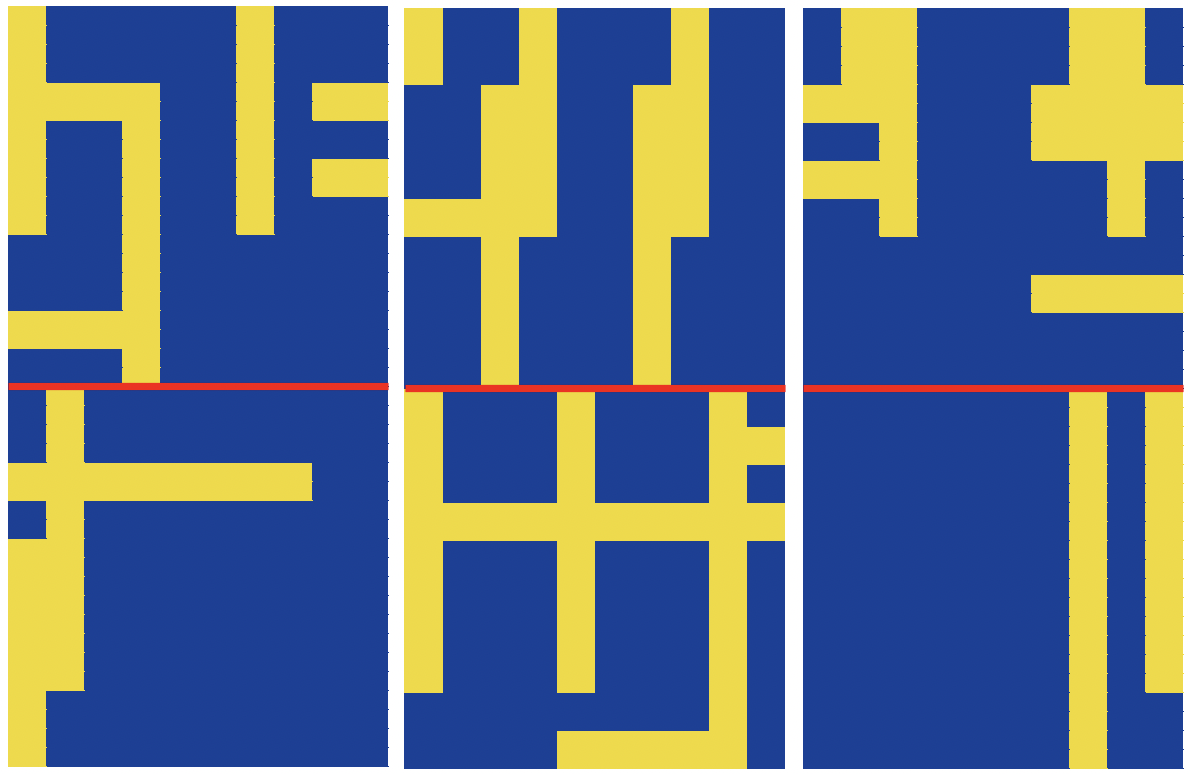}
\includegraphics[width=0.45\textwidth]{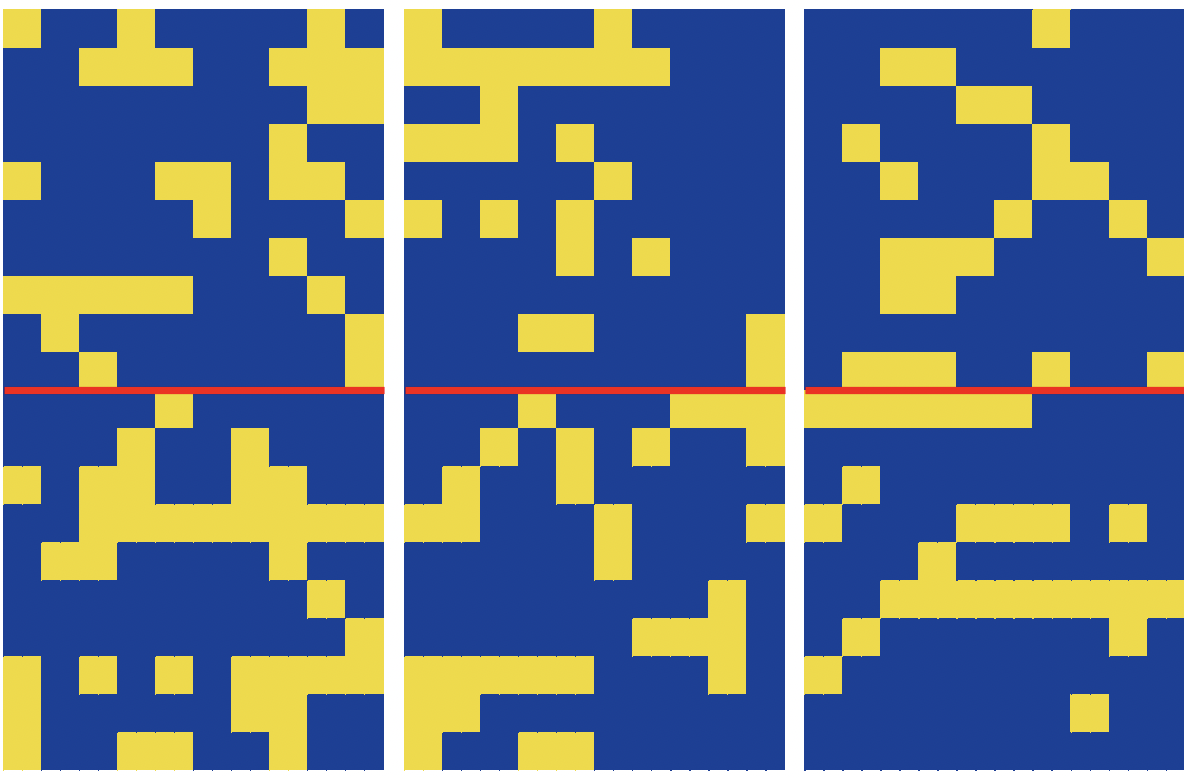}
\caption{Exemplary training data with a smart random coefficient distribution.}
\label{fig:train_data}
\end{figure}

The horizontal pattern is built by partitioning the reference grid with $H/h=10$ into consecutive horizontal strips of width $\ell_h$; each strip is independently declared as high-contrast with probability $p_h$. The vertical pattern is constructed analogously with $\ell_v$ and $p_v$, and the two patterns are combined by a logical OR (a cell is high-contrast if it is high in either orientation). This yields structured, random heterogeneous media rather than pixel-wise independent coefficients. We refer to this set of training data as \textit{smart-random data}.
Six exemplary coefficient distributions generated by this technique are visualized in Fig.~\ref{fig:train_data}.
For all training data configurations, we always set the high coefficient to $\rho_1 = 1e6$ in the yellow pixels and $\rho_2=1$ otherwise.

To obtain the ground truth and output data for the network models, we solve the eigenvalue problem of AGDSW and AGDSW-slab for all $36\,000$ data configurations and save the number of necessary constraints for robustness as well as the concrete adaptive edge constraints for both coarse spaces. For the selection of the constraints we have used the tolerance $tol=0.01$ in both cases  and for AGDSW-slab the width of the area around the edge is two rows of finite elements in both directions. An overview of the number of necessary constraints required for the different coefficient distributions within the training and validation data is summarized in Table~\ref{tab:no_train}.

To monitor the generalization properties of the trained networks, the available data to train each network $N_l,l\leq k,$ have always been split into $70\%$ training and $30\%$ validation data. 

\begin{table}
\centering
\begin{tabular}{l|c|r|r}
 &  & no. train. data& no. train. data \\
class & no. constraints &AGDSW & AGDSW-slab \\\hline
class 0 & constant GDSW function & - & - \\
class 1 & +1 constraint &$19\,957$& $22\,934$ \\
class 2 & +2 constraints &$7\,636$& $10\,480$ \\
class 3 & +3 constraints &$1\,531$& $2\,229$ \\
class 4 & > 3 constraints & - & - \\
\end{tabular}
\caption{Overview of the number of training data for the prediction of the first, second, and third constraint of AGDSW. }
\label{tab:no_train}
\end{table}

\begin{figure}[ht]
\centering
\begin{subfigure}[t]{0.45\textwidth}
\centering
\includegraphics[width=0.9\textwidth]{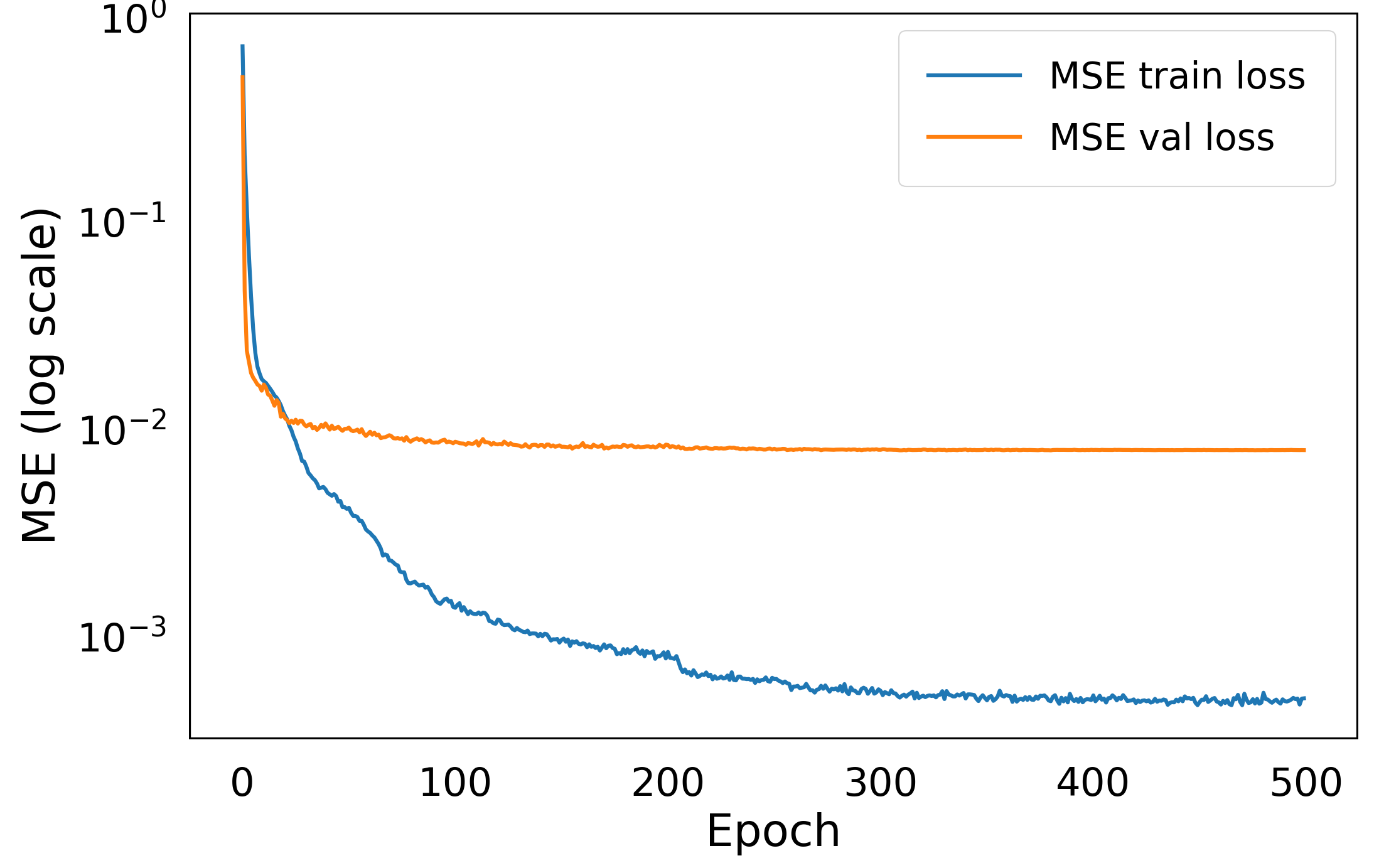}
\caption{Loss curve for training with data from AGDSW, resulting in train-loss=$5.18e$-$04$, val-loss=$8.22e$-$03$. }
\end{subfigure}
\begin{subfigure}[t]{0.45\textwidth}
\centering
\includegraphics[width=0.9\textwidth]{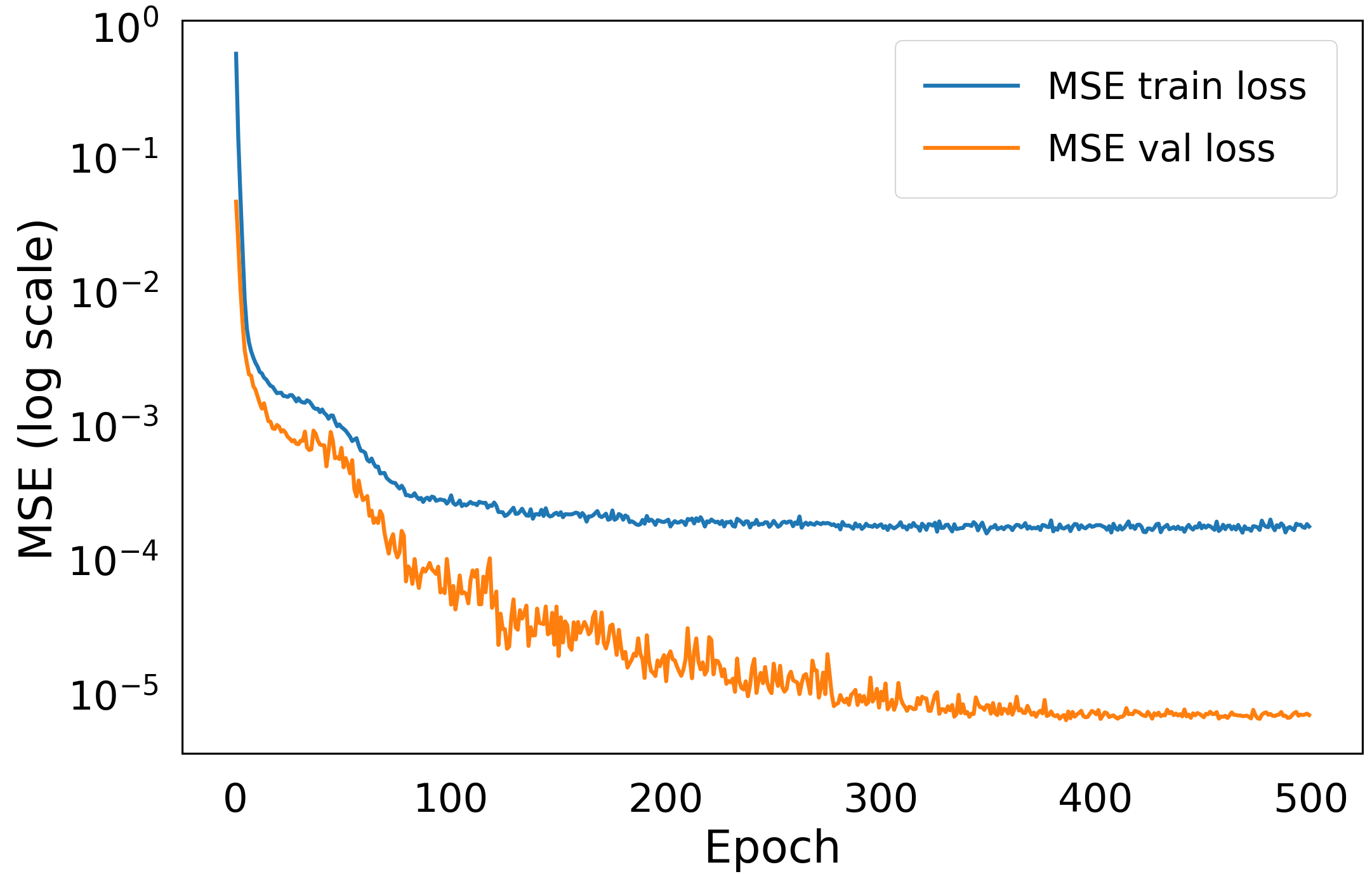}
\caption{Loss curve for training with data from AGDSW-slab, resulting in train-loss=$1.99e$-$04$, val-loss=$7.89e$-$06$.}
\end{subfigure}
\caption{Comparison of different loss curves for training and validation loss for the regression model to predict the second constraint for different AGDSW variants, i.e., full AGDSW and AGDSW-slab with a slab of two rows of finite elements around the edge.}
\label{fig:loss_curves}
\end{figure}

\subsection{Classifier for the number of constraints}
\label{sec:clf}

With the aim to fully automate the construction of an AGDSW coarse space, in this work, we combine the prediction of adaptive edge constraints with a classifier to predict the number of necessary constraints for the first time. 
In order to train a suitable classification neural network, we extend our previous work~\cite{HKLW:2018:ML_acc} on adaptive FETI-DP. 
Here, we train a classification neural network using the same $36\,000$ datapoints as described in Section~\ref{sec:nn_traindata}. In extension to our original results in~\cite{HKLW:2018:ML_acc}, where we exclusively differentiated between two or three classes of edges, respectively, we now extend our neural network-based surrogate model to categorize the following five classes of edges: 

\begin{itemize}
\item[i)] class 0: only constant GDSW function necessary for robust coarse space
\item[ii)] class 1: +1 necessary constraint
\item[iii)] class 2: +2 necessary constraints
\item[iv)] class 3: +3 necessary constraints
\item[v)] class 4: $>3$ additional necessary constraints.
\end{itemize}

The distribution of the five classes within the $36\,000$ datapoints is summarized in Table~\ref{tab:no_train}. To control the generalization properties of the classification model, these data have been split into $70\%$ training and $30\%$ validation data. To mitigate the unequal distribution of the training and validation data among the five classes, class weights were computed using the inverse-frequency weighting scheme and incorporated during the training process, thereby assigning greater importance to underrepresented classes in the loss function. 

\begin{figure}
\centering
\includegraphics[width=0.98\textwidth]{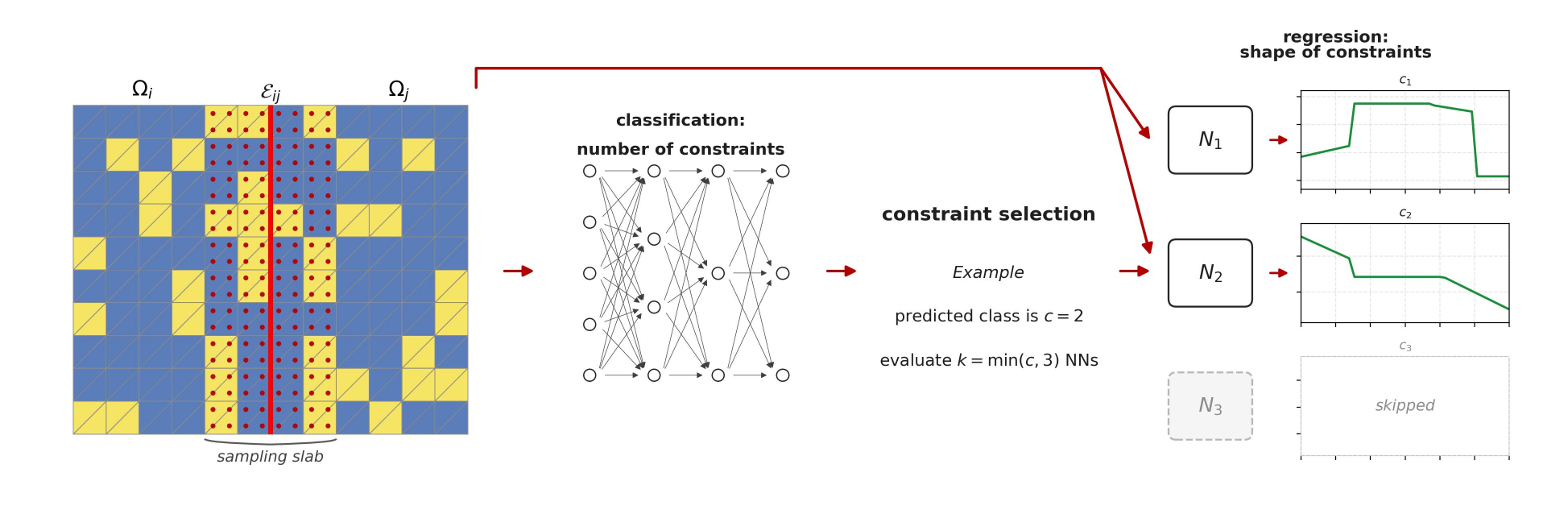}
\caption{Visualization of our SciML model for the prediction of discrete adaptive edge constraints in LAGDSW. First, a classifier predicts the number of necessary coarse basis functions for the considered edge (middle). Depending on the predicted class $c$, the first $c$ regression neural networks $N_l, l=1,\ldots, c$ are evaluated to predict the respective edge constraints (right). }
\label{fig:ml-adaptive-blackbox}
\end{figure}

The trained neural network classifier consists of a fully connected feedforward neural network with three hidden layers of $512$, $256$, and $128$ neurons, respectively, with the architecture optimized via a Bayesian gridsearch. For all hidden layers, the GELU activation function is used, as well as batch normalization, and a dropout rate of $0.3$ after each hidden layer. The output layer of the network comprises five neurons with a softmax activation function corresponding to the five assigned  classes. A schematic representation of our resulting SciML model, combining the classifier with the regression networks defined in sect.~\ref{sec:nn_model}, is shown in Fig.~\ref{fig:ml-adaptive-blackbox}.

\begin{figure}[ht]
\centering
\begin{subfigure}[t]{0.48\textwidth}
\includegraphics[width=0.95\textwidth]{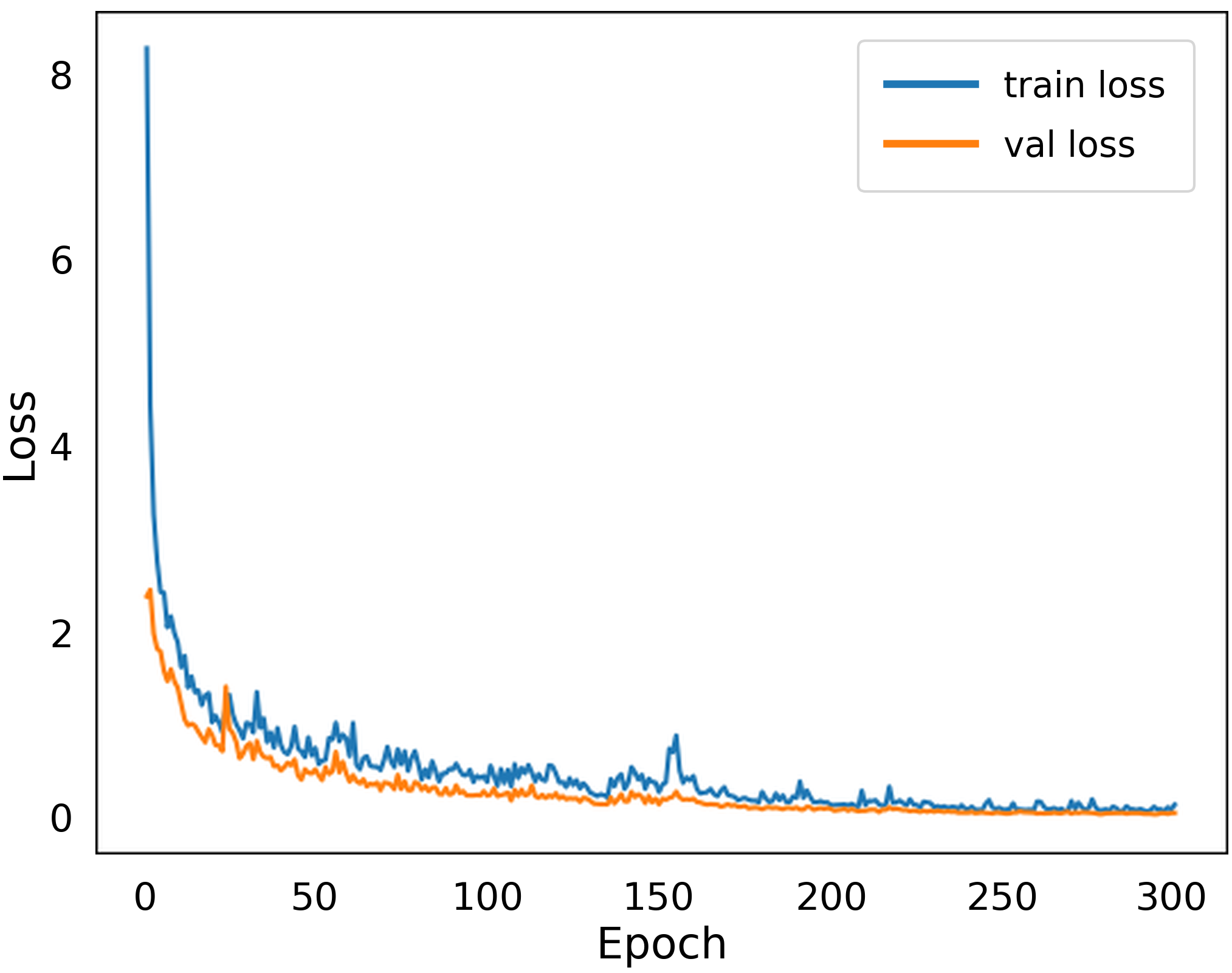}
\caption{Loss curve for training and validation data from AGDSW-slab.}
\end{subfigure}
\hfill
\begin{subfigure}[t]{0.48\textwidth}
\includegraphics[width=0.9\textwidth]{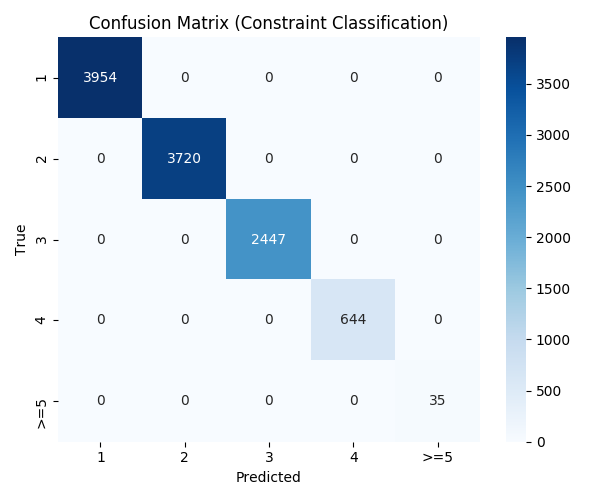}
\caption{Confusion matrix for validation data for AGDSW-slab. }
\end{subfigure}
\caption{Loss curves for training and validation data as well as confusion matrix for the validation data for the classification neural network for AGDSW-slab with a slab of two rows of finite elements around the edge. }
\label{fig:res_class}
\end{figure}

The loss curves during training for the training and validation data as well as a multi-class confusion matrix for the resulting classification model are visualized in Fig.~\ref{fig:res_class}. As we can observe from Fig.~\ref{fig:res_class}, the classification model shows a very accurate performance with respect to the validation data and can efficiently distinguish between the five defined classes of edges.

\section{Numerical results for stationary diffusion}
\label{sec:num}

Before we start to present numerical results for a robust and flexible learned coarse space without solving any eigenvalue problem, we first discuss two necessary modifications. For simplicity, we discuss both in detail only for the second edge constraint, that is, the first constraint that is not constant. Later, we use the best configuration for all trained models  $N_l$ for all constraints  and, additionally, in combination with the classifier. The first modification, as already mentioned, is to use a slab-variant of AGDSW to train our coarse space. As mentioned earlier, the advantage is a smaller input size compared to a full sampling of both subdomains and thus a simpler learning task. The second, less obvious modification, is a post-processing routine which is necessary to obtain a robust solver. In the next subsection, we motivate and describe the post-processing, before we afterwards find our best combination of options.

\subsection{Noise affecting adaptive constraints and post-processing}

\begin{figure}
\centering
\includegraphics[width=1.0\textwidth]{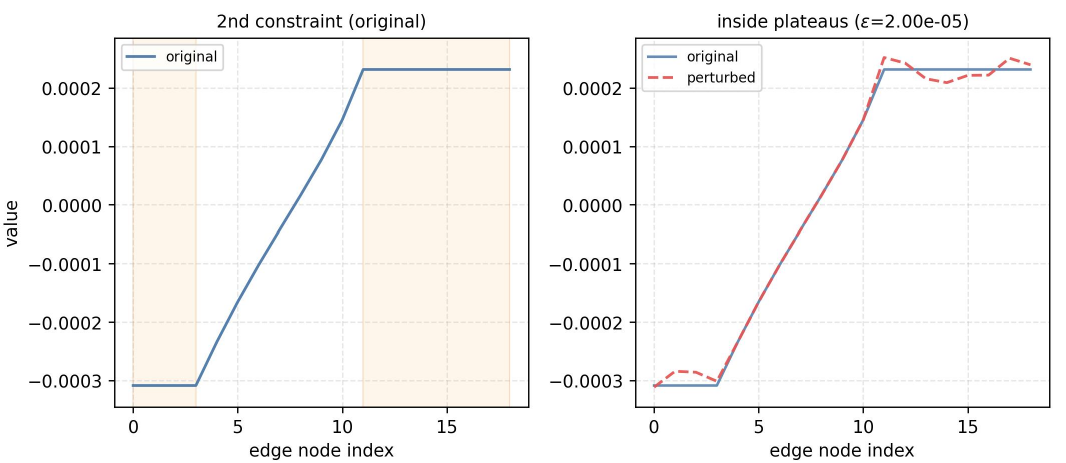}\\
\includegraphics[width=1.0\textwidth]{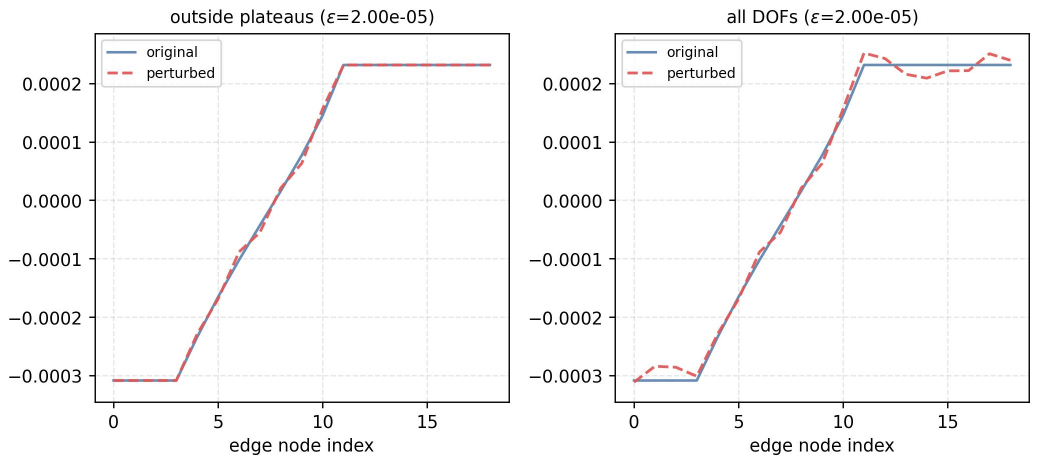}
\caption{{\bf Top left:} Example of an original AGDSW (slab) constraint with marked plateaus; {\bf top right:} perturbed constraint adding noise within plateaus; {\bf bottom left:} perturbed constraint adding noise outside the plateaus; {\bf bottom right:} perturbed constraint adding noise everywhere.}
\label{fig:noise_cosntraint}	
\end{figure}

\begin{figure}
\centering
\includegraphics[width=0.75\textwidth]{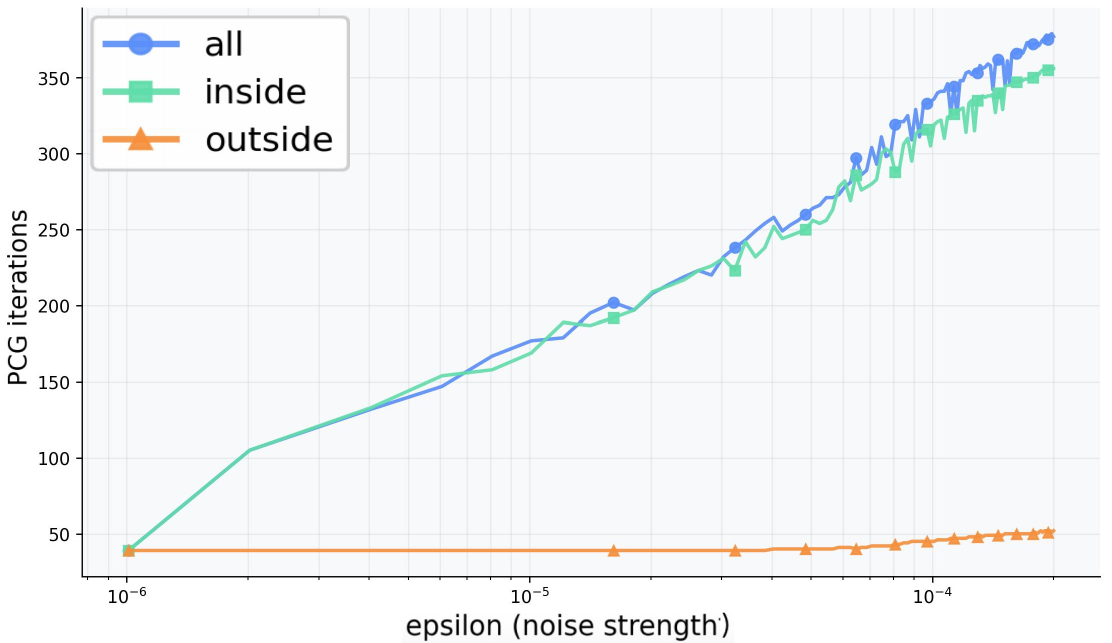}
\caption{Iteration count versus noise for all three modes: perturbation within plateaus, outside of them, and everywhere.}
\label{fig:noise}	
\end{figure}

In general, constraints coming from the AGDSW or AGDSW-slab eigenvalue problem tend to have plateaus, that is, parts with a constant value. Often, these plateaus are aligned with high coefficients. A general problem of neural network predictions using an MSE-type loss is that they tend to be smooth with sometimes larger curvature or sometimes to have noise in a very small range. Consequently, reproducing the plateaus exactly is a hard task for a machine learning model, since the gradient of noisy or smooth functions is usually not zero. Unfortunately, it seems to be crucial to have a zero gradient within the plateau areas. There is currently no theory which supports this statement, but we provide some numerical evidence in the following. Therefore, we take all adaptive constraints obtained by solving the associated eigenvalue problems, for an example with 25 subdomains (stationary diffusion with coefficient jumps) and add some weighted random noise to the adaptive constraints. We consider three different perturbations: only adding noise in the plateaus, only adding noise outside the plateaus, and adding it everywhere. An exemplary adaptive constraint and the three noise modes are shown in Fig~\ref{fig:noise_cosntraint}. Increasing $\epsilon$, where $\epsilon$ is the weighting factor of the noise, deteriorates convergence of PCG with an additive Schwarz preconditioner and AGDSW coarse space quite drastically and already for very small perturbations, at least, if the noise is added inside of the plateaus or everywhere; see Fig.~\ref{fig:noise}. In contrast, if noise is only added outside the plateaus, the iteration count stays nearly constant. Consequently, the most important task for our prediction models is to accurately predict the plateaus and the right constant value within. 

\begin{figure}
\centering
\includegraphics[width=1.0\textwidth]{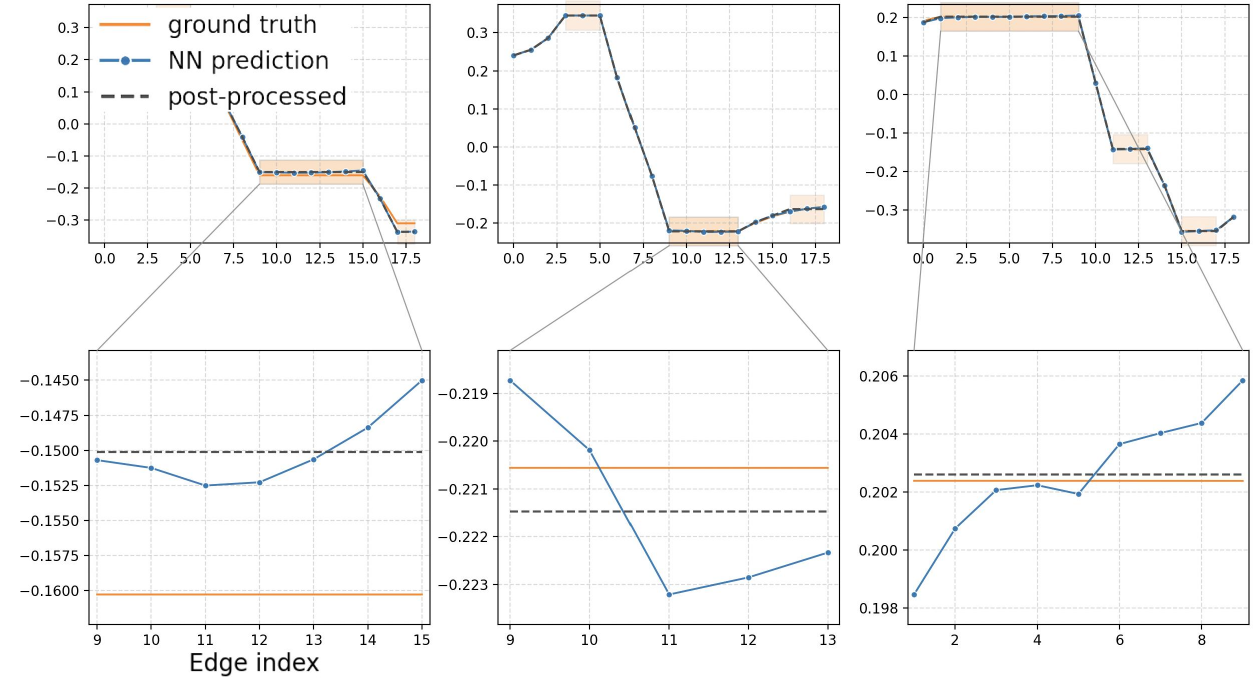}
\caption{Three different adaptive constraints (ground truth), the neural network prediction, and the post-processed neural network prediction reconstructing the plateaus; bottom row shows a zoom into the largest found plateau.}
\label{fig:zoom}	
\end{figure}

If we now look at the prediction of AGDSW-slab constraints, the plateaus seem to be well captured; at least, there is no visible noise here. However, if you zoom in, you can see small deviations from the ground truth; see Fig.~\ref{fig:zoom}. In order to remove this noise or, in other words, to reconstruct exactly constant plateaus, we decided to run a simple post-processing procedure on all neural network outputs before including them into the coarse space. 

\paragraph{Plateau post-processing.}
Let $\mathbf v=(v_0,\ldots,v_{m-1})^{\top}$ be the neural network prediction on an edge. We split $\{0,\ldots,m-1\}$ into maximal contiguous index sets (runs) $I$ such that, for each $k\in I\setminus\{m-1\}$,
\[
|v_{k+1}-v_k|\le \varepsilon_{\mathrm{abs}}.
\]
If $|I|\ge 2$, replace $(v_k)_{k\in I}$ by the constant $\bar v_I:=|I|^{-1}\sum_{k\in I} v_k$; shorter runs are unchanged. The effect of the post-processing can also be seen in~Fig.~\ref{fig:zoom}. 

Of course, this procedure adds a further hyper-parameter $\varepsilon_{\mathrm{abs}}$ which has to be chosen. Since we already have a large amount of training data that we used to train the regression neural networks, we can use this data to optimally tune $\varepsilon_{\mathrm{abs}}$. Generally speaking, if the parameter is too small, plateaus will likely not be fully identified. If the parameter is too large, the plateaus will appear larger than in the ground truth, or two or even more plateaus that are actually separate in the ground truth will be merged. The first effect is obviously the worst, since it introduces non-zero gradients in plateaus. On the other hand, plateaus which are slightly too large might be tolerable. Consequently, we performed the following data analysis: We tested an increasing sequence of values of $\varepsilon_{\mathrm{abs}}$ on the neural network predictions for all available training data until 100\% of the plateaus in the ground truth data were covered by those found in the predictions during post-processing. We did this separately for all three trained AGDSW-slab networks $N_1$ to $N_3$, specialized on the second, third, and fourth adaptive constraints, and their specific training data. The three epsilon values found are $\varepsilon_{\mathrm{abs}}=0.024$ for the second constraint, $\varepsilon_{\mathrm{abs}}=0.034$ for the third one, and finally $\varepsilon_{\mathrm{abs}}=0.071$ for the third. Throughout the rest of the paper we use exactly these data-motivated parameters in the post-processing. We also tested different choices varying between 0.01 and 0.1, which performed slightly worse on average, but still achieved decent results. The choice of this hyperparameter is therefore not that critical for the overall performance.

\subsection{Best possible combination: Slab-variant and post-processing}
\begin{figure}[h!]
\centering
\includegraphics[width=0.8\textwidth]{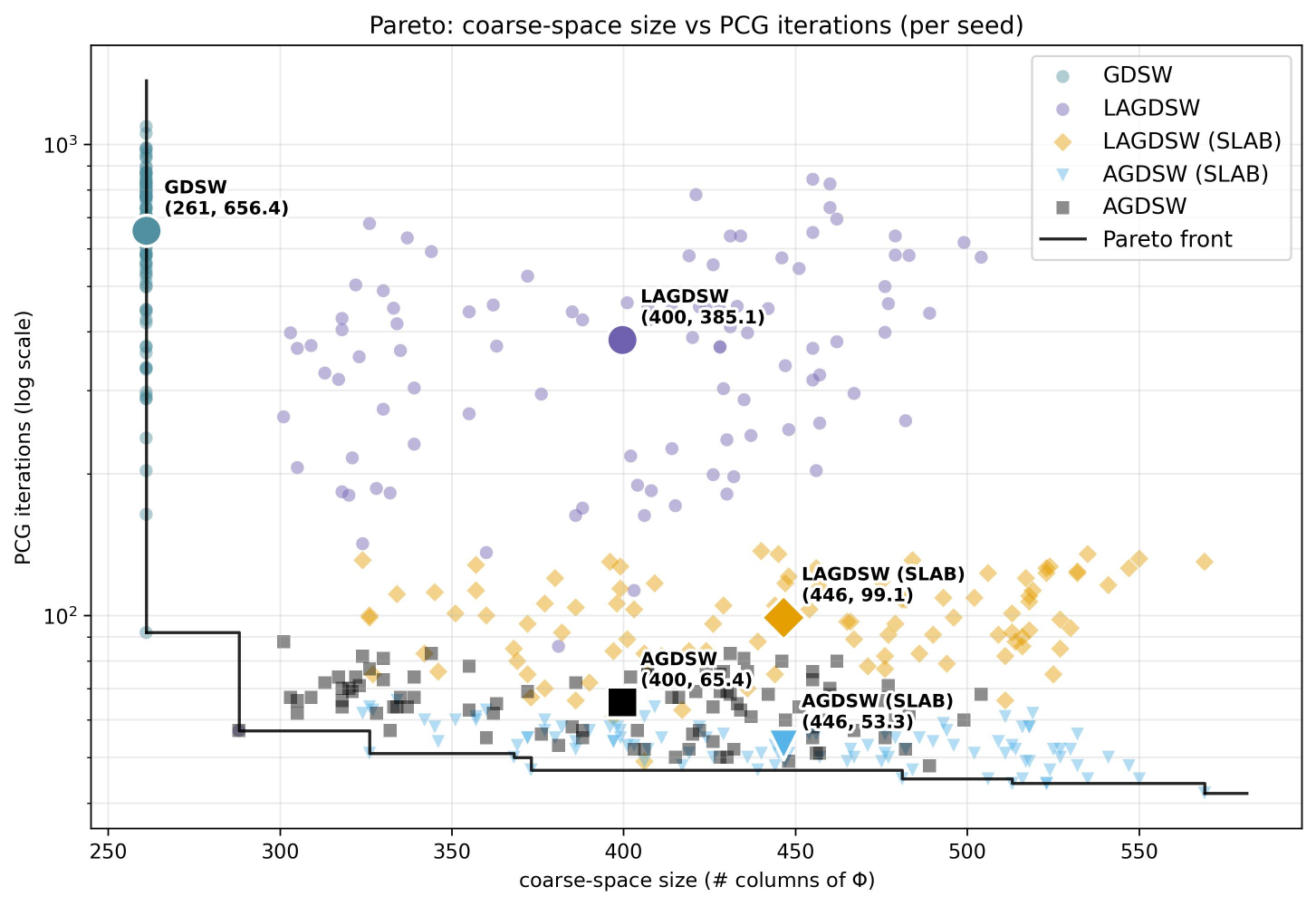}	
\includegraphics[width=0.8\textwidth]{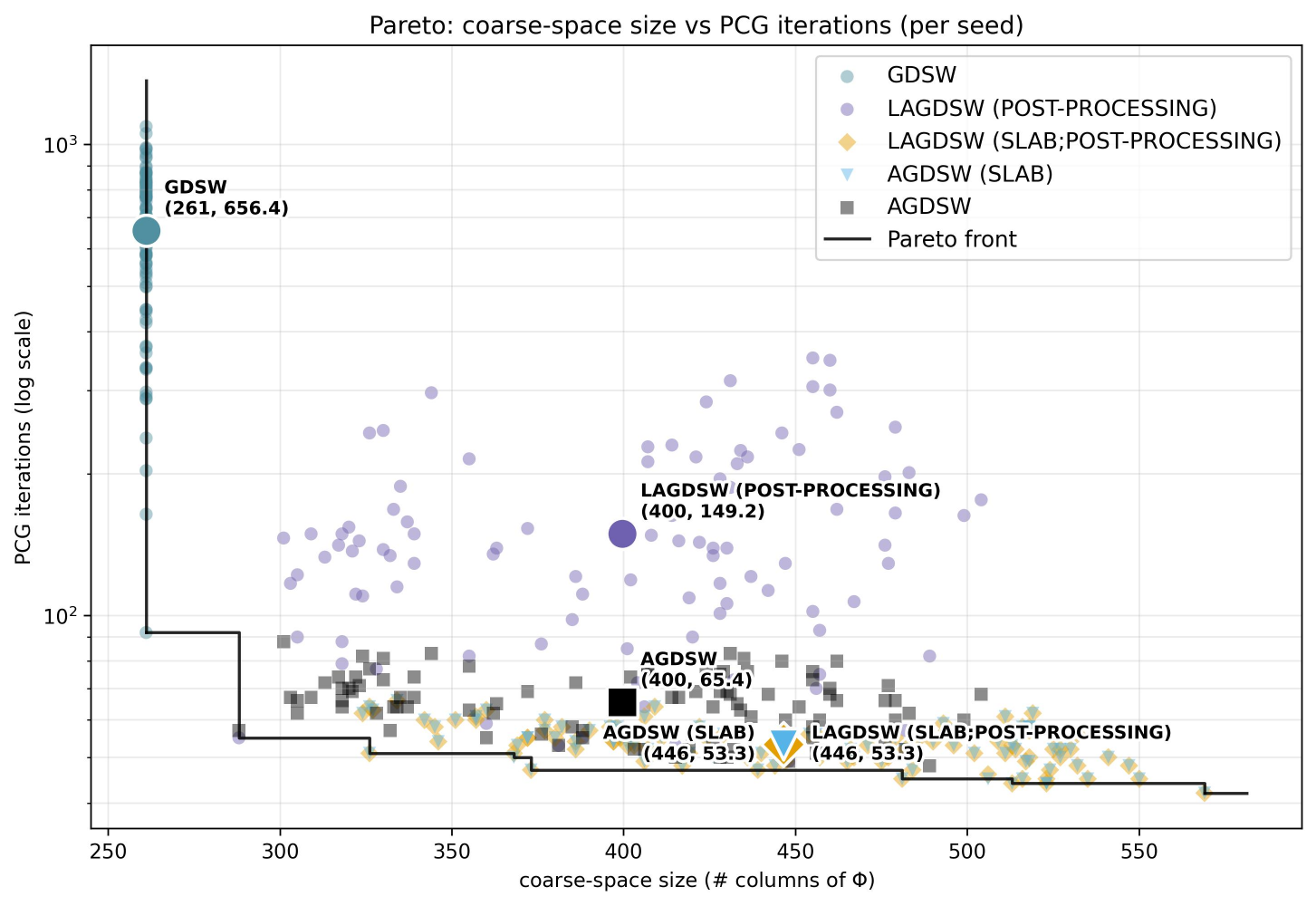}	
\caption{{\bf Pareto front plot:} Comparison of different setups (AGDSW-slab versus AGDSW; clean NN outputs versus post-processed ones) for 100 subdomains and  100 different coefficient distributions (randomly generated; see Fig.~\ref{fig:train_data} for similar examples used as training data). While the $x$-axis represents the size of the coarse space, the $y$-axis shows the number of PCG iterations needed for convergence; coarse spaces: GDSW, AGDSW, AGDSW-slab, and LAGDSW. Here, LAGDSW means that the second adaptive constraint in AGDSW/AGDSW-slab is replaced by a prediction while all other constraints are obtained from the eigenvalue problem. {\bf Top:} without post-processing; {\bf Bottom:} with post-processing.}
\label{fig:slab}
\end{figure}

Having trained a neural network predicting the second adaptive constraint of either AGDSW and AGDSW-slab, we can combine both with the post-processing described above. In this subsection, we will only replace the second constraint of either AGDSW or AGDSW-slab, respectively, by one of these four options to see which one is the best and with which one we will proceed in the following sections. As one can see from Fig.~\ref{fig:slab}, using the network trained on data from the AGDSW-slab coarse space and combine it with the post-processing is clearly the best choice. Therefore, we will proceed using this approach to also predict the third and fourth constraint and apply the learned coarse space to more complex problems. 

\begin{figure}
\centering
\includegraphics[width=0.9\textwidth]{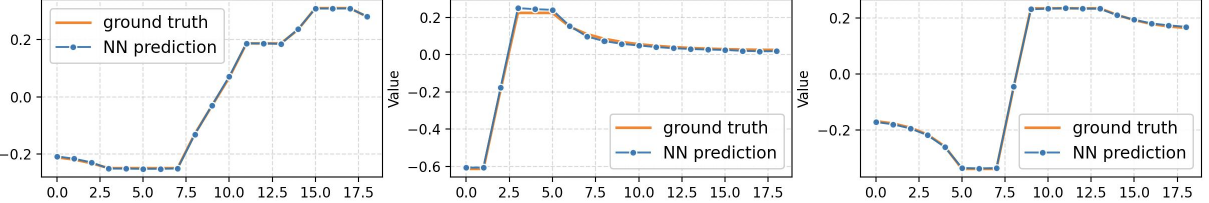}	
\includegraphics[width=0.9\textwidth]{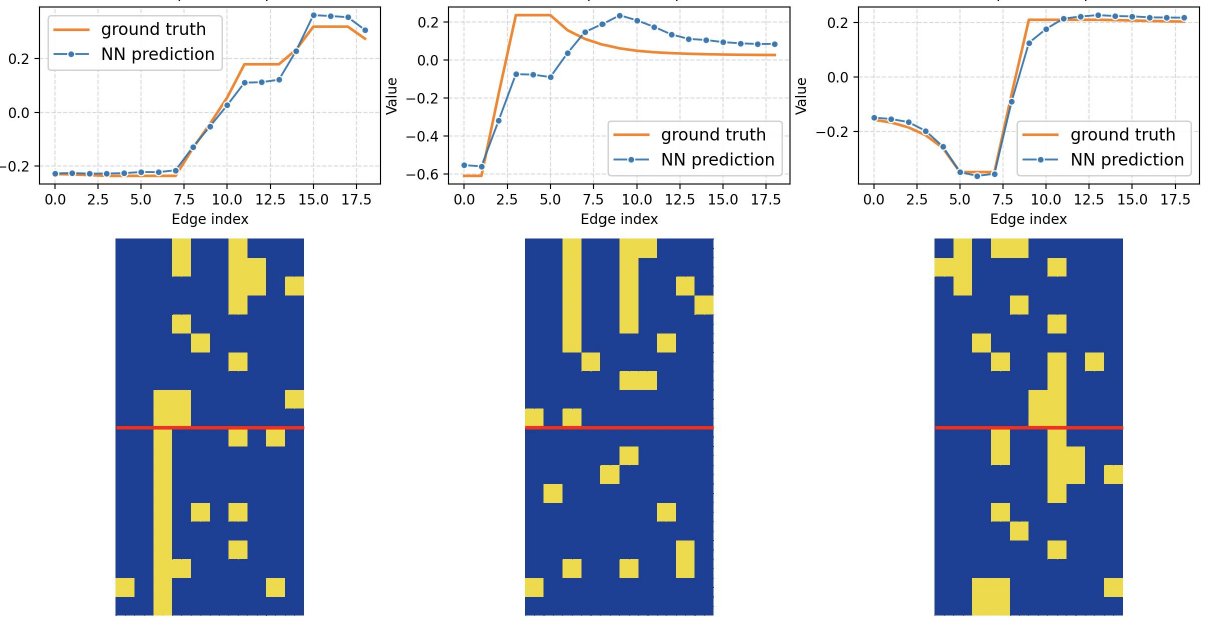}	
\caption{Three exemplary edges; from top to bottom: prediction with slab-based neural network getting area around the edge as sampling-input, prediction with neural network seeing complete input, coefficient distributions.}
\label{fig:plateau42}
\end{figure}

\begin{figure}
\centering
\includegraphics[width=0.8\textwidth]{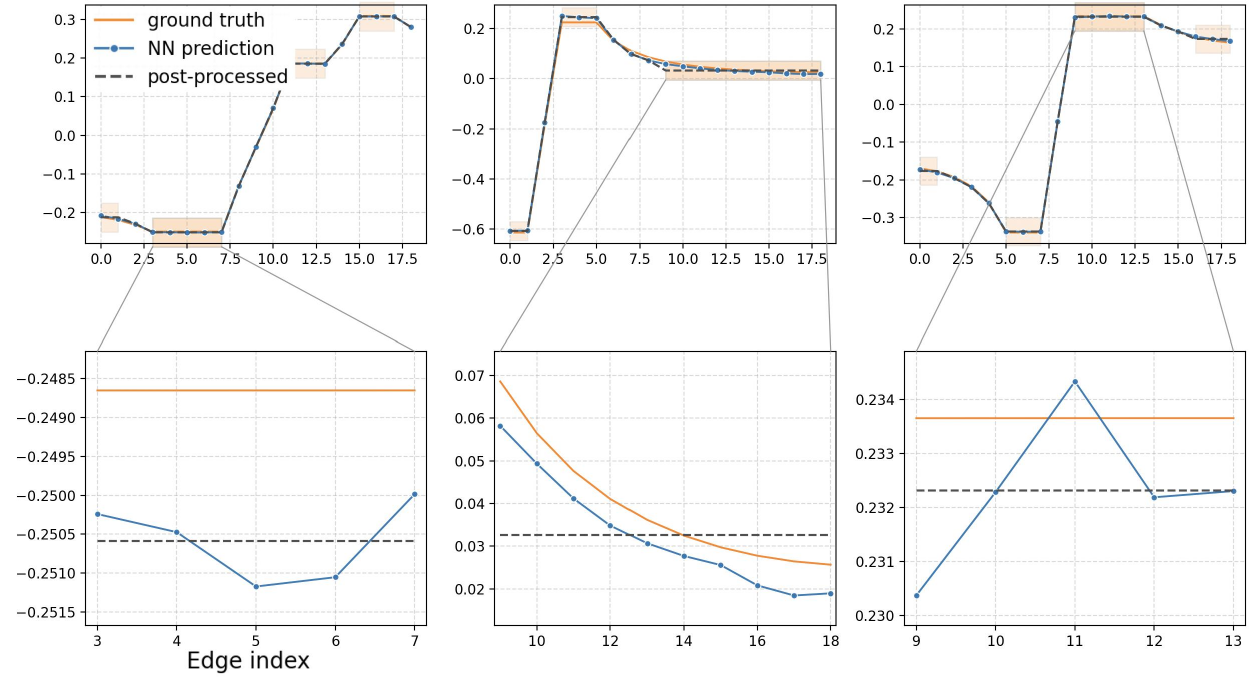}	
\caption{Same three edges as in Fig.~\ref{fig:plateau42} and same prediction with network acting on the slab; zoom to the largest  plateau found and modified during post-processing. {\bf Left and Right:} The small deviations from the ground truth can be seen; the plateau is identified successfully. {\bf Middle:} Due to the chosen tolerance, the largest plateau found does not actually correspond to any plateau in the ground truth; this can happen, but does not affect the performance of the coarse space much.}
\label{fig:plateau43}
\end{figure}

We would like to mention a few more details that can be learned from Fig.~\ref{fig:slab}. First, the AGDSW-slab coarse space tends to be slightly larger than AGDSW, which is to be expected since the harmonic extensions in the eigenvalue problem do not cover the entire subdomain. Second, without post-processing, the learned constraint corresponding to the AGDSW-slab coarse space is significantly more effective than that of AGDSW, which actually suggests a simpler learning task. If we now add post-processing, the learned constraint can achieve results nearly identical to those of AGDSW-slab. Post-processing also improves the predicted constraint corresponding to the classical AGDSW problem, though not to the same extent. The output of the neural network is not good enough to effectively find the plateaus during post-processing. To illustrate that, we provide the predictions on three edges for both networks including a zoom to a plateau area; see Fig.~\ref{fig:plateau42} and Fig.~\ref{fig:plateau43}.

\subsection{Robustness tests}
After choosing the best combination, that is, using the AGDSW-slab coarse space and sampling on the slab only plus post-processing, we use this setup in the remainder of the article. We extend it to the constraints three and four, combine it with a classifier, and test it for different problems and setups in this subsection.

\subsubsection{First tests with different coefficient distributions and different coarse spaces}

\begin{figure}[h!]
\centering
\includegraphics[width=0.9\textwidth]{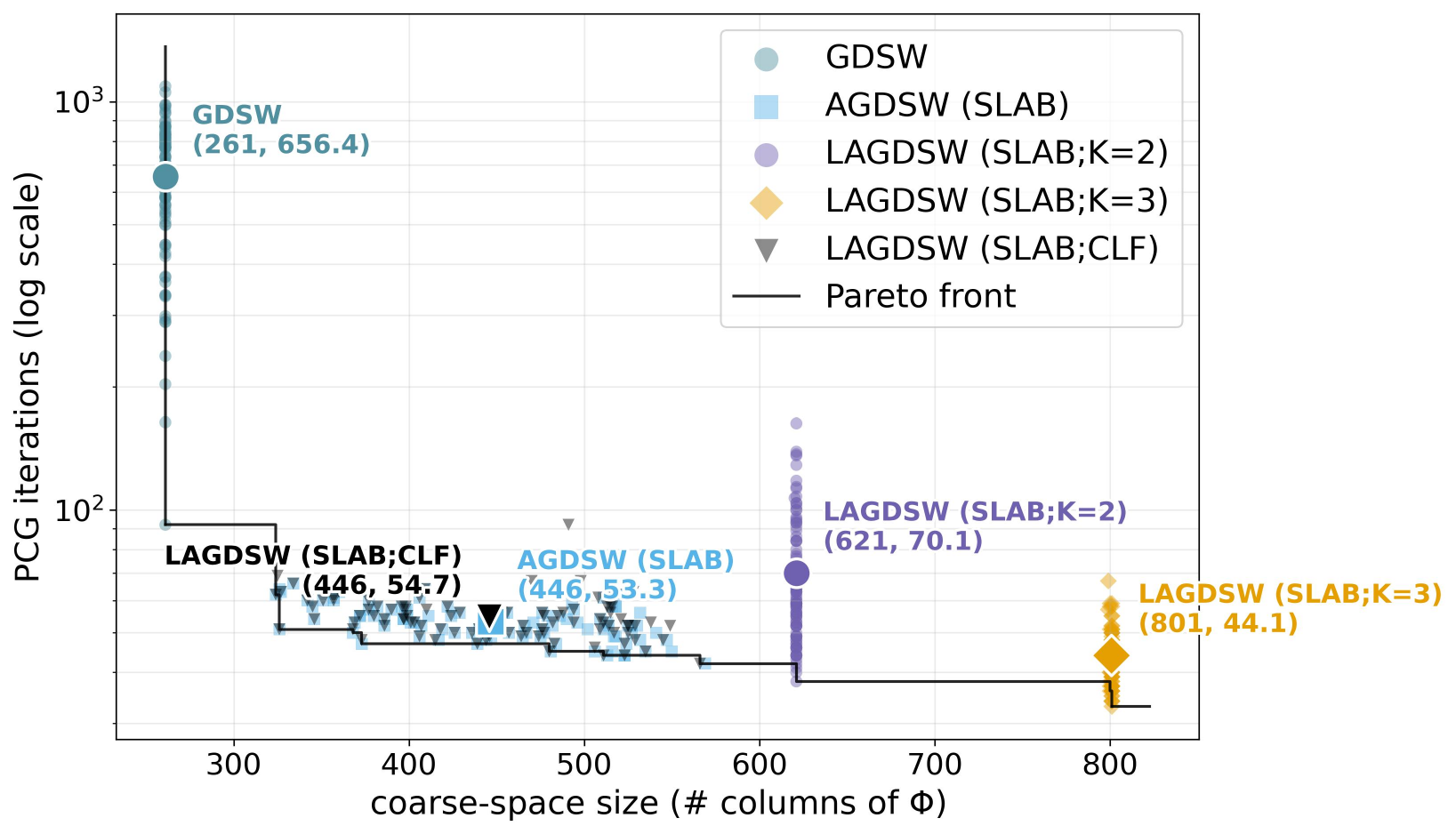}	
\includegraphics[width=0.9\textwidth]{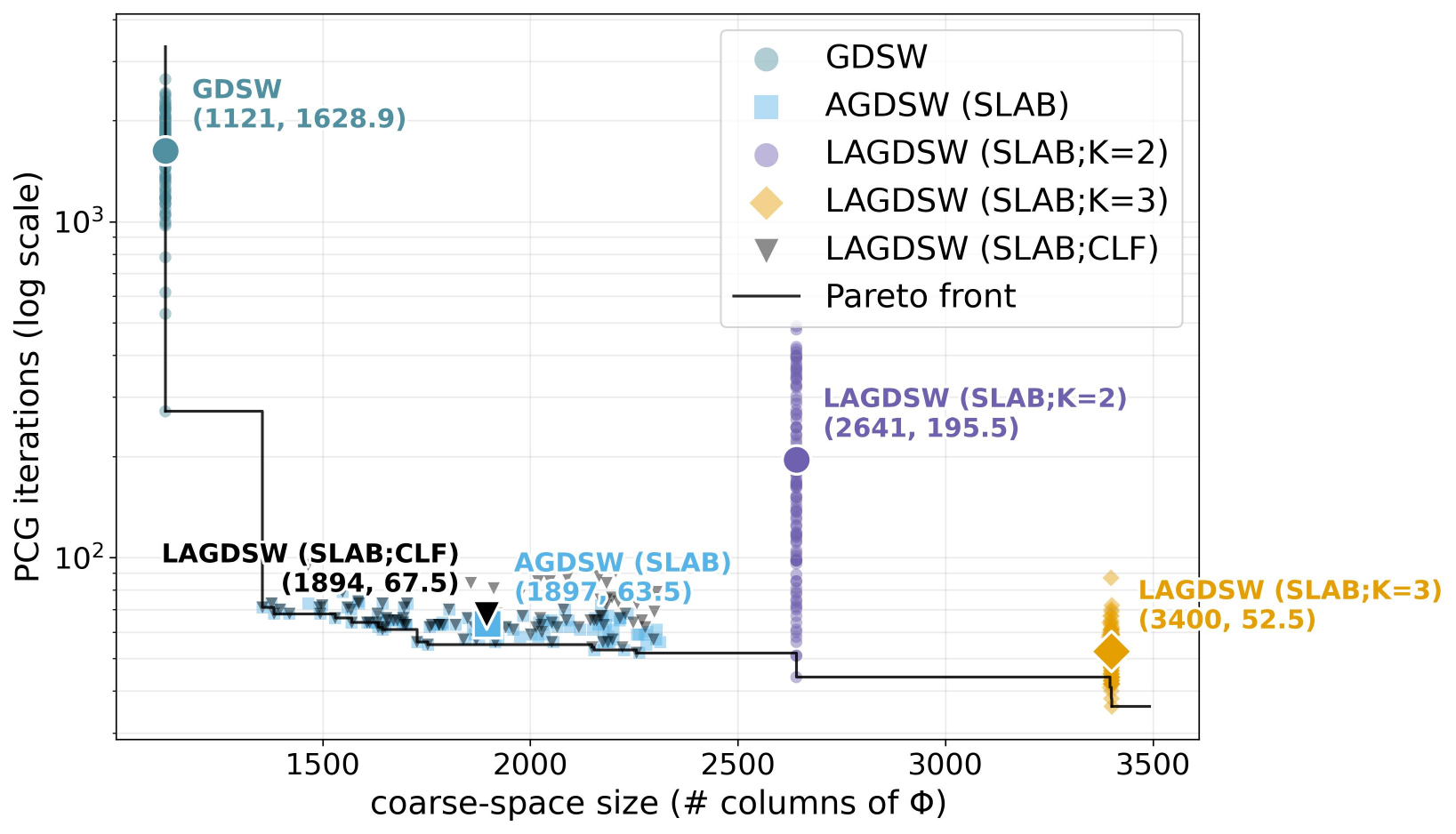}	
\caption{{\bf Pareto front plots:} coarse space size versus PCG iterations; comparison of different coarse spaces for 100 (top) and 400 (bottom) subdomains. While the $x$-axis represents the size of the coarse space, the $y$-axis shows the number of PCG iterations; $100$ different coefficient distributions (randomly generated; see Fig.~\ref{fig:train_data} for similar examples used as training data) are considered and solved with GDSW, AGDSW-slab, and LAGDSW-slab in three different variants using a fixed number of two/three constraints per edge and finally combining LAGDSW with a classifier. The means on the results over the 100 samples are visualized with the larger markers.}
\label{fig:pareto_diff}
\end{figure}

As we have seen, the first neural network can predict the second adaptive constraint, at least with additional post-processing.  So far, we just analyzed the quality of the constraints by replacing the second of the original adaptive constraints but did not suggest a usable solver setup. Such a solver without the need for any expensive eigenvalue problems can easily be built by adding the constant GDSW-type constraint plus a fixed number of $k=1$, $k=2$, or $k=3$ additional constraints predicted by the trained neural networks on each edge. Unfortunately, $k=3$ delivers a large coarse space while $k<3$ might result in a non-robust coarse space for some cases.  To overcome the difficulty to choose the right number of constraints, as described in Section~\ref{sec:clf}, we trained a classification model, which takes the same sampling input as the networks used to predict the constraints and classifies the considered edge into the classes 0 (GDSW-type constant constraint), 1 to 3 (one to three additional adaptive constraints, and 4 (more than 3 additional constraints); see again Fig~\ref{fig:ml-adaptive-blackbox} for an illustration of the full setup. Instead of using a fixed number of $k$ predicted constraints on each edge, we thus can also evaluate the classifier on each edge and only add the necessary number of predicted constraints. This potentially results in smaller coarse spaces which should still be similarly robust as in the $k=3$ case. In Fig.~\ref{fig:pareto_diff}, we compare GDSW and AGDSW-slab with three approaches with learned adaptive coarse spaces using $k=2$, $k=3$, or an adaptive number of predicted constraints per edge. In the latter case, the classifier is used as described. In the following, we always refer to learned AGDSW coarse spaces with LAGDSW. While GDSW obviously has the smallest coarse space, it is not robust for the complicated heterogeneous problems. Using two learned constraints per edge seems also not to be enough. In contrast, with $k=3$ learned constraints per edge, the number of PCG iterations is, on average, even lower than the one of the AGDSW-slab variant, but at a high cost of a large coarse space. The combination with the classifier performs best, mimicking the behavior of AGDSW-slab in terms of coarse space size and iteration count.    

\subsubsection{Scaling number of subdomains}

\begin{figure}
\centering
\includegraphics[width=0.75\textwidth]{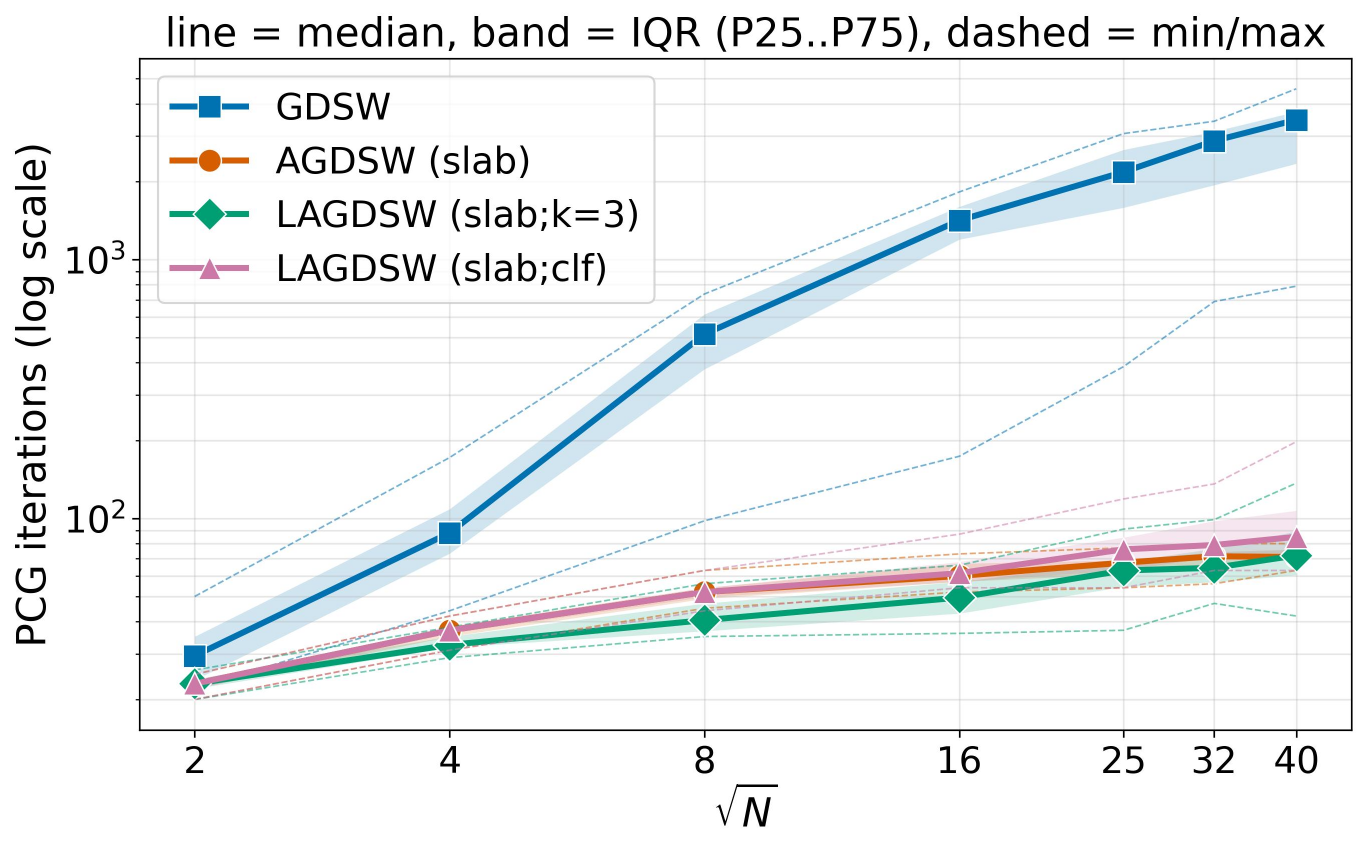}	
\includegraphics[width=0.75\textwidth]{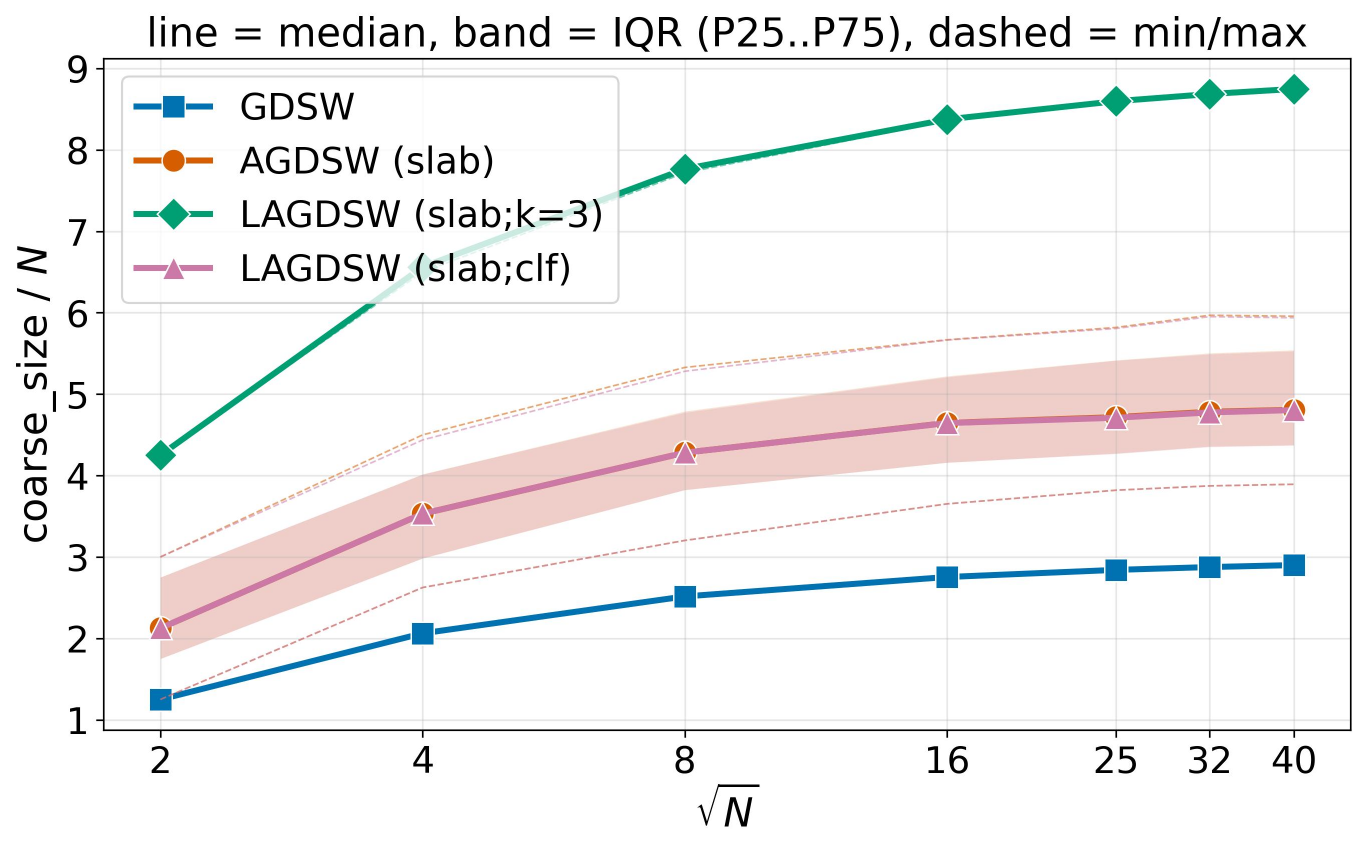}	
\caption{Comparison of scaling behavior of different coarse spaces; scalability from $N=4$ to $N=1\,600$ subdomains; $H/h=20$ and $\delta=2$. We tested $20$ different coefficient distributions (randomly generated; see Fig.~\ref{fig:train_data} for similar examples used as training data) and solved with GDSW, AGDSW-slab, and LAGDSW in two different variants using a fixed number of three constraints per edge and combining LAGDSW with a classifier. {\bf Top:} Iteration counts; {\bf Bottom:} average number of coarse constraints per subdomain, that is, number of coarse constraints divided by $N$.}
\label{fig:scaling_diff}
\end{figure}

\begin{figure}
\centering
\includegraphics[width=0.7\textwidth]{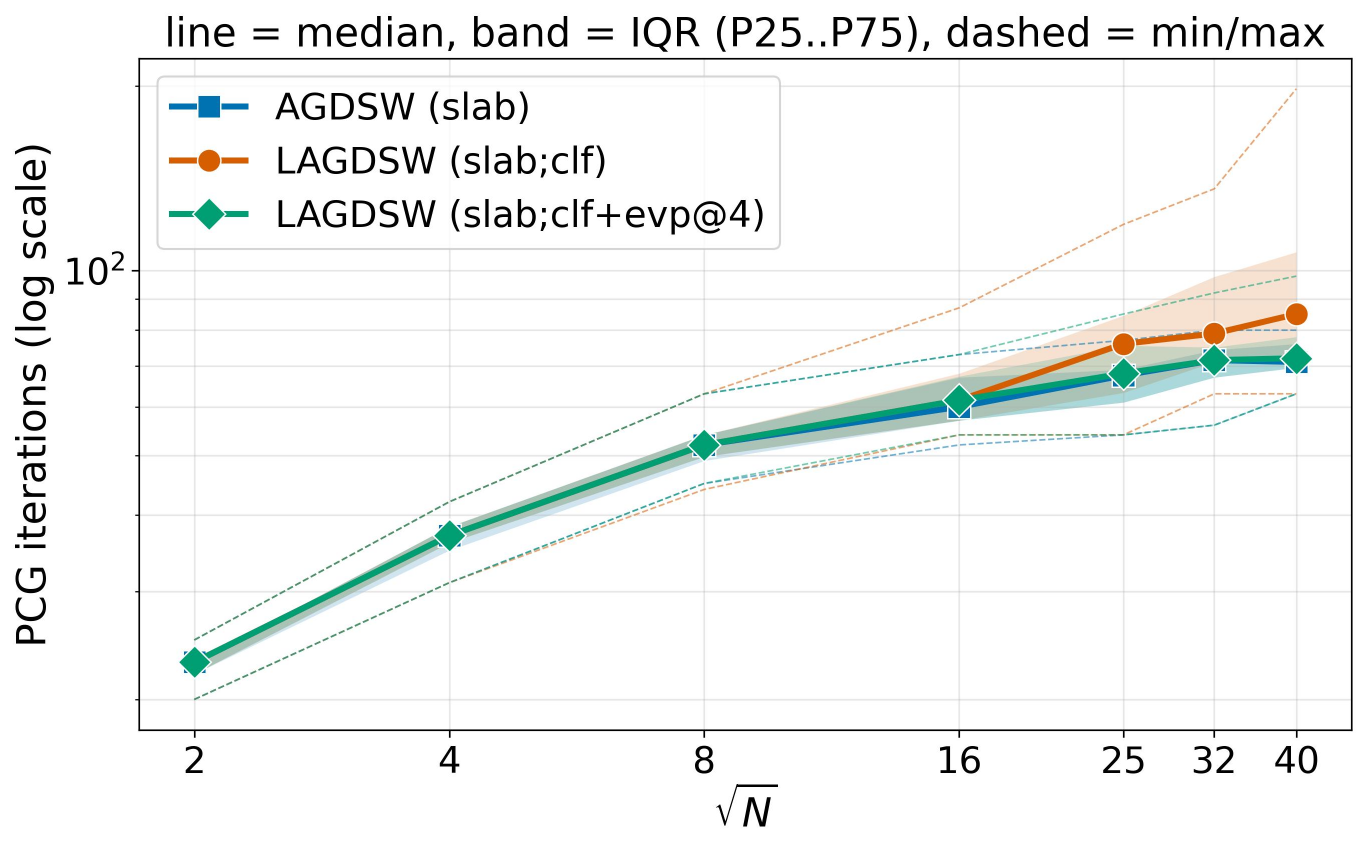}	
\caption{Comparison of scaling behavior of LAGDSW with classifier; scalability from $N=4$ to $N=1\,600$ subdomains; $H/h=20$ and $\delta=2$. The only difference of the two LAGDSW variants is the behavior on the few edges classified to class 4. While in the (slab;clf) case, class 4 edges are treated as class 3 edges, in the (slab;clf+evp@4) case the adaptive eigenvalue problem is solved on these edges and the exact adaptive constraints are used.}
\label{fig:scaling_diff_var}
\end{figure}

We provide also a small scalability study, where we fix the size of the subdomain to $H/h=20$ with an overlap of $\delta=2$ and increase the number of subdomains from $4$ to $1\,600$; see Fig.~\ref{fig:scaling_diff} for the results. While GDSW does not scale, the iteration counts for the adaptive approach stay below 80 iterations, also for the largest problems and all tested coefficient distributions. Also for the LAGDSW coarse space combined with a classifier, the median of $123.75$ for $1\,600$ iterations is still very good, especially compared with the one of GDSW, which is beyond $3\,000$ iterations. There are only a few outliers in the case of LAGDSW, and the worst one is 281 PCG iterations, which is still acceptable. Using always $k=3$ learned constraints for each edge performs slightly better, but the coarse space is, as already mentioned, much larger. All in all, the behavior of LAGDSW is quite robust and the method can be used also for large subdomain counts. Let us additionally remark that the reason for the outliers is not necessarily the quality of the neural networks or the post-processing. If the classifier puts an edge into class 4, that is, more than 4 adaptive constraints are needed, we nonetheless only use the first 4 constraints, that is, the constant one plus the prediction of our three neural networks. Therefore, we essentially use too few constraints for theoretical robustness. As an alternative, one could solve the adaptive eigenvalue problem on these few edges belonging to class $4$. If we do so, the maximum iteration count on $1\,600$ subdomains drops to 136 from 281; see also Figure~\ref{fig:scaling_diff_var}. That means, capping the coarse space on edges of class 4 has a bigger impact on the performance than the quality of the neural network. In practice, three possible strategies can be used: 1) solving the eigenvalue problem on edges of class 4, which is robust but might be expensive, 2) training even more models for the fourth, fifth constraint, 3) keeping the cap on the expense of slightly higher iteration counts. At the moment, we prefer the third  and most simple option.  

\subsubsection{Scaling subdomain size and coefficient jump}
\begin{figure}
\centering
\includegraphics[width=0.8\textwidth]{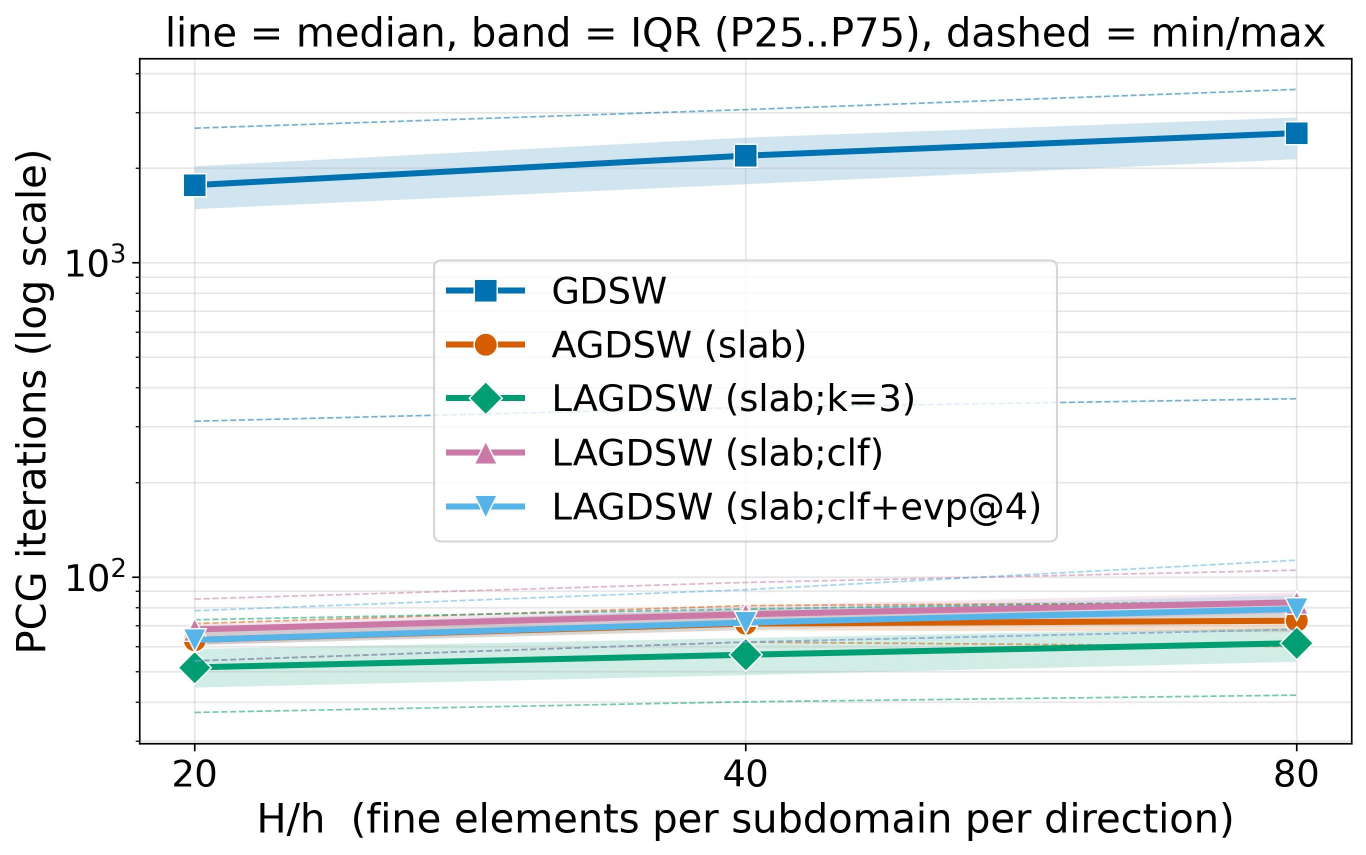}	
\includegraphics[width=0.8\textwidth]{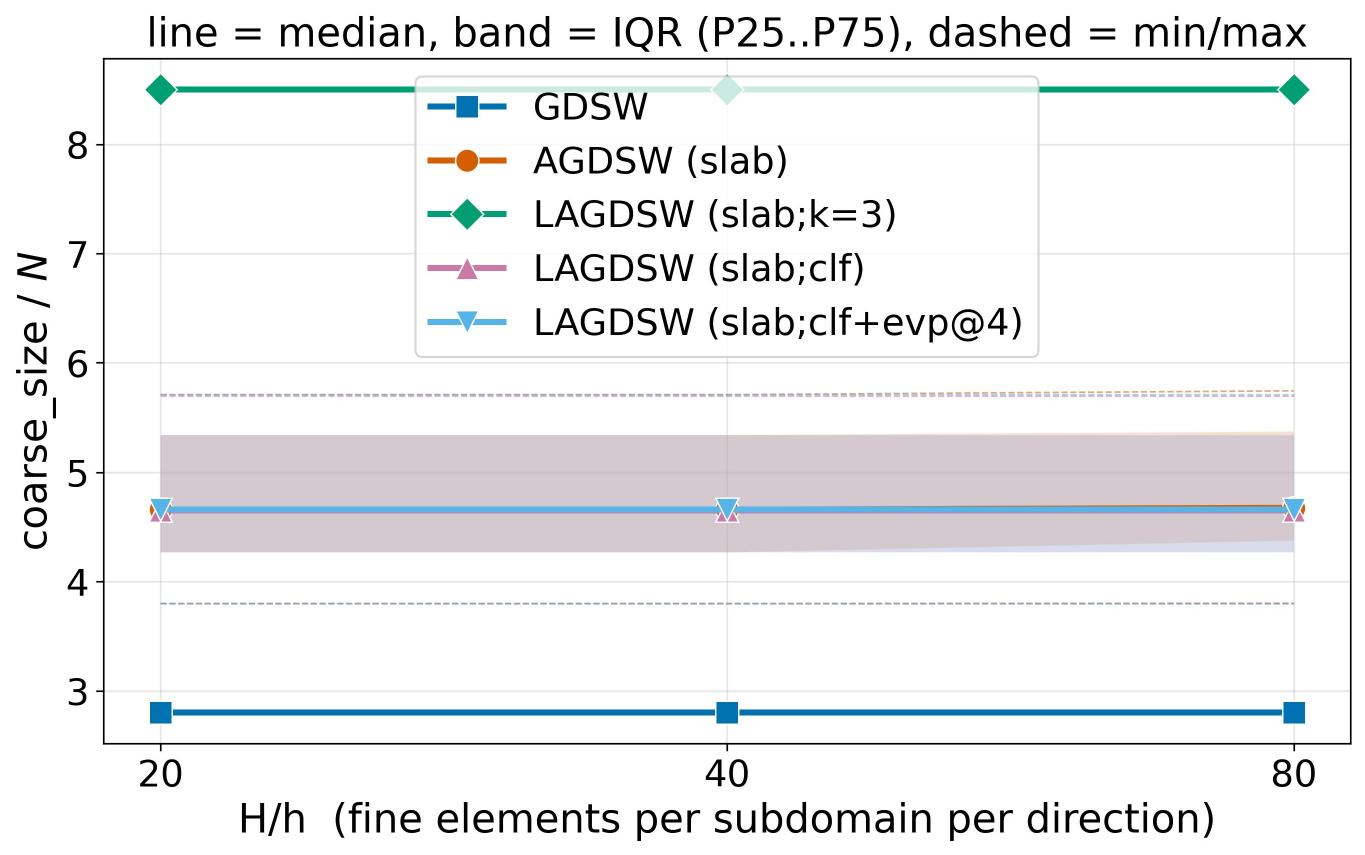}	
\caption{Increasing the subdomain size, that is, increasing $H/h$ from $20$ to $80$. We tested $20$ different coefficient distributions (randomly generated; see Fig.~\ref{fig:train_data} for similar examples used as training data) and solved with GDSW, AGDSW-slab, and LAGDSW-slab in different variants using a fixed number of three constraints per edge, combining LAGDSW with a classifier, and additionally solving the adaptive eigenvalue problem on edges assigned to class 4 (clf+evp@4). {\bf Top:} Iteration counts; {\bf Bottom:} average number of coarse constraints per subdomain, that is, number of coarse constraints divided by $N$.}
\label{fig:Hh}
\end{figure}

To obtain a complete picture, we also provide results for a test with increasing subdomain sizes starting from $H/h=20$ up to $H/h=80$; see~Fig.~\ref{fig:Hh}. All learned coarse spaces behave well and are robust against increasing subdomain sizes. Since the output of the neural networks has the length of an edge with 19 inner nodes, an interpolation on a finer grid is necessary. We choose a very simple linear interpolation approach with two special features: 1) all new finite element nodes at the start or the end of the edge are set to the same values as the starting- and, respectively, endpoint on the coarser mesh; 2) No linear interpolation is allowed to shorten an identified plateau. This can happen if one endpoint of a plateau (on the coarser mesh) lies inside a finite element of the finer mesh. In that case, the plateau is extended instead by that finite element of the finer grid. In case of our study this case does not occur, since we interpolate from 19 to 39 and finally to 79 nodes, which always matches perfectly. 

\begin{figure}
\centering
\includegraphics[width=0.75\textwidth]{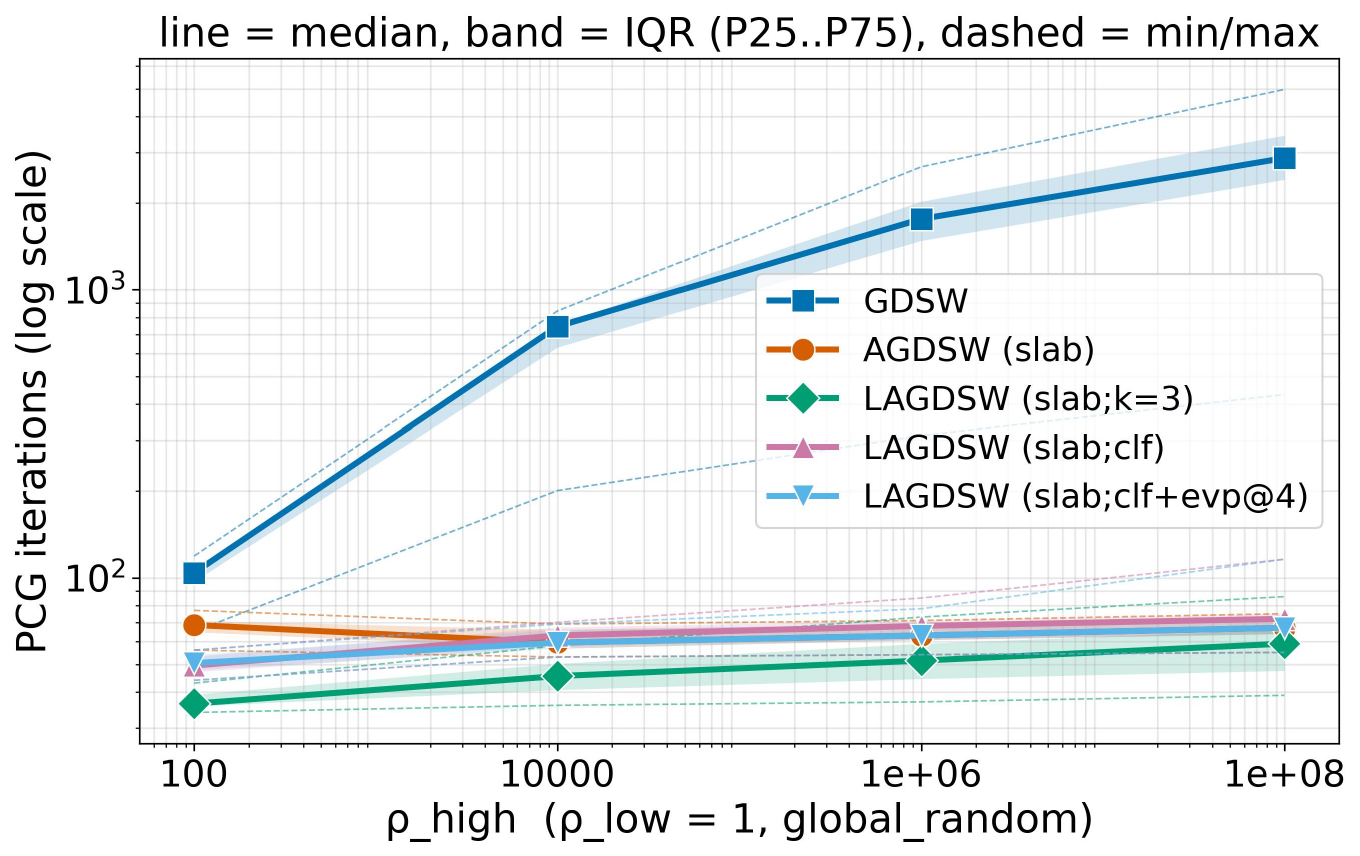}	
\includegraphics[width=0.75\textwidth]{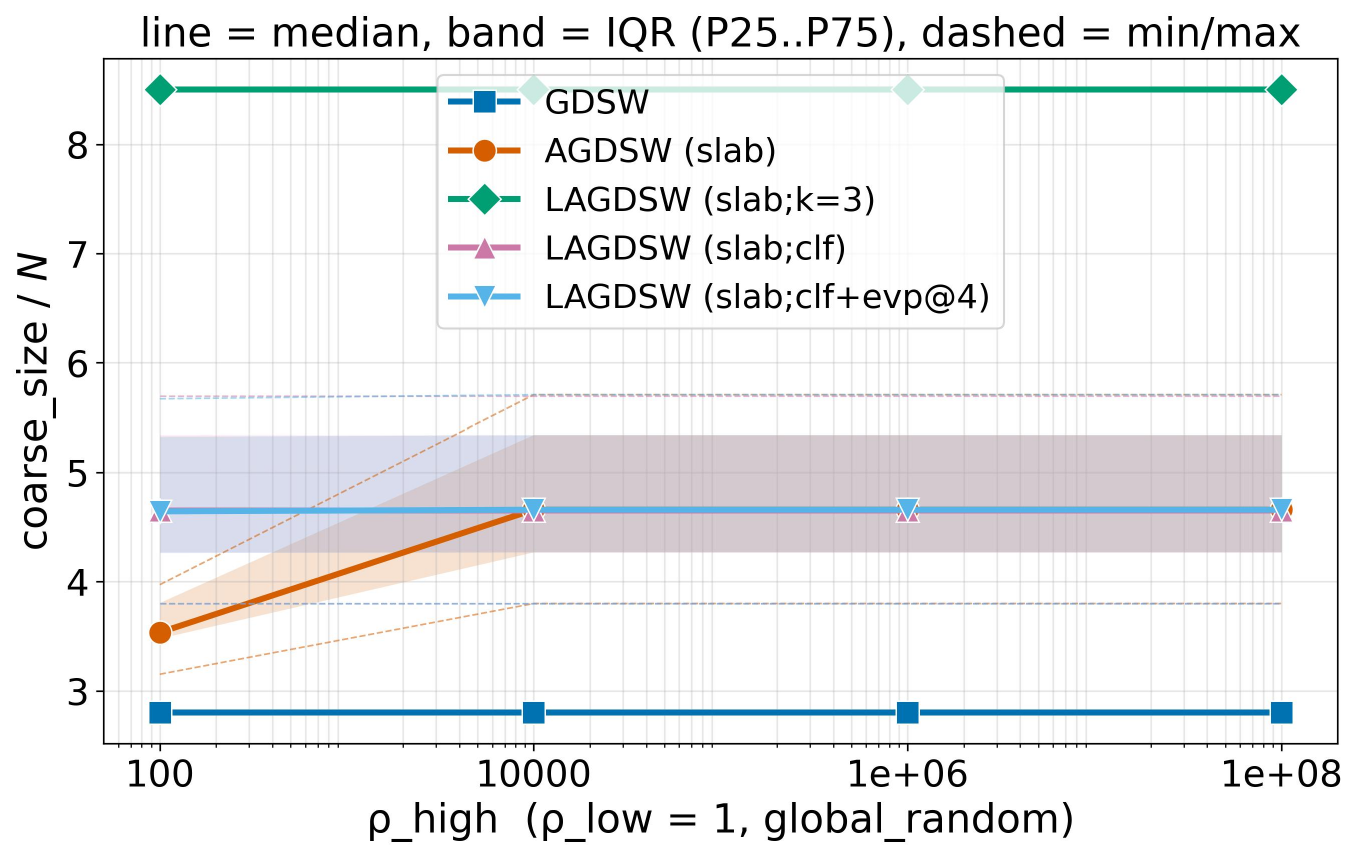}	
\caption{Increasing the high coefficient $\rho_{\rm max}$ from $100$ to $1e8$. We considered $20$ different coefficient distributions and solved with GDSW, AGDSW-slab, and LAGDSW-slab in different variants using a fixed number of three constraints per edge, combining LAGDSW with a classifier, and additionally solving the adaptive eigenvalue problem on edges assigned to class 4 (clf+evp@4). {\bf Top:} Iteration counts; {\bf Bottom:} Average number of coarse constraints per subdomain, that is, number of coarse constraints divided by $N$.}
\label{fig:rho_scaling}
\end{figure}

Additionally, we assume that the learned coarse spaces can also be used for different coefficient jumps, although only trained on data with a jump of $1e6$. We present results for increasing coefficient jumps in Fig.~\ref{fig:rho_scaling}. Only for small jumps, the size of the learned coarse space with classifier is larger than the AGDSW-slab coarse space, but still very effective.  

\section{Numerical results for linear elasticity}
\label{sec:num_elas}
 In this section, we provide results for linear elasticity problems, without retraining any machine learning model. Let us first introduce the linear elasticity test problem in detail.
\subsection{Linear elasticity}

Let $\Omega=(0,1)^2$, $\mathbf{u}:\Omega\to\mathbb{R}^2$ the displacement,
$E:\Omega\to\mathbb{R}_+$ Young's modulus, which we consider to be piecewise constant on each finite element,
$\nu=0.3$, and consequently the plane-strain Lam\'e parameters given by 
\[
\lambda=\frac{\nu E}{(1+\nu)(1-2\nu)},\qquad
\mu=\frac{E}{2(1+\nu)}.
\]
With
\[
\varepsilon(\mathbf{u})=\frac{1}{2} \left( \nabla u + \nabla u^T\right) {\rm and\; the\; stress}
\qquad
\boldsymbol{\sigma}(\mathbf{u})=\lambda\,\mathrm{tr}(\varepsilon(\mathbf{u}))\mathbf{I}
+2\mu\,\varepsilon(\mathbf{u}),
\]
we obtain the {\it linear elasticity problem} by

\begin{equation}\label{eq:elas-strong}
\begin{aligned}
-\operatorname{div}\boldsymbol{\sigma}(\mathbf{u}) &= \mathbf{f}
&& \text{in }\Omega,\\
\mathbf{u} &= 0
&& \text{on }\partial\Omega.
\end{aligned}
\end{equation}
Here, we always use the body force $\mathbf{f}=(0,-100)^\top$.

 We define
\[
\left(H_0^1(\Omega)\right)^2:=\{\mathbf{v}\in H^1(\Omega)^2:\ \mathbf{v}|_{\partial\Omega}=\mathbf{0}\}.
\]
The weak form is: Find $\mathbf{u}$ with $\mathbf{u}|_{\partial\Omega}=0$ such that
\begin{equation}\label{eq:elas-weak}
a(\mathbf{u},\mathbf{v})=F(\mathbf{v})\qquad\forall\mathbf{v}\in \left(H_0^1(\Omega)\right)^2,
\end{equation}
\begin{equation}\label{eq:elas-forms}
a(\mathbf{u},\mathbf{v}):=\int_\Omega \boldsymbol{\sigma}(\mathbf{u}):\varepsilon(\mathbf{v})\,dx,
\qquad
F(\mathbf{v}):=\int_\Omega \mathbf{f}\cdot\mathbf{v}\,dx.
\end{equation}

We again use linear finite elements for discretization and the resulting linear system writes
\begin{equation}\label{eq:elas-disc}
A\mathbf{u}=\mathbf{b},
\end{equation}
where $\mathbf{u},\mathbf{b}$ have two dofs per node, that is, the displacements in $x$- and $y$-direction $(u_x,u_y)$.

\subsection{Coarse spaces for linear elasticity}
For the problem of linear elasticity, it is well-known that the three rigid body modes, that is, two translations and one rotation, form the null-space of the operator. Consequently, the GDSW coarse space is defined slightly different compared to the stationary diffusion case, where only the constant function is in the kernel. Let us note, that in each node of the finite element mesh, we now have two degrees of freedom, one for the $x$- and one for the $y$-direction. For each vertex of the domain decomposition, we now consider two coarse basis functions. The first one is set to $1$ in the $x$-dof of the vertex and $0$ in the $y$-dof as well as in all DOFs of the remaining interface. The second coarse basis function is set to $1$ only in the $y$-dof of the vertex. Both basis functions are extended to the interior of the subdomains by a discrete harmonic extension using the full stiffness matrix. For each edge, we define three coarse basis functions. The two translations $t_x$ and $t_y$ on the edge are defined as $1$ in each $x$- or $y$-dof of the edge, respectively, and zero otherwise. The rotation $r$ on the edge is defined as $(-y,x)$ in each node on the edge, where $(x,y)$ are the coordinates of the respective node. Again, all three functions are extended using the full stiffness matrix.

 With the aim to extend our approach of LAGDSW to linear elasticity, the question is how can we exploit the learned constraints, which are originally made for scalar problems? To explain our heuristic approach, we recapitulate how the GDSW edge-functions can be lifted from the scalar to the elasticity case. The GDSW coarse basis function on an edge is the constant function $c$ with value one. In principle, we obtain the edge functions for elasticity by a node-wise multiplication with $[1,0]$, $[0,1]$, and $[-x,y]$, where the values for the two DOFs in a node are obtained by a multiplication of the nodal value of the scalar function $c$. We obtain the entries $2i-1$ and $2i$ of the constraints belonging to the $i$-th node of the edge by $t_x(2i-1,2i) = [1,0] * c(i)$, $t_y(2i-1,2i) = [0,1] * c(i)$, and $r(2i-1,2i)=[-y,x]*c(i)$, where $[x,y]$ are the coordinates of the specific node. Now, using the same principles, we can  lift any other scalar constraint, for example, any learned constraint $l$, by $l_x(2i-1,2i) = [1,0] * l(i)$, $l_y(2i-1,2i) = [0,1] * l(i)$, and $l_r(2i-1,2i)=[-y,x]*l(i)$. As a consequence, we build the LAGDSW coarse spaces for elasticity as follows: 
\begin{itemize}
\item[1)] sampling is done for the coefficient $E$ and the input is min-max scaled as for scalar PDEs,
\item[2)] we decide for the number of constraints used per edge, either choosing a fixed number or using the same classifier as before,
\item[3)] the regression NNs for the scalar case are evaluated,
 \item[4)] the post-processing is applied as in the scalar case,
 \item[5)] all constraints are lifted as described using the rigid body modes, and
 \item[6)] a final orthogonalization is done with QR factorization - eventually removing linearly dependent constraints on the edge.  	
\end{itemize}

In contrast, the definition of the AGDSW coarse space is rather simple. It is completely analogous to the scalar case using the same eigenvalue problem. Of course, the underlying linear operator is now the stiffness matrix and the eigenvalue problem is two times larger and more expensive. Nonetheless, the coarse space is again provably robust in both variants, that is, AGDSW and AGDSW on a slab. 

{\it Remark.} Let us note that there is no theory that supports the assumption that a lifted adaptive coarse space is robust for elasticity problems. This can also be seen in the numerical results, where for extremely large coefficient jumps and complicated micro-structures the lifted learned constraints deteriorate, while AGDSW is still robust. In these cases the AGDSW coarse space tends to be larger, which shows that some important coarse modes cannot be found and lifted easily in the scalar case. To overcome this issue, it is possible to train new neural networks with training data taken from the elasticity case. We do not choose this option here, since our goal is to demonstrate how far we can get with the surrogate models exclusively  trained on stationary diffusion data and how different problems can be accelerated with it.

\subsection{Comparison of coarse spaces}

We would like to answer two questions in this section: 1) can a good coarse space for elasticity problems be obtained by solving the adaptive eigenvalue problem on a scalar stationary diffusion problem with the same coefficient distribution and then lifting the scalar constraints by a multiplication with the rigid body modes? 2) If so, can we also use the lifted learned constraints?

\begin{figure}[h!]
\centering
\begin{minipage}[c]{0.29\textwidth}
    \centering
    \includegraphics[width=\linewidth]{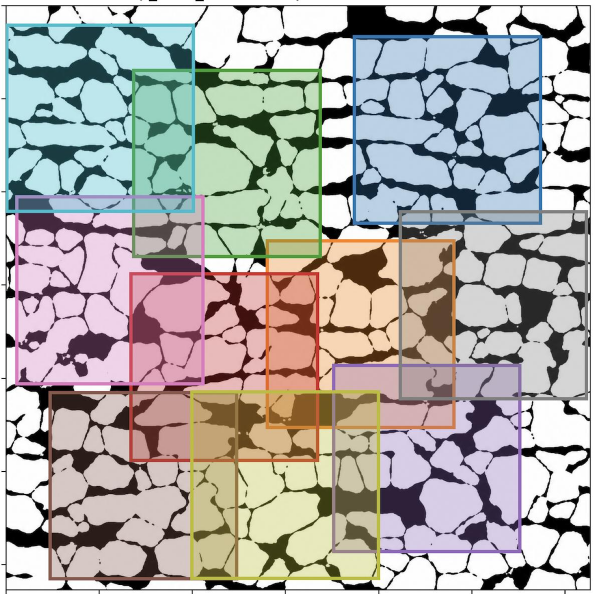}
\end{minipage}
\hfill
\begin{minipage}[c]{0.7\textwidth}
    \centering
    \includegraphics[width=\linewidth]{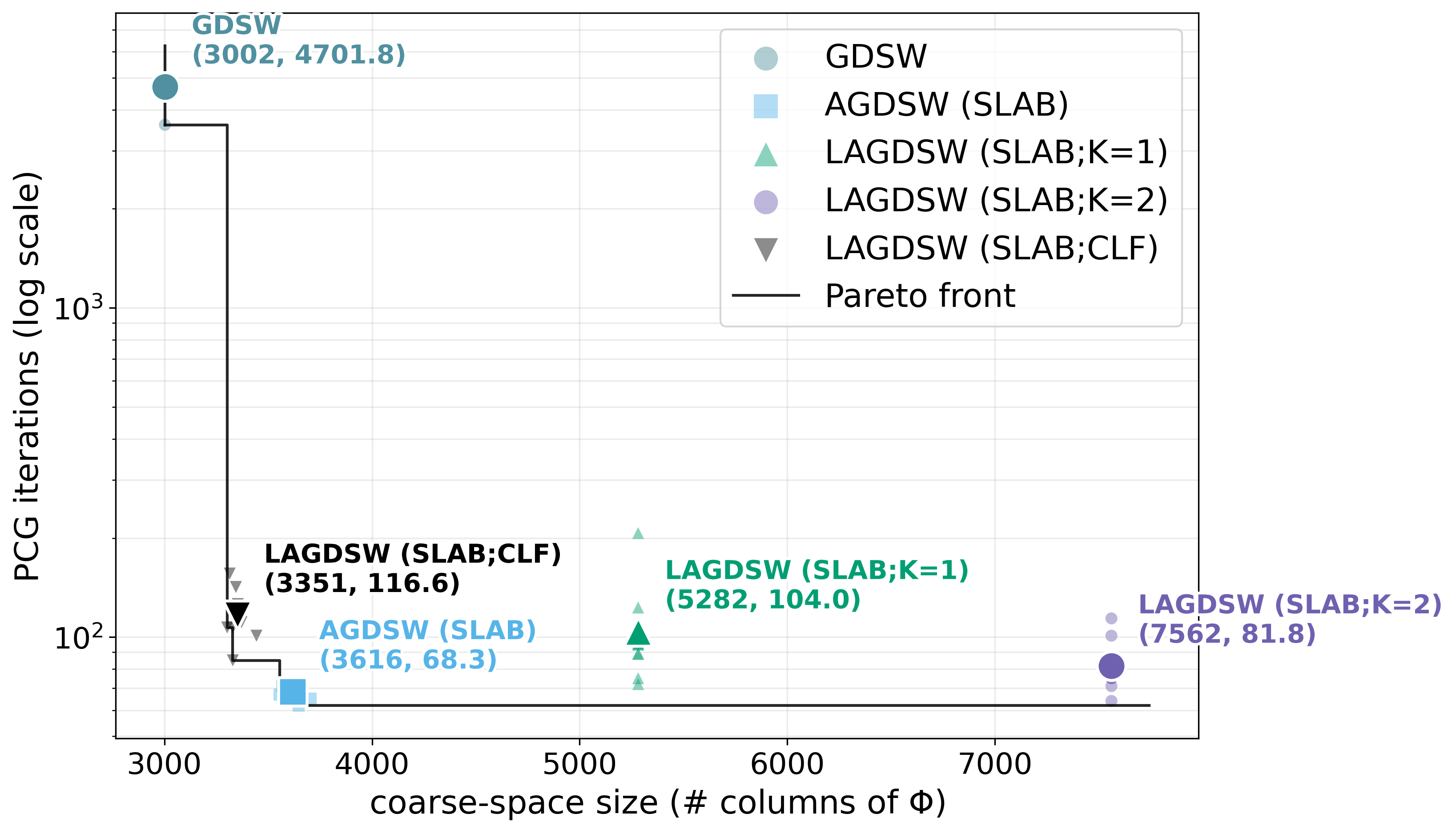}
\end{minipage}
\caption{{\bf Pareto front plot:} Coarse space size versus PCG iterations; comparison of different coarse spaces for linear elasticity and 400 subdomains; $E_{\rm low} =210$ and $E_{\rm high} =2.1e8$. While the $x$-axis represents the size of the coarse space, the $y$-axis shows the number of PCG iterations; $10$ different randomly chosen parts of a microsection of a realistic dual-phase steel are considered (see left image) and solved with GDSW, AGDSW-slab, and LAGDSW-slab in different variants. The means on the results over the 10 samples are visualized with the larger markers.}
\label{fig:el_dp}
\end{figure}

\begin{figure}[h!]
\centering
\includegraphics[width=0.9\textwidth]{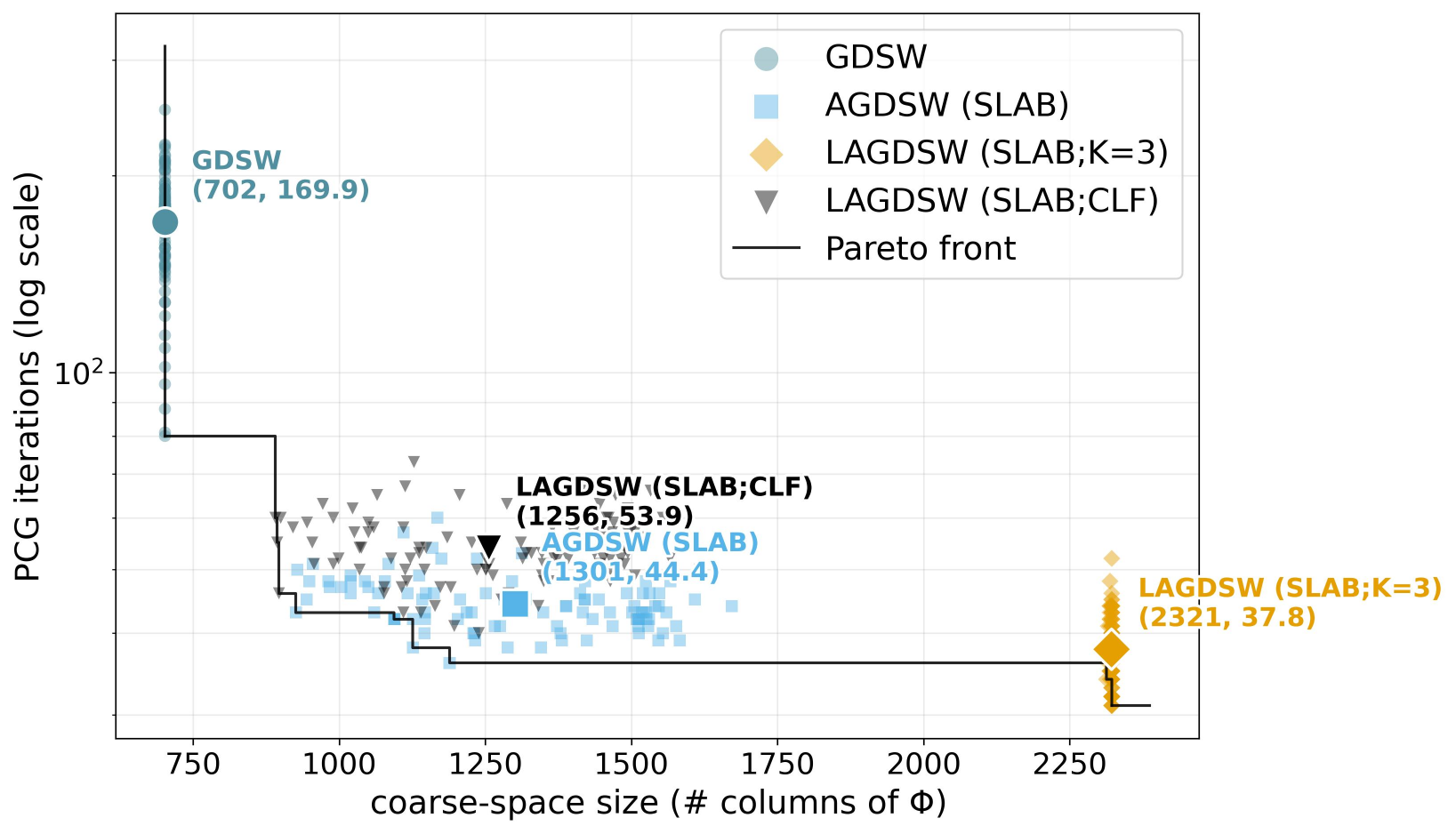}	
\includegraphics[width=0.9\textwidth]{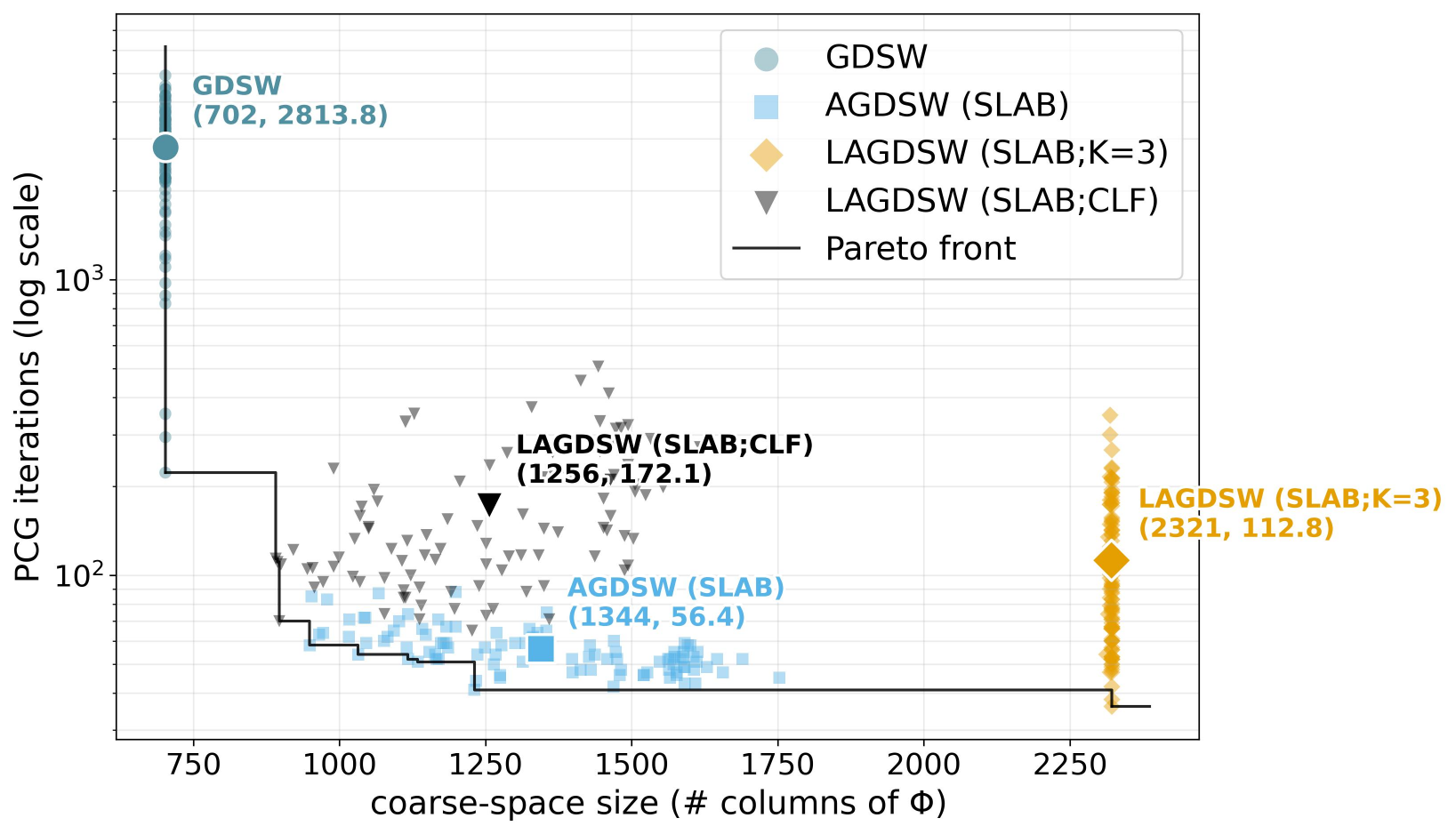}	
\caption{{\bf Pareto front plots:} Coarse space size versus PCG iterations; comparison of different coarse spaces for 100 subdomains for an elasticity problem. {\bf Top:} $E_{\rm low} =210$ and $E_{\rm high} =210\,000$; {\bf Bottom:} $E_{\rm low} =210$ and $E_{\rm high} =2.1e8$. While the $x$-axis represents the size of the coarse space, the $y$-axis shows the number of  PCG iterations; $100$ different coefficient distributions (randomly generated; see Fig.~\ref{fig:train_data} for similar examples used as training data) are considered and solved with GDSW, AGDSW-slab, and LAGDSW-slab in two different variants using a fixed number of three constraints per edge and combining LAGDSW with a classifier. The means on the results over the 100 samples are visualized with the larger markers.}
\label{fig:pareto_elas}
\end{figure}

\begin{figure}
\centering
\includegraphics[width=0.8\textwidth]{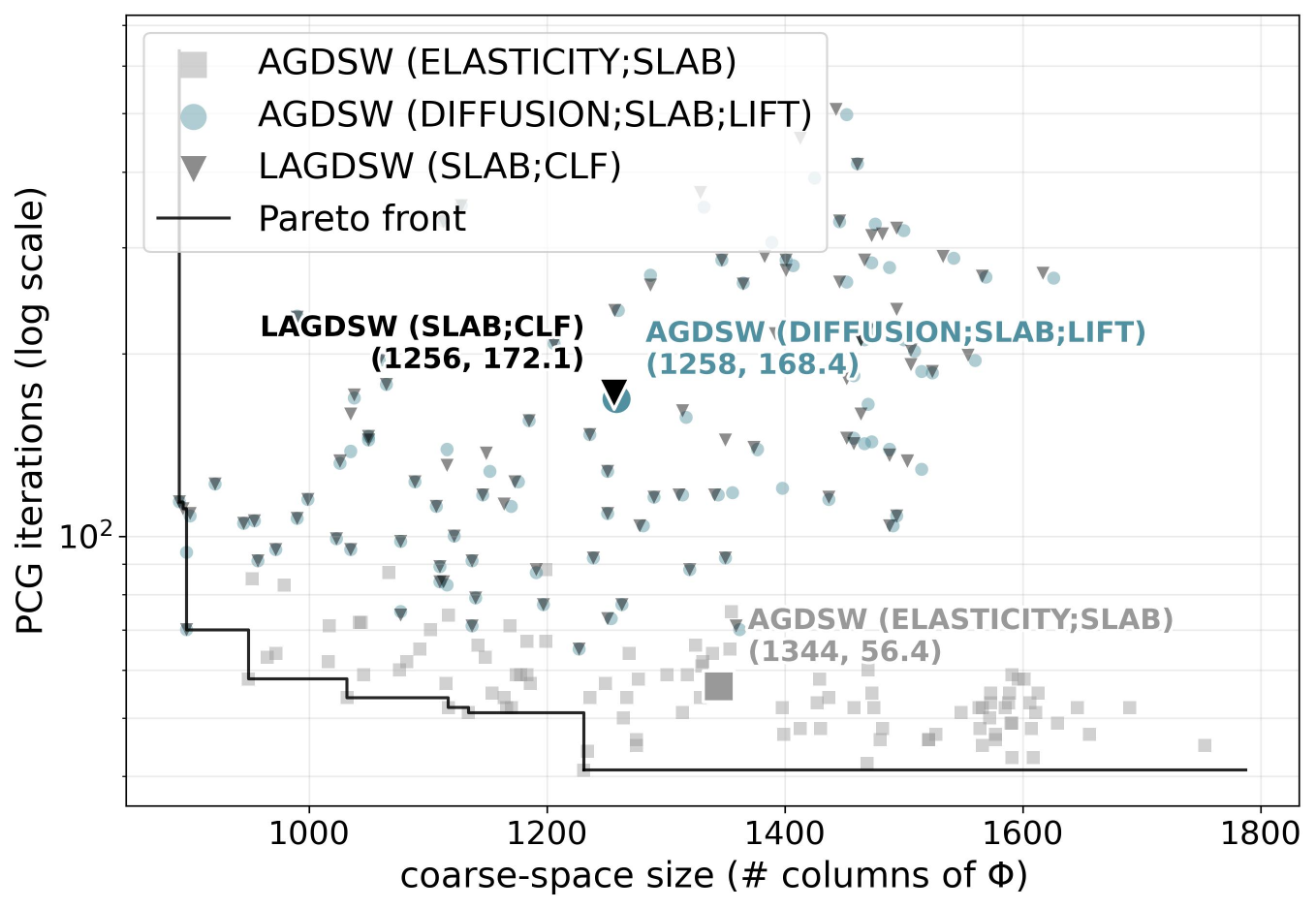}		
\caption{{\bf Pareto front plot:} Comparison of different coarse spaces for 100 subdomains for an elasticity problem; $E_{\rm low} =210$ and $E_{\rm high} =2.1e8$. While the $x$-axis represents the size of the coarse space, the $y$-axis shows the number of PCG iterations; $100$ different coefficient distributions (randomly generated; see Fig.~\ref{fig:train_data} for similar examples used as training data) are considered and solved with robust AGDSW-slab applied to the elasticity problem directly, heuristic AGDSW-slab applied to the stationary diffusion problem followed by a lift, and LAGDSW-slab with classifier and with lift afterwards. The means on the results over the 100 samples are visualized with the larger markers.}
\label{fig:pareto_lift}
\end{figure}

We present first results in Fig.~\ref{fig:pareto_elas}, where we considered the same 100 coefficient distributions, but different coefficient jumps of $1\,000$ and $1e6$. The first interesting observation is that the AGDSW-slab coarse space tends to be larger for the higher jump, that is, on average 1344 constraints instead of 1301. This is in contrast to the typical behavior in stationary diffusion problems and already a hint that the approach of lifting scalar constraints cannot catch all important adaptive modes. Nonetheless, for a physically meaningful jump of 1000, which is realistic, for example, for dual-phase steels, the lifted LAGDSW coarse space performs very well and at least outperforms GDSW by far. For the larger jump, it can be clearly seen that some modes are not caught and the lifted LAGDSW coarse space deteriorates, still outperforming GDSW. Let us note, that for a more realistic dual-phase steel micro-structure, for $400$ subdomains, LAGDSW performs very well also for the larger coefficient jump; see Fig. \ref{fig:el_dp}. To summarize, the heuristic approach of lifting scalar adaptive or learned coarse spaces works well for realistic coefficient distributions and realistic coefficient jumps, but cannot catch all important modes found by the theoretically robust eigenvalue problem for elasticity. This is also supported by the comparison presented in Fig.~\ref{fig:pareto_lift}, where we compare three modes: the theoretically robust AGDSW (slab) coarse space applied to the elasticity problem, the heuristic approach of building the AGDSW-slab coarse space for stationary diffusion and lifting it, and the lifted learned coarse space. It can be seen that the neural network predictions of the scalar constraints are pretty good, but both lifted scalar coarse spaces cannot span the same coarse space as the theoretically sound AGDSW-slab for elasticity.   

\clearpage 
\section{Numerical results for nonlinear problems}
\label{sec:num_elas_nonlin}
We also consider the solution of discretized nonlinear PDEs using two-level Schwarz methods, both, as a preconditioner for the tangent system in Newton's method and as a nonlinear preconditioner for Newton's method. For completeness, also one-level nonlinear Schwarz methods are compared. Again, we will exclusively use the learned coarse spaces trained with data from the linear stationary diffusion case. 

\subsection{Nonlinear scalar problems}

We consider a nonlinear PDE based on the $p$-Laplace operator. All nonlinear experiments in this work are posed on the unit square
$\Omega=(0,1)^2$ with a structured  triangulation, as already used in the linear studies.
We seek a scalar unknown $u:\overline\Omega\to\mathbb{R}$ and the heterogeneous permeability $\rho$
is piecewise constant on each finite element.

\paragraph{$p$-Laplace equation}\label{sec:plaplace}
Let $p>2$ and $\rho:\Omega\to(0,\infty)$ be piecewise constant.
The \textit{$p$-Laplace} problem is
\begin{equation}\label{eq:plaplace-strong}
  -\nabla\!\cdot\!\bigl(\rho\,|\nabla u|^{\,p-2}\,\nabla u\bigr)=f
  \quad\text{in }\Omega,
  \qquad u=0\quad\text{on }\partial\Omega.
\end{equation}
For $p=2$ one recovers the linear diffusion equation.
In our benchmarks we take $f\equiv 1$ and $p=4$.

Considering a weak formulation, discretization in the finite element space $V^h$, and collecting all terms on the lefthand side finally leads to a discrete nonlinear problem 
$$
F(u)=0,
$$ 
with $u \in V^h$; see~\cite{nl_right1} for details.

\subsection{Newton-Krylov-Schwarz and Nonlinear Schwarz}

We use two different classes of nonlinear solvers and investigate their performance when LAGDSW is used as coarse space. The first nonlinear solver is the classical {\it Newton-Krylov-Schwarz} approach, that is, for a discrete nonlinear problem
$$
F(u)=0,
$$
we solve iteratively $u^{(k+1)} = u^{(k)} + \lambda \delta u^{(k)}$, where $\delta u^{(k)}$ is obtained by solving the linearized system
$$
DF(u^{(k)}) \delta u^{(k)} = -F(u^{(k)}).
$$
Here, $\lambda$ is an optional step length and $DF(\cdot)$ the Jacobian of $F$. The linearized system is solved iteratively with a Krylov subspace method (PCG/GMRES) and we use a Schwarz-type preconditioner for acceleration. The preconditioner is built using the matrix $A := DF(u^{(k)})$. Let us note that we can either compute the adaptive coarse basis functions, that is, the columns of $\Phi$, using $DF(u^{(k)})$ in each Newton iteration or compute $\Phi$ once using $DF(u^{(0)})$ and recycle $\Phi$  in all following Newton steps. We always use the latter approach in all our computations. All other parts of the Schwarz preconditioner are always rebuilt in each Newton step using $DF(u^{(k)})$, especially the coarse matrix $\Phi^T DK(u^{(k)}) \Phi$ and its factorization. We use two-level Schwarz preconditioners with GDSW, AGDSW, and variants of LAGDSW coarse spaces.   All coarse spaces affect the number of Krylov subspace iterations but, of course, not the number of Newton iterations, unless some rounding errors cause small deviations. 

This is completely different for our second nonlinear solver. We consider the nonlinear additive one- or two-level Schwarz method; see~\cite{nl_left3,nl_adaptive} for details. Here, the domain decomposition method is used to form a nonlinear preconditioner $G$ which is applied to the residual $F$, resulting in the preconditioned nonlinear system $\mathcal{F}(u) = G(F(u))=0$. Afterwards, the nonlinearly preconditioned problem is solved using Newton's method. Depending on its exact choice, the application of the nonlinear Schwarz preconditioner results in the solution of local and coarse nonlinear problems in each outer Newton iteration. As a result, all parameters of the method, for example, the size of the overlap, the tolerances used in the nonlinear solves of the nonlinear preconditioner, the choice of the coarse space, and the coupling between levels, affect the nonlinear convergence and the linear one. The rule of thumb is that, if all parameters are chosen appropriately, nonlinear Schwarz methods converge rapidly. Here, we will investigate if LAGDSW trained for linear stationary diffusion problems can be used in the context of nonlinear Schwarz methods without deteriorating convergence. We will not go into a detailed discussion of the other parameters. We will also not give a detailed description of nonlinear Schwarz methods itself, since we want to focus on the learned coarse spaces and their effects. Let us just note that, for the first level of the two-level nonlinear Schwarz method, we use the RASPEN approach introduced in~\cite{nl_left2} with a multiplicity scaling in the overlap to obtain a partition of unity for the local prolongation operators. The same RASPEN method also constitutes the one-level nonlinear approach in all comparisons, that is, using no coarse space. The nonlinear coarse space is coupled additively and spanned by $\Phi$ as in the Newton-Krylov-Schwarz case; see~\cite{nl_adaptive} for details. As stated above, we can recycle coarse space basis functions from the first step and will do so in all experiments. In each Newton iteration, a linear system with left hand side $D \mathcal{F}(\cdot) = DG(F(\cdot))\,DF(\cdot)$ has to be solved. This matrix has the same structure as $M^{-1} DF$, where $M^{-1}$ is the corresponding linear Schwarz preconditioner - either two- or one-level. There is only one difference to the Newton-Krylov-Schwarz approach: the local and coarse matrices of the preconditioner are linearized in different points, which is a result of the application of the nonlinear preconditioner. This might affect linear iteration counts slightly, but, in our experience, not significantly. Let us note that GMRES is always applied directly to the system $D\mathcal{F}(\cdot)$ and no additional linear preconditioner is necessary. 

\begin{figure}[h!]
\centering
\includegraphics[width=1.0\textwidth]{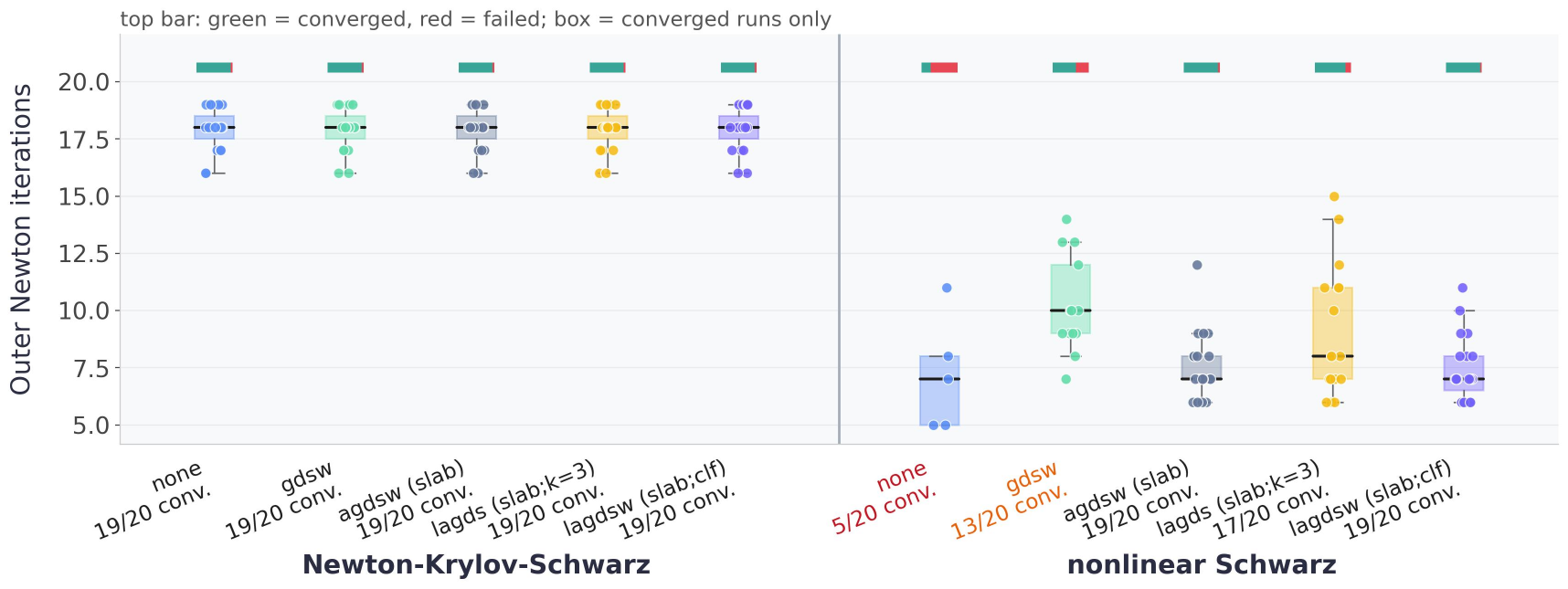}
\includegraphics[width=1.0\textwidth]{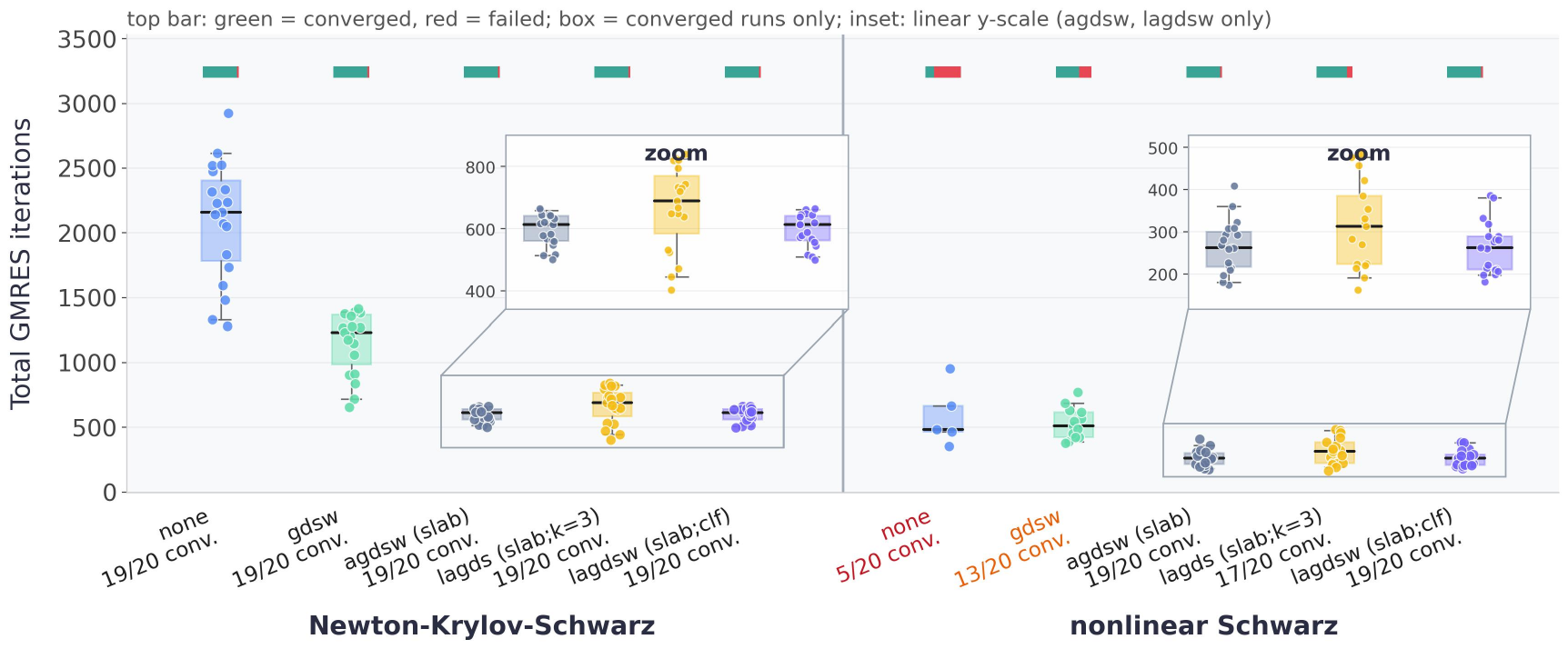}
\caption{Comparison of Newton-Krylov-Schwarz and nonlinear Schwarz methods with different coarse spaces for 20 randomly generated coefficient distributions; Newton iterations {\bf (top)} and total number GMRES iterations {\bf (bottom)}; high coefficient of $\rho=1e4$ and overlap of $\delta=2$.}
\label{fig:42}	
\end{figure}

\begin{figure}
\centering
\includegraphics[width=1.0\textwidth]{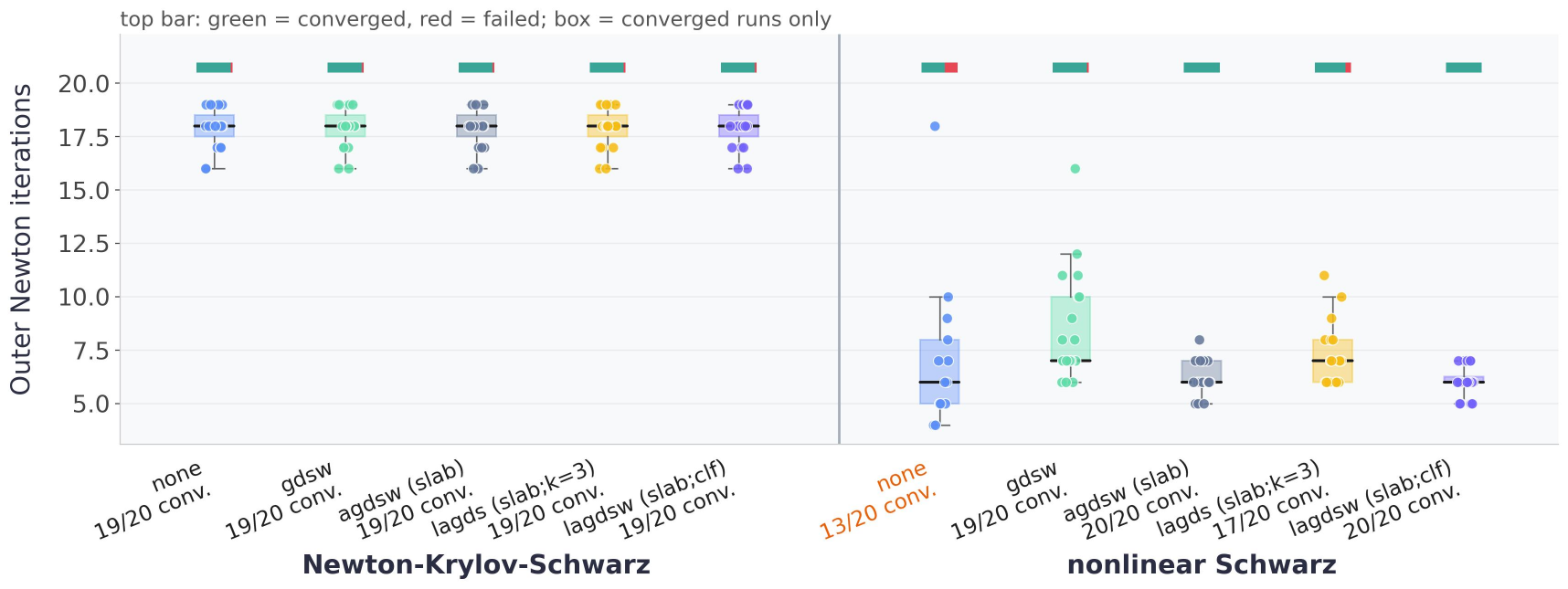}
\includegraphics[width=1.0\textwidth]{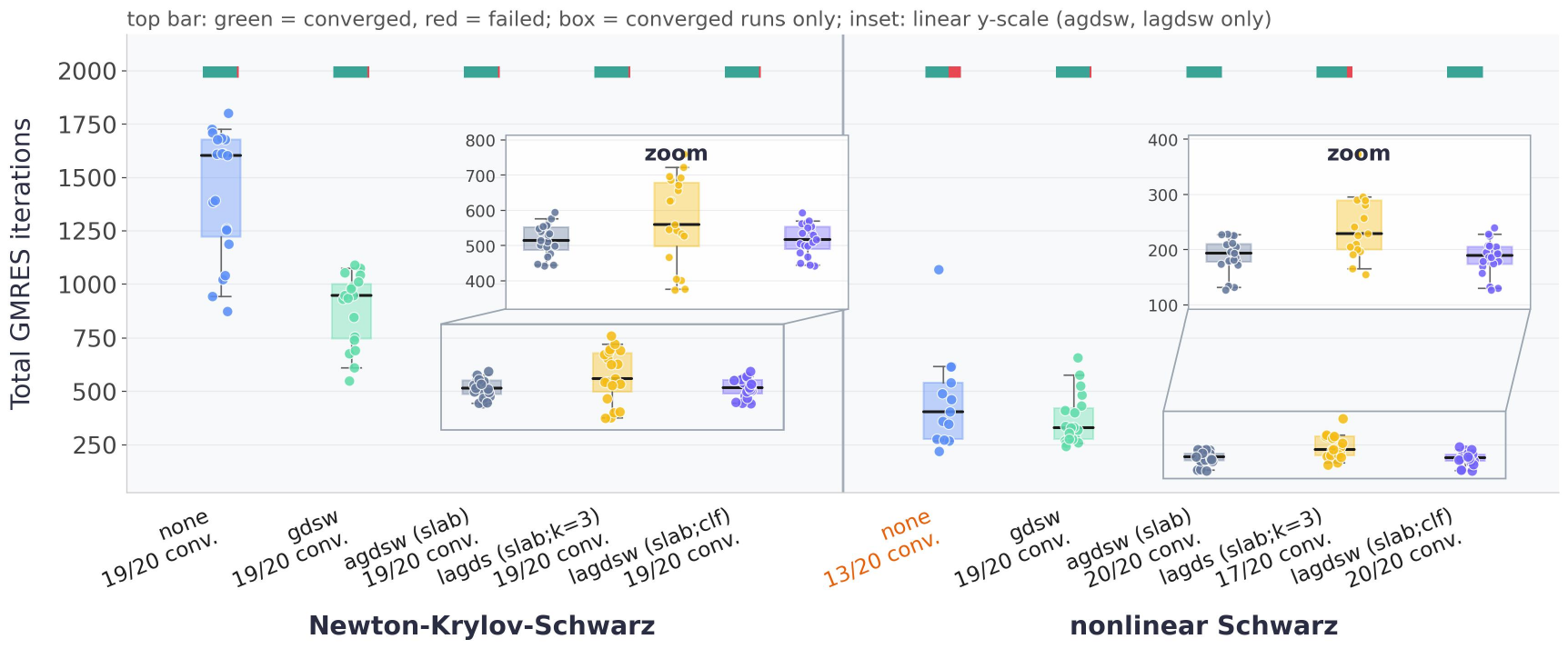}
\caption{Comparison of Newton-Krylov-Schwarz and nonlinear Schwarz methods with different coarse spaces for 20 randomly generated coefficient distributions; Newton iterations {\bf (top)} and total number GMRES iterations {\bf (bottom)}; high coefficient of $\rho=1e4$ and overlap of $\delta=4$.}
\label{fig:44}	
\end{figure}

\subsection{Results for the $p$-Laplace}
In general, nonlinear Schwarz methods only pay off if they can reduce the number of outer Newton iterations compared to classical Newton's method, since the application of the nonlinear preconditioner is easily parallelizable but certainly not for free. Therefore, to present the full picture, we use 20 different random coefficient distributions with a high coefficient of $\rho=10\,000$ or $\rho=1\,000\,000$ and a low coefficient of $\rho=1$, and compare nonlinear and linear iteration counts for Newton-Krylov-Schwarz and nonlinear Schwarz methods. We consider GDSW, AGDSW-slab, and different variants of LAGDSW in both approaches and also add results for one-level RASPEN, that is, nonlinear Schwarz without coarse space. We consider always an additive coupling between levels as well as recycling of $\Phi$. Let us note that a hybrid coupling approach can sometimes reduce the number of nonlinear iterations more drastically, but is more sensitive to the right choice of the coarse basis functions, especially for hard problems with complicated coefficient distributions as considered here. Consequently and for brevity, we only consider the robust additive approach.
For reproducibility, we provide the list of parameters here. We use a relative stopping tolerance of $1e-5$ for Newton's method in all outer loops considering the non-preconditioned residual with a maximum number of 20 iterations. GMRES stops if a relative residual reduction of $1e-8$ is reached, with a maximum number of $2\,000$ iterations and a restart every $100$ iterations. All inner nonlinear solves in nonlinear Schwarz (local and coarse) use a relative stopping tolerance of $1e-6$ and an absolute one of $1e-8$. We vary the overlap in the Schwarz methods between $\delta=2$ and $\delta=4$. For stability, we use a simple backtracking enforcing the Armijo condition in the local nonlinear solves, but we do not use any globalization strategy in the coarse or outer solves, neither in nonlinear Schwarz nor in Newton-Krylov-Schwarz. Both approaches can be made more robust by using, e.g., a cubic backtracking, but we want to focus on the coarse spaces and their effect on the nonlinear and linear convergence, especially if the learned one can compete with AGDSW-slab.

\begin{figure}
\centering
\includegraphics[width=1.0\textwidth]{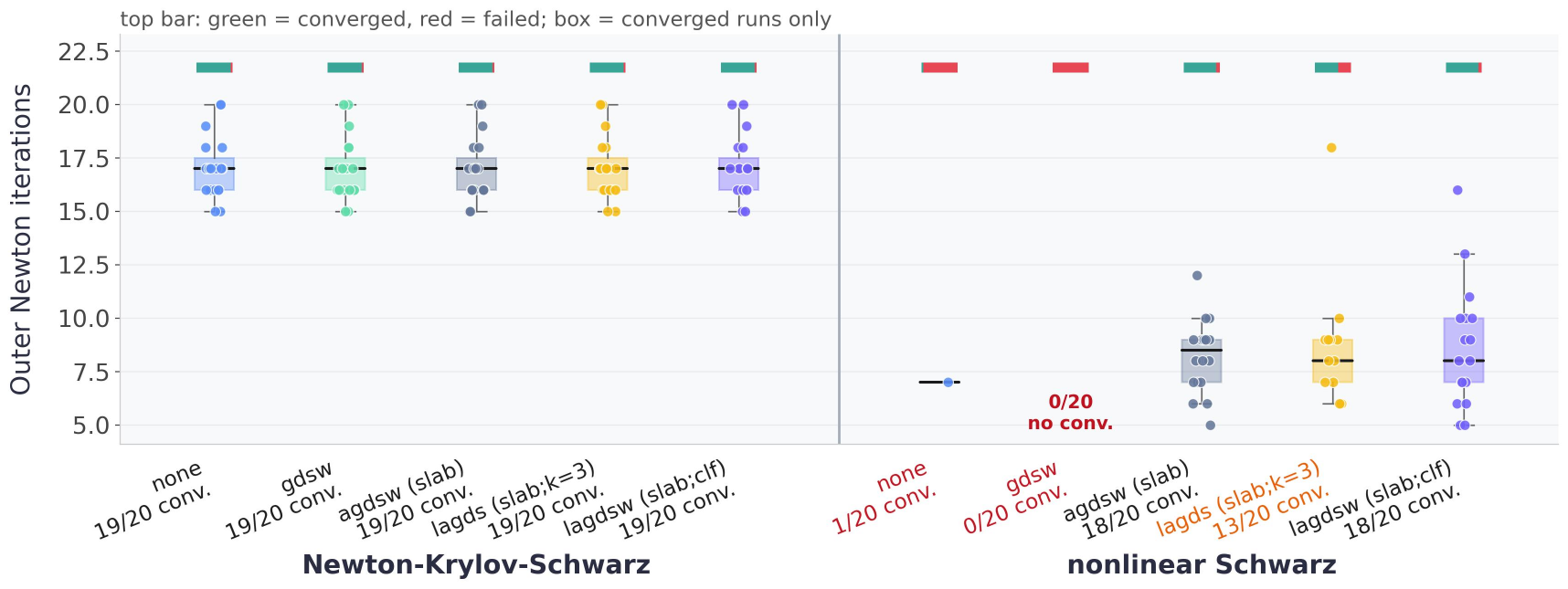}
\includegraphics[width=1.0\textwidth]{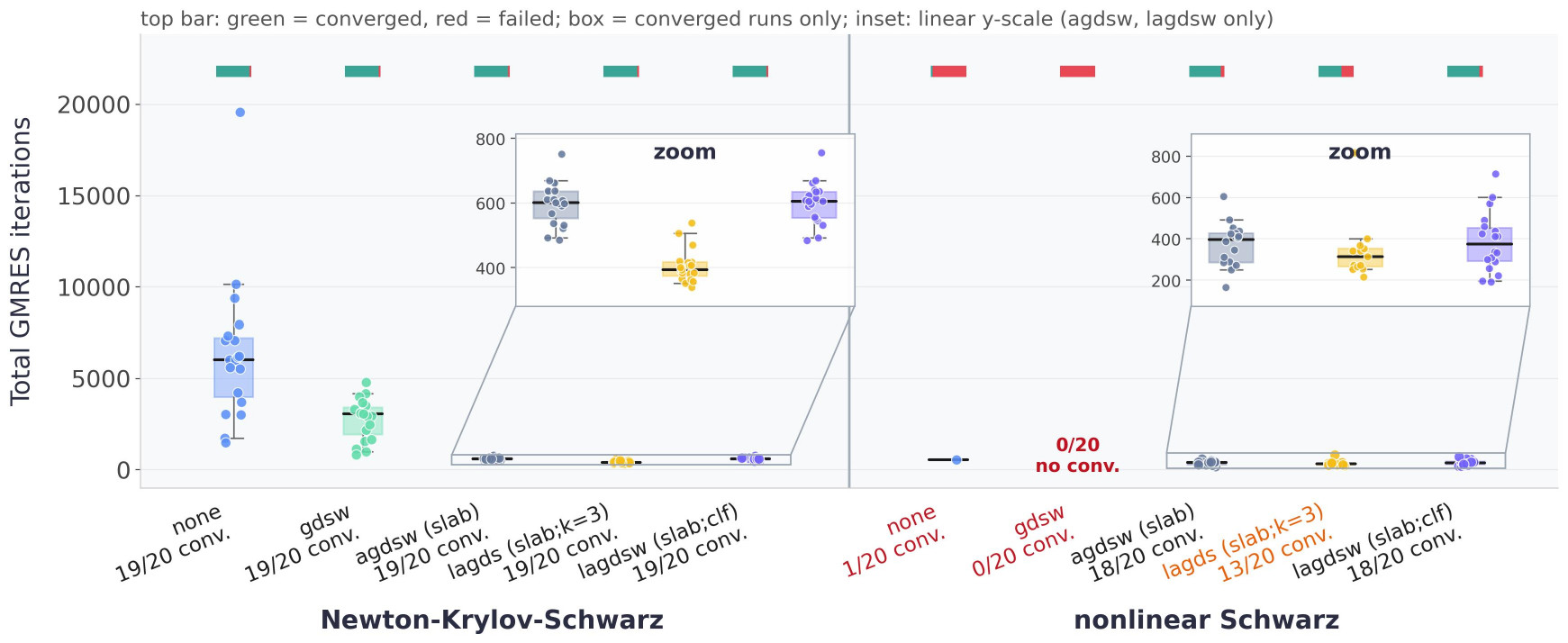}
\caption{Comparison of Newton-Krylov-Schwarz and nonlinear Schwarz methods with different coarse spaces for 20 randomly generated coefficient distributions; Newton iterations {\bf (top)} and total number GMRES iterations {\bf (bottom)}; high coefficient of $\rho=1e6$ and overlap of $\delta=2$.}
\label{fig:62}	
\end{figure}

\begin{figure}
\centering
\includegraphics[width=1.0\textwidth]{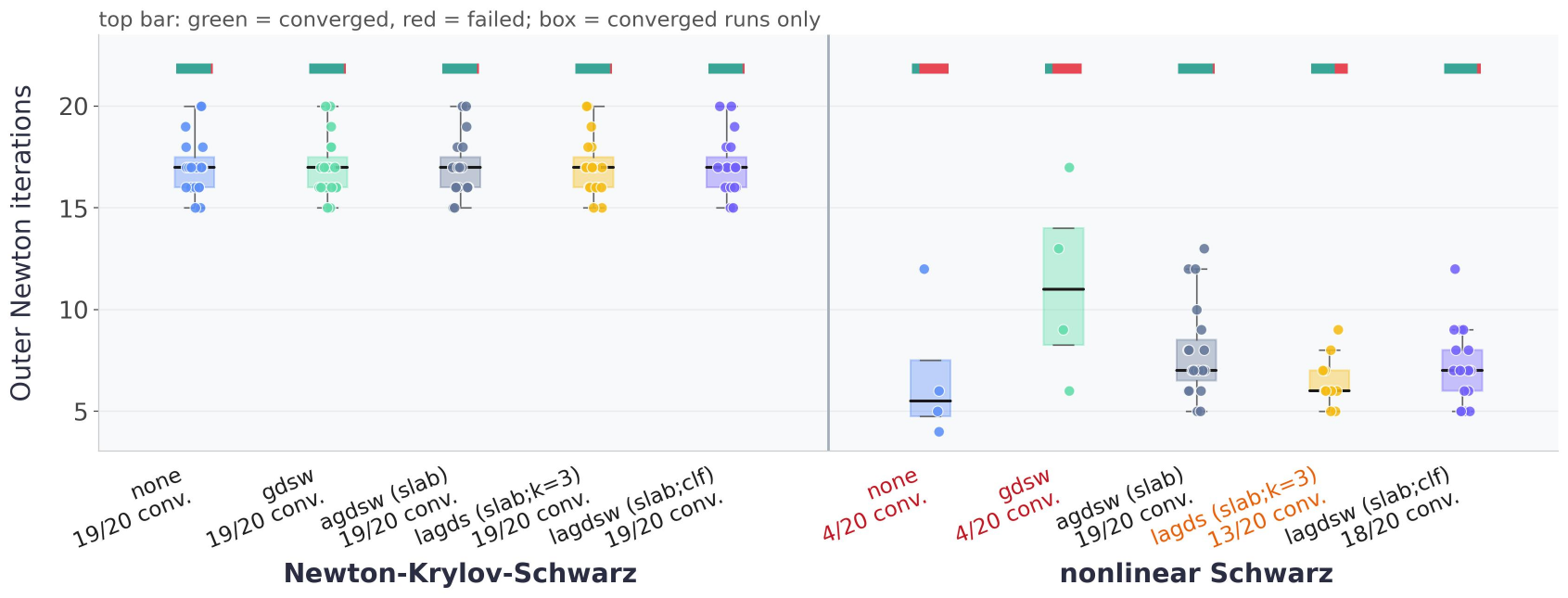}
\includegraphics[width=1.0\textwidth]{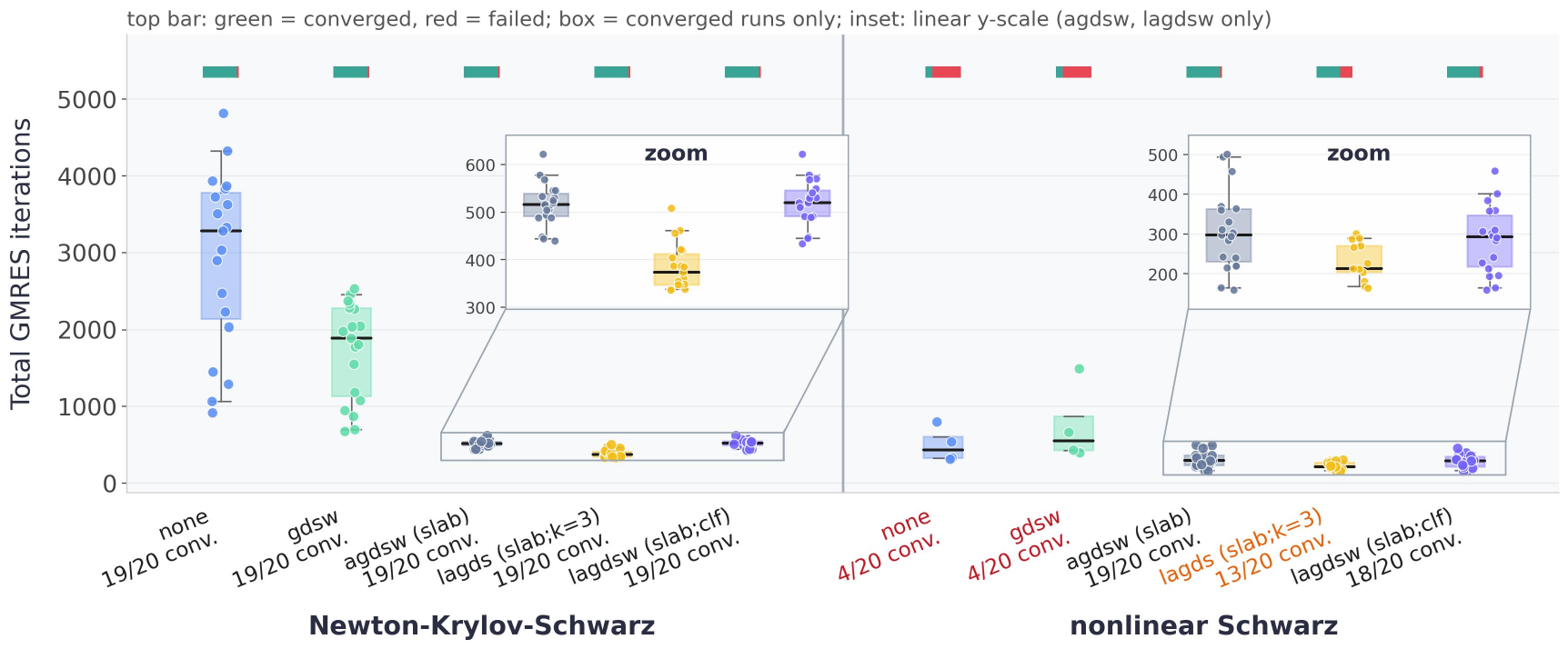}
\caption{Comparison of Newton-Krylov-Schwarz and nonlinear Schwarz methods with different coarse spaces for 20 randomly generated coefficient distributions; Newton iterations {\bf (top)} and total number GMRES iterations {\bf (bottom)}; high coefficient of $\rho=1e6$ and overlap of $\delta=4$.}
\label{fig:64}	
\end{figure}

We present the results for Newton iterations and GMRES iterations in Figures~\ref{fig:42},~\ref{fig:44},~\ref{fig:62}, and~\ref{fig:64}. We can conclude that for Newton-Krylov-Schwarz, the coarse space does not effect the nonlinear convergence and that the learned approach LAGDSW with  classifier delivers comparable GMRES iteration counts to AGDSW-slab. It can be used without any concern and outperforms GDSW by far. While one-level RASPEN (no coarse space) and GDSW are not sufficient to guarantee convergence of the nonlinear Schwarz approach, it is as reliable as Newton-Krylov-Schwarz when using an AGDSW-slab coarse space. Additionally, the number of outer Newton iterations is significantly lower in most cases. The same effect can be seen for LAGDSW with classifier and again, the learned constraints can easily be used also in nonlinear Schwarz methods. Let us note that LAGDSW with three constraints per edge is less reliable and often fails to converge, probably caused by a coarse space which is too large and which does not capture the most important nonlinear features precisely enough. To summarize, nonlinear Schwarz with adaptive coarse spaces can converge reliably and fast. The expensive adaptive coarse space can easily be replaced by our learned constraints, without retraining with data obtained from linearized $p$-Laplace problems.

\clearpage

 \section{Refined Models and Training Data}
 As described, the coefficient distribution for the training data is generated on a mesh of finite elements with $H/h=10$ and then mapped to a finer one with $H/h=20$, where the adaptive eigenvalue problems are solved. The sampling grid is even finer with 40 sampling points in each row. That actually means that  $4 \times 4$ sampling points always have the same value for the current training data. Here, we loose some potential of our approach and the models are not able to predict constraints for distributions with finer details. Therefore, in all numerical results so far, we made sure that the resolution of all coefficient distributions is not finer than the training data, since our models did not have the chance to learn finer details. To show the further potential and investigate the performance of our approach for finer coefficient distributions, we created an additional larger training data set resolving finer details. Since we are not interested in having smaller or more inclusions with high coefficients and instead just want to obtain more fine-grained and more realistic   distributions with more details and possibly frayed boundaries, we decided to take the current generation process and modify the resulting distributions slightly in a second step as follows. After mapping the coarse structure of the training data to the mesh with $H/h=20$, we randomly choose fine pixels at the boundary of regions with high coefficients and set them to the low coefficient. As a result, we keep the global structure of the distribution but obtain fine-grained, frayed boundary shapes of the inclusions. An example of refined training data is shown in Fig.~\ref{fig:rf}. Let us note that our sampling approach with 40 points in a row could still handle further refinements and resolve even more details in the future, if necessary.
 
 With the new, refined smart random data, we retrained all three regression models for constraints 2,3, and 4 and also the classifier. While previous results were performed on coarsened microsections of, for example, dual-phase steels, we now use the original data and compare the new refined models and the old ones from the previous sections. For a visualization how the microsection was coarsened so far, we refer to Fig.~\ref{fig:rf}. 
 
\begin{figure}
\centering
\includegraphics[width=0.7\textwidth]{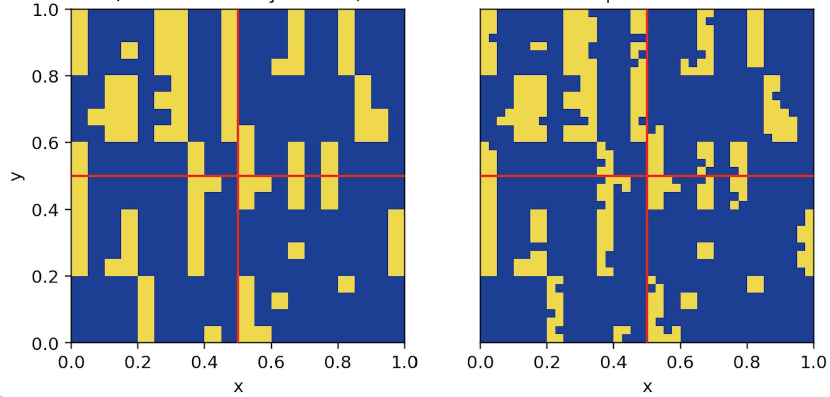}
\caption{{\bf Left:} Training data generated on a mesh with $H/h=10$ and mapped to $H/h=20$. {\bf Right:} Same data with randomly generated refinement at the boundary of high coefficient areas.}
\label{fig:rf}	
\end{figure} 

\begin{figure}
\centering
\includegraphics[width=0.7\textwidth]{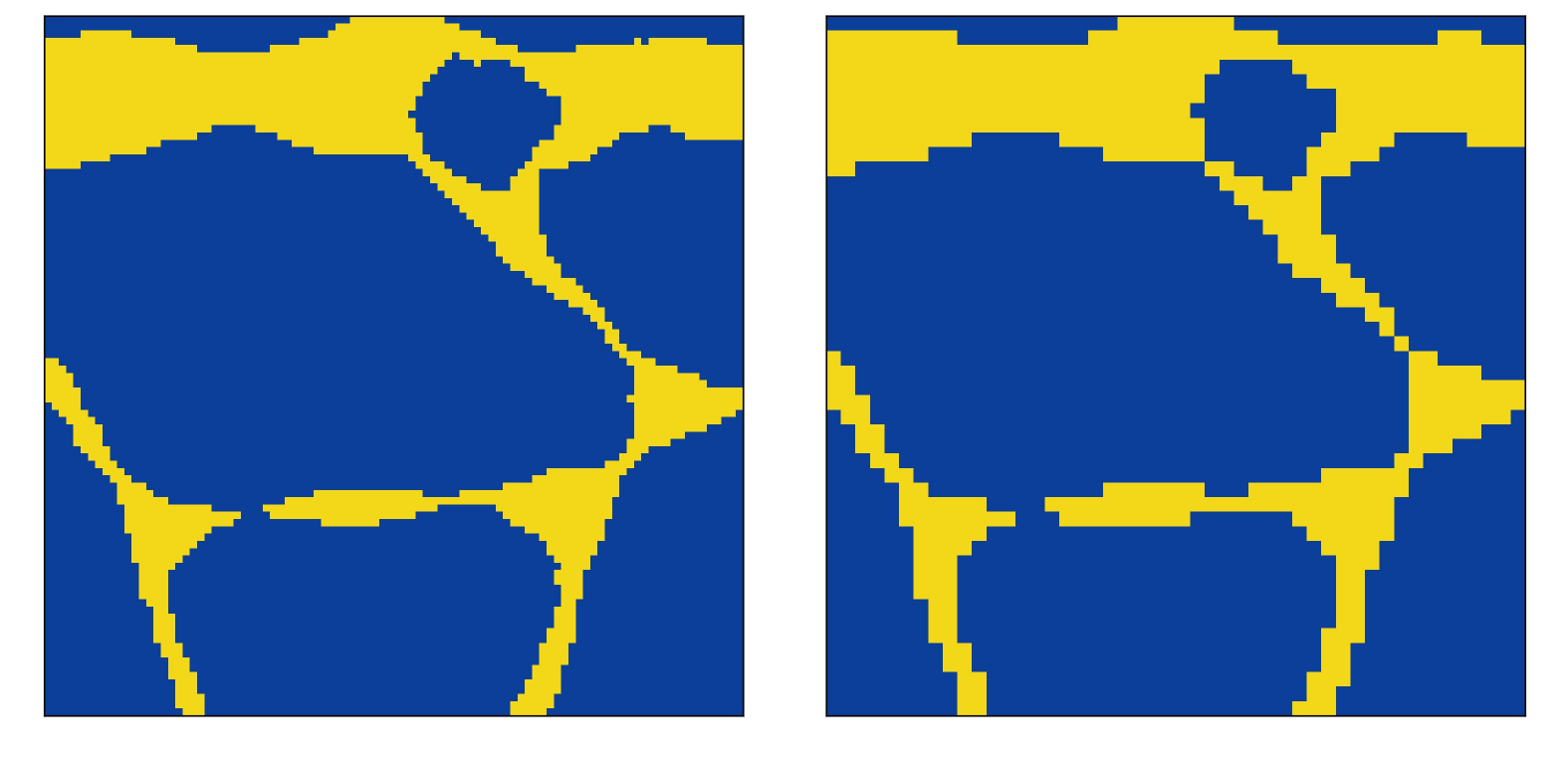}
\caption{{\bf Left:} Part of the original DP steel microsection. {\bf Right:} Coarsening using a majority voting on blocks of $2 \times 2$ pixels is applied to match the resolution of non-refined smart training data.}
\label{fig:rf}	
\end{figure}

\begin{figure}[h!]
\centering
\begin{minipage}[c]{0.29\textwidth}
    \centering
    \includegraphics[width=\linewidth]{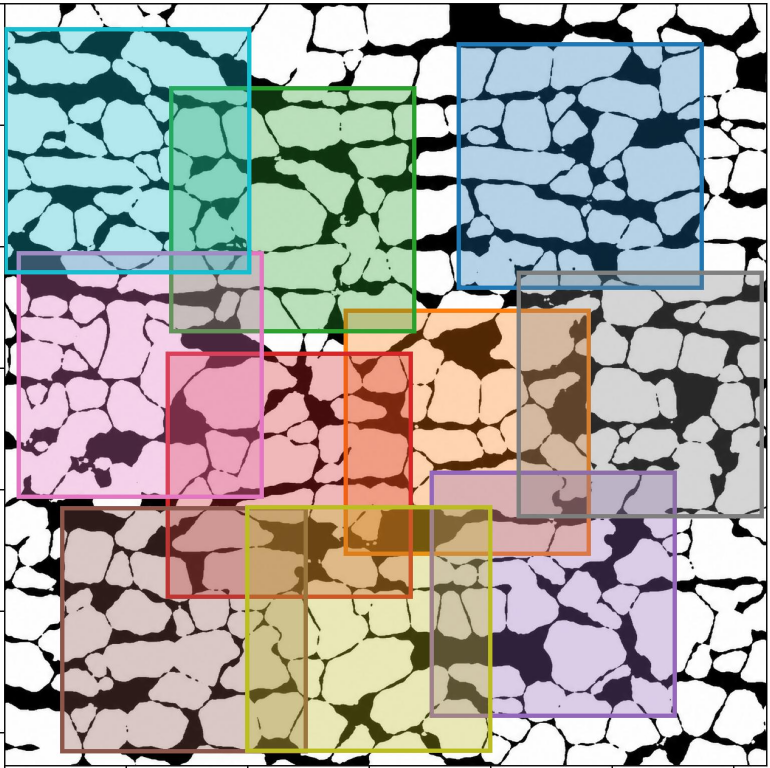}
\end{minipage}
\hfill
\begin{minipage}[c]{0.7\textwidth}
    \centering
    \includegraphics[width=\linewidth]{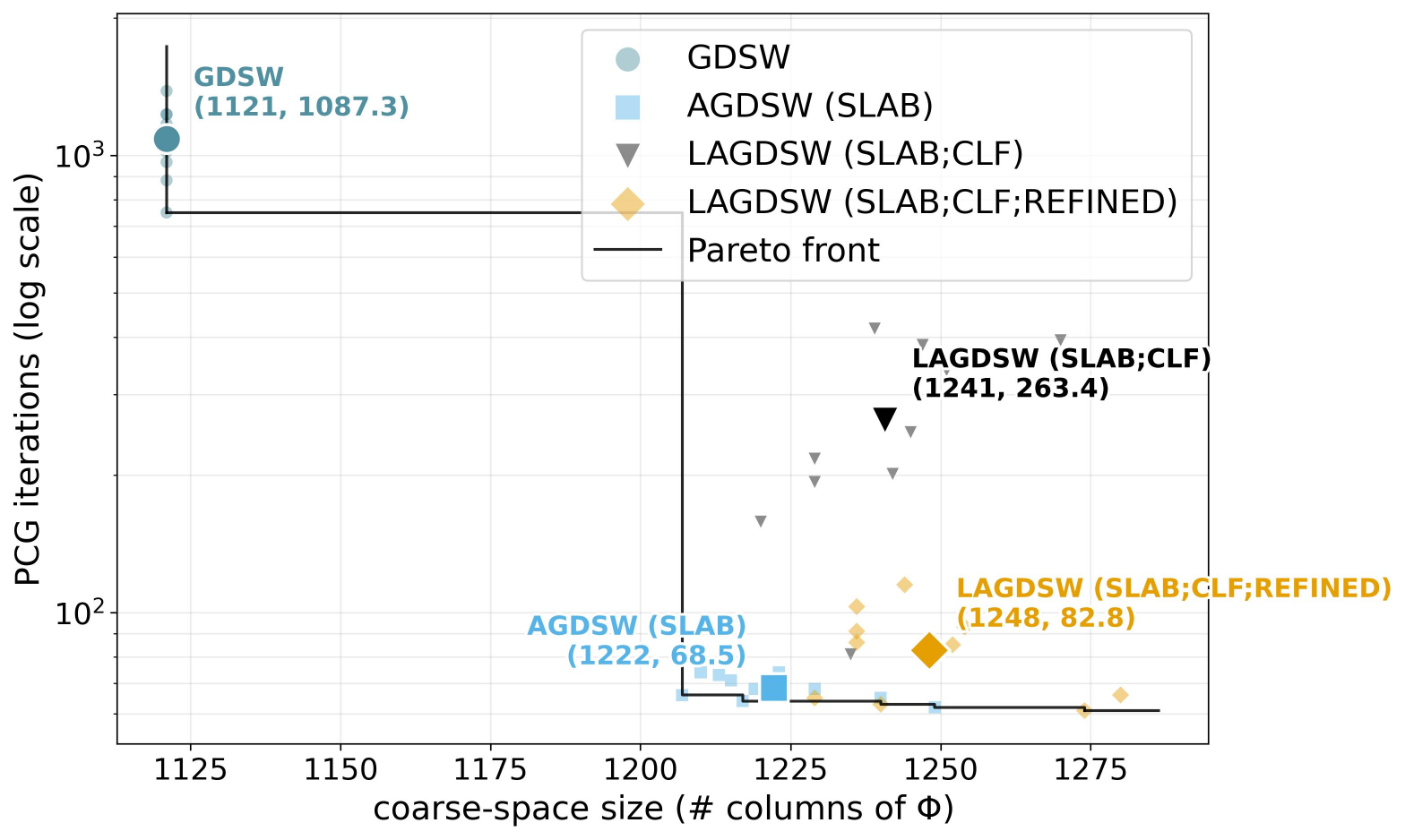}
\end{minipage}
\caption{{\bf Pareto front plot (stationary diffusion):} Coarse space size versus PCG iterations; comparison of different coarse spaces for 10 different subsections of a dual-phase steel microstructure (see left image); stationary diffusion with 400 subdomains each; $\rho_{\rm low} =1$ and $\rho_{\rm high} =1e6$. While the $x$-axis represents the size of the coarse space, the $y$-axis shows the number of PCG iterations; coarse spaces: robust AGDSW-slab, LAGDSW-slab with classifier, trained on smart data or refined smart data. The mean values of the results over the 10 samples are visualized with the larger markers.}
\label{fig:ref1}
\end{figure}

\begin{figure}[h!]
\centering
\begin{minipage}[c]{0.29\textwidth}
    \centering
    \includegraphics[width=\linewidth]{pictures2/dp_steel_fine.pdf}
\end{minipage}
\hfill
\begin{minipage}[c]{0.7\textwidth}
    \centering
    \includegraphics[width=\linewidth]{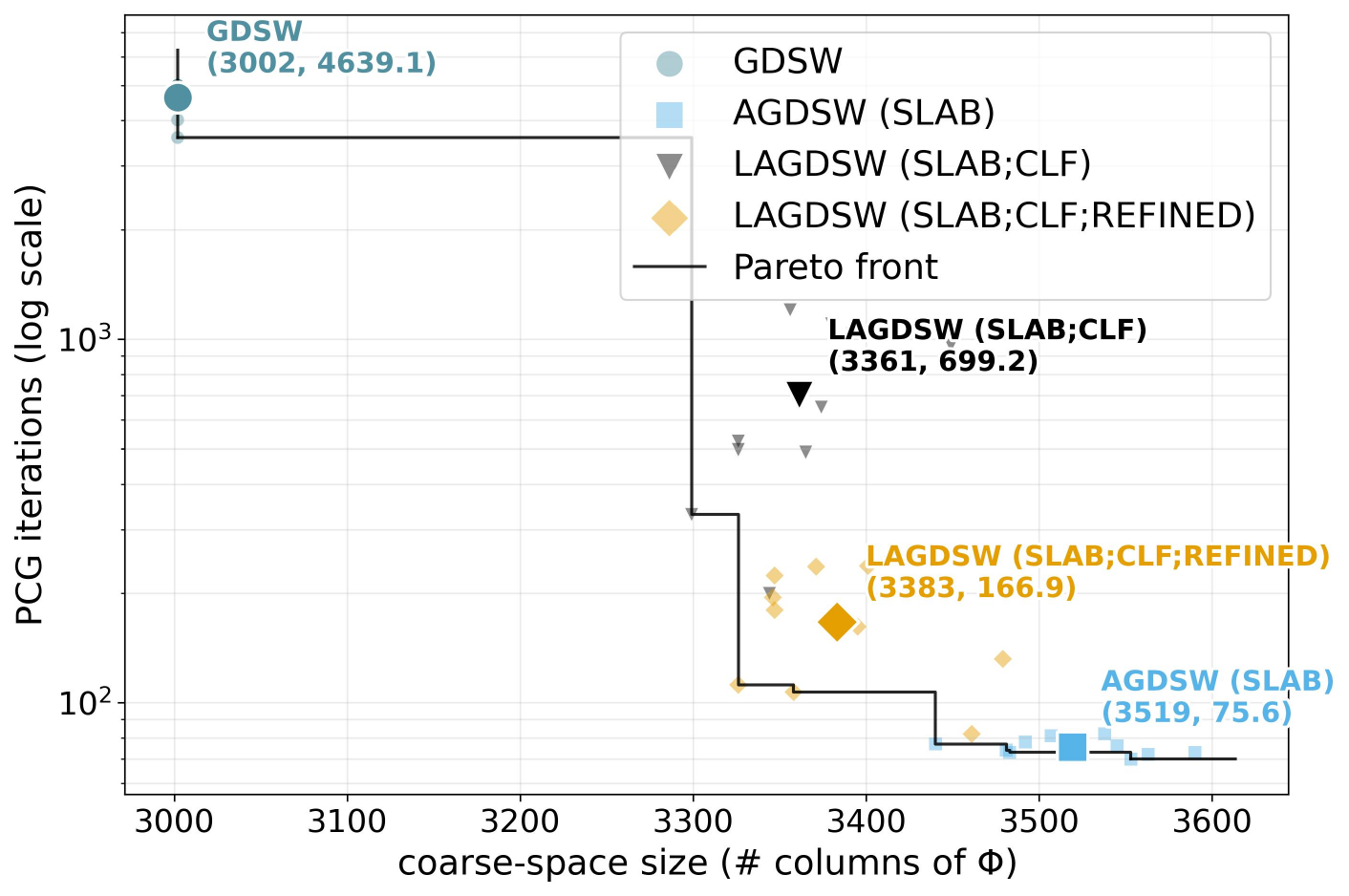}
\end{minipage}
\caption{{\bf Pareto front plot (linear elasticity):} Coarse space size versus PCG iterations; comparison of different coarse spaces for 10 different subsections of a dual-phase steel microstructure (see left image); linear elasticity with 400 subdomains each; $E_{\rm low} =210.0$ and $E_{\rm high} =2.1e8$. While the $x$-axis represents the size of the coarse space, the $y$-axis shows the number of PCG iterations; coarse spaces: robust AGDSW-slab, LAGDSW-slab with classifier, trained on smart data or refined smart data. The mean values of the results over the 10 samples are visualized with the larger markers.}
\label{fig:ref2}
\end{figure}

\begin{figure}[h!]
\centering
\begin{minipage}[c]{0.29\textwidth}
    \centering
    \includegraphics[width=\linewidth]{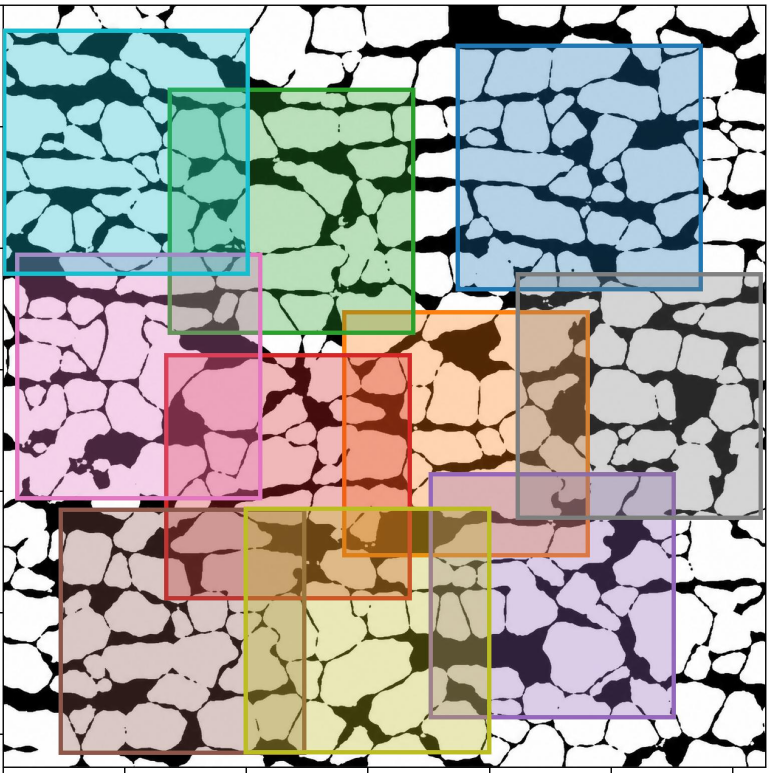}
\end{minipage}
\hfill
\begin{minipage}[c]{0.7\textwidth}
    \centering
    \includegraphics[width=\linewidth]{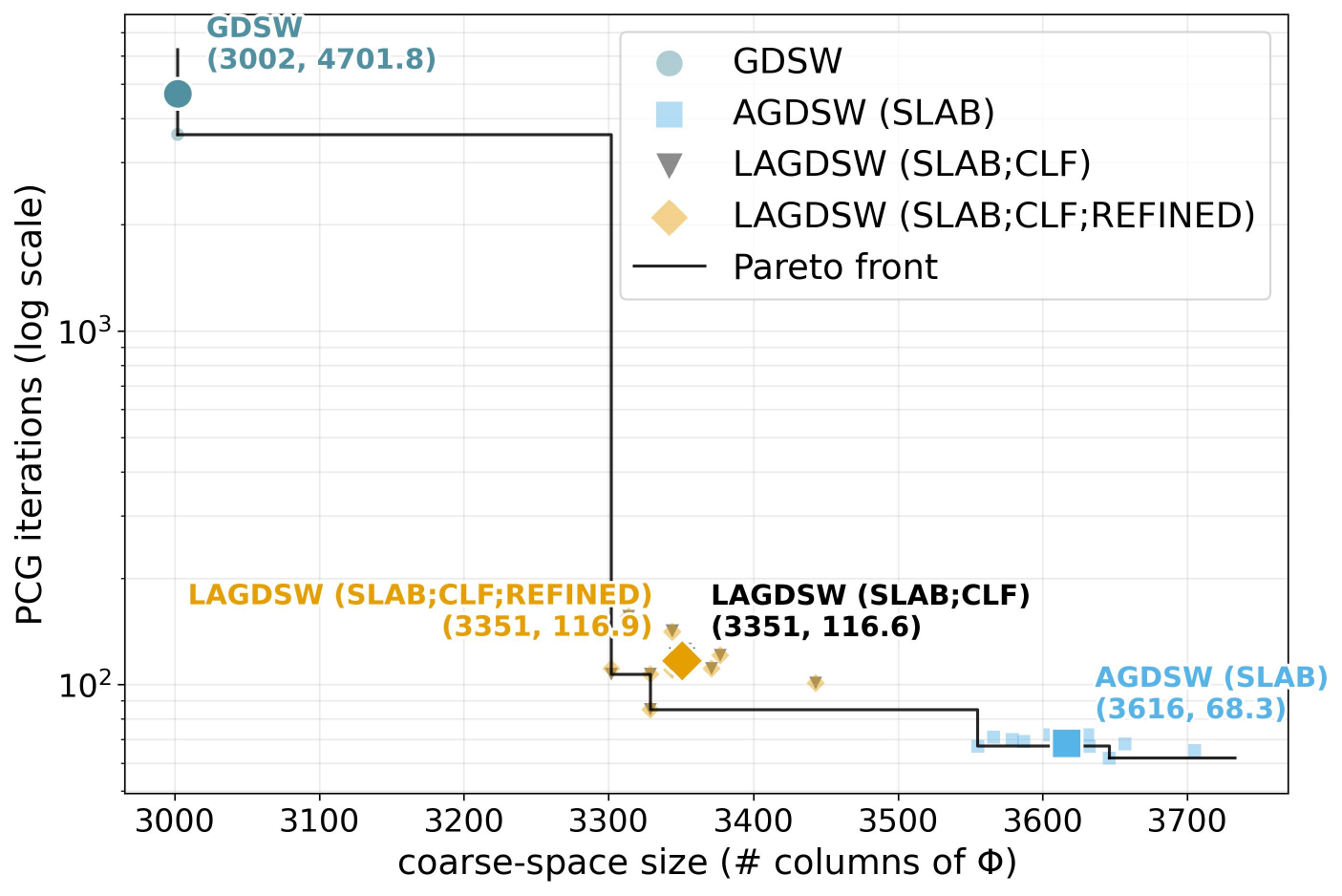}
\end{minipage}
\caption{{\bf Pareto front plot (linear elasticity):} Coarse space size versus PCG iterations; comparison of different coarse spaces for 10 different subsections of a dual-phase steel microstructure, which is coarsened by a $2 \times 2$ pixel majority voting (see left image); linear elasticity with 400 subdomains each; $E_{\rm low} =210.0$ and $E_{\rm high} =2.1e8$. While the $x$-axis represents the size of the coarse space, the $y$-axis shows the number of PCG iterations; coarse spaces: robust AGDSW-slab, LAGDSW-slab with classifier, trained on smart data or refined smart data. The means on the results over the 10 samples are visualized with the larger markers.}
\label{fig:ref3}
\end{figure}

We consider both, linear diffusion and linear elasticity for 10 different subsections of the dual-phase steel coefficient distribution with a coefficient jump of one million, either in $\rho$ or in $E$. In both cases, the refined models perform significantly better than the old models; see Fig.~\ref{fig:ref1} and Fig.~\ref{fig:ref2}. Let us remark that in case of elasticity, the learned coarse spaces are smaller than the adaptive one. As before, not all adaptive constraints can be interpreted as a lifted scalar adaptive constraint. Still, the iteration counts for the LAGDSW approach with refined training data are reasonably small especially compared to GDSW. As a final sanity check, we run the same elasticity tests with the coarsened DP steel microstructure, where both models, the old and the refined one, perform equally well; see Fig.~\ref{fig:ref3} for the results.\\

{\bf Use of Artificial Intelligence Tools.} 
 All authors used artificial intelligence (AI) tools 
for editorial support during the final preparation of the manuscript, including the identification of typographical errors, consistency checks, and linguistic suggestions. All AI-generated suggestions were critically reviewed and approved by all authors. 

AI tools for software development were used only by the author Martin Lanser and only for specific software engineering tasks. In particular, AI assistance was used to develop a Python-based software test environment for linear and nonlinear Schwarz methods with various coarse spaces. This implementation was derived from an existing MATLAB code base that had been developed entirely without AI assistance in recent years. In addition, Martin Lanser extended the resulting Python software by implementing several additional modules and functionalities beyond the original MATLAB implementation, especially modules for the visualization of results. AI assistance was also used during the development of the software pipeline for the generation of the training data used in this work.

All scientific aspects of the research, including the formulation of the research questions, methodological decisions, mathematical derivations, experimental design, analysis and interpretation of the results, and the conclusions, were carried out without AI assistance. The remaining authors did not use AI tools for any research-related activities.

{\bf Acknowledgements.} We gratefully acknowledge the use of the computational
facilities of the Center for Data and Simulation Science (CDS) at the University of
Cologne.


\bibliographystyle{siamplain}
\bibliography{feti,ml,BDDC_bib,agdsw} 

@article {robust,
    AUTHOR = {Klawonn, Axel and Rheinbach, Oliver},
     TITLE = {Robust {FETI}-{DP} methods for heterogeneous three dimensional
              elasticity problems},
   JOURNAL = {Comput. Methods Appl. Mech. Engrg.},
  FJOURNAL = {Computer Methods in Applied Mechanics and Engineering},
    VOLUME = {196},
      YEAR = {2007},
    NUMBER = {8},
     PAGES = {1400--1414},
      ISSN = {0045-7825},
   MRCLASS = {74S05 (65N55 74E05 74G15)},
  MRNUMBER = {2277025},
       DOI = {10.1016/j.cma.2006.03.023},
       URL = {https://doi.org/10.1016/j.cma.2006.03.023},
}

@article{gdsw1,
  author  = {Clark R. Dohrmann and Axel Klawonn and Olof B. Widlund},
  title   = {Domain Decomposition for Less Regular Subdomains:
             Overlapping {S}chwarz in Two Dimensions},
  journal = {SIAM Journal on Numerical Analysis},
  volume  = {46},
  number  = {4},
  pages   = {2153--2168},
  year    = {2008},
  doi     = {10.1137/070685841}
}

@article{agdsw1,
	author = {Heinlein, Alexander and Klawonn, Axel and Knepper, Jascha and Rheinbach, Oliver},
	title = {{A}daptive {GDSW} {C}oarse {S}paces for {O}verlapping {S}chwarz {M}ethods in {T}hree {D}imensions},
	journal = {SIAM Journal on Scientific Computing},
	volume = {41},
	number = {5},
	pages = {A3045-A3072},
	year = {2019},
}

@incollection{agdsw2,
  author    = {Alexander Heinlein and Axel Klawonn and
               Jascha Knepper and Oliver Rheinbach},
  title     = {An Adaptive {GDSW} Coarse Space for Two-Level Overlapping
               {S}chwarz Methods in Two Dimensions},
  booktitle = {Domain Decomposition Methods in Science and Engineering XXIV},
  series    = {Lecture Notes in Computational Science and Engineering},
  volume    = {125},
  pages     = {373--382},
  publisher = {Springer},
  year      = {2018},
  doi       = {10.1007/978-3-319-93873-8_37}
}

@phdthesis{Knepper2022,
  author       = {Jascha Knepper},
  title        = {Adaptive Coarse Spaces for the Overlapping Schwarz Method and Multiscale Elliptic Problems},
  school       = {Universit{\"a}t zu K{\"o}ln},
  year         = {2022},
  type         = {PhD thesis},
  url          = {https://kups.ub.uni-koeln.de/62002/},
  note         = {Open Access Dissertation, University of Cologne}
}

@Book{toselli_widlund,
  author	= "Andrea Toselli and Olof Widlund",
  title		= "Domain Decomposition Methods - Algorithms and Theory",
  publisher	= "Springer",
  year		= 2004,
  volume	= 34,
  series	= "Springer Series in Computational Mathematics"
}

@article {kl_p1,
    AUTHOR = {Klawonn, Axel and Lanser, Martin and Rheinbach, Oliver},
     TITLE = {Toward extremely scalable nonlinear domain decomposition
              methods for elliptic partial differential equations},
   JOURNAL = {SIAM J. Sci. Comput.},
  FJOURNAL = {SIAM Journal on Scientific Computing},
    VOLUME = {37},
      YEAR = {2015},
    NUMBER = {6},
     PAGES = {C667--C696},
      ISSN = {1064-8275,1095-7197},
   MRCLASS = {65N30 (65F08 65F10 65N55 65Y05)},
  MRNUMBER = {3432150},
MRREVIEWER = {Marius\ Ghergu},
       url = {https://doi.org/10.1137/140997907},
}

@article {kl_p3,
    AUTHOR = {Klawonn, Axel and Rheinbach, Oliver},
     TITLE = {A parallel implementation of dual-primal {FETI} methods for
              three-dimensional linear elasticity using a transformation of
              basis},
   JOURNAL = {SIAM J. Sci. Comput.},
  FJOURNAL = {SIAM Journal on Scientific Computing},
    VOLUME = {28},
      YEAR = {2006},
    NUMBER = {5},
     PAGES = {1886--1906},
      ISSN = {1064-8275,1095-7197},
   MRCLASS = {74S05 (65N30 74B05 74G15)},
  MRNUMBER = {2272193},
       url = {https://doi.org/10.1137/050624364},
}

@article {kl_p4,
    AUTHOR = {Klawonn, Axel and Lanser, Martin and Rheinbach, Oliver and
              Uran, Matthias},
     TITLE = {Nonlinear {FETI}-{DP} and {BDDC} methods: a unified framework
              and parallel results},
   JOURNAL = {SIAM J. Sci. Comput.},
  FJOURNAL = {SIAM Journal on Scientific Computing},
    VOLUME = {39},
      YEAR = {2017},
    NUMBER = {6},
     PAGES = {C417--C451},
      ISSN = {1064-8275,1095-7197},
   MRCLASS = {65N55 (65F08 65N22 65N30 65Y05 68U20 68W10)},
  MRNUMBER = {3724259},
MRREVIEWER = {Christos\ Kravvaritis},
       url = {https://doi.org/10.1137/16M1102495},
}

@article{fetidp1,
  author  = {Charbel Farhat and Michel Lesoinne and Kendall Pierson},
  title   = {A Scalable Dual-Primal Domain Decomposition Method},
  journal = {Numerical Linear Algebra with Applications},
  volume  = {7},
  number  = {7--8},
  pages   = {687--714},
  year    = {2000},
  doi     = {10.1002/1099-1506(200010/12)7:7/8<687::AID-NLA219>3.0.CO;2-0}
}

@article{fetidp2,
  author  = {Charbel Farhat and Michel Lesoinne and Patrick Le Tallec
             and Kendall Pierson and Daniel Rixen},
  title   = {{FETI-DP}: A Dual-Primal Unified {FETI} Method -- Part {I}:
             A Faster Alternative to the Two-Level {FETI} Method},
  journal = {International Journal for Numerical Methods in Engineering},
  volume  = {50},
  number  = {7},
  pages   = {1523--1544},
  year    = {2001},
  doi     = {10.1002/1097-0207(20010310)50:7<1523::AID-NME167>3.0.CO;2-B}
}

@article{fetidp3,
  author  = {Axel Klawonn and Maksymilian Dryja and Olof B. Widlund},
  title   = {Dual-Primal {FETI} Methods for Three-Dimensional Elliptic Problems with Heterogeneous Coefficients},
  journal = {SIAM Journal on Numerical Analysis},
  volume  = {40},
  number  = {1},
  pages   = {159--179},
  year    = {2002},
  doi     = {10.1137/S0036142900365982}
}

@article{fetidp4,
  author  = {Axel Klawonn and Olof B. Widlund},
  title   = {Dual-Primal {FETI} Methods for Linear Elasticity},
  journal = {Communications on Pure and Applied Mathematics},
  volume  = {59},
  number  = {11},
  pages   = {1523--1572},
  year    = {2006},
  doi     = {10.1002/cpa.20124}
}

@article{bddc1,
  author  = {Clark R. Dohrmann},
  title   = {A Preconditioner for Substructuring Based on Constrained
             Energy Minimization},
  journal = {SIAM Journal on Scientific Computing},
  volume  = {25},
  number  = {1},
  pages   = {246--258},
  year    = {2003},
  doi     = {10.1137/S1064827502412887}
}

@article{bddc2,
  author  = {Jan Mandel and Clark R. Dohrmann},
  title   = {Convergence of a Balancing Domain Decomposition by
             Constraints and Energy Minimization},
  journal = {Numerical Linear Algebra with Applications},
  volume  = {10},
  number  = {7},
  pages   = {639--659},
  year    = {2003},
  doi     = {10.1002/nla.338}
}

@inproceedings{ZampiniTuLi2016,
  author    = {Stefano Zampini and Xiao-Chuan Tu and Jed Brown and
               Barry F. Smith and Hong Zhang and Jianjun Li},
  title     = {{PCBDDC}: A Class of Robust Dual-Primal Methods in {PETSc}},
  booktitle = {Domain Decomposition Methods in Science and Engineering XXIII},
  series    = {Lecture Notes in Computational Science and Engineering},
  volume    = {116},
  pages     = {535--543},
  publisher = {Springer},
  year      = {2017},
  doi       = {10.1007/978-3-319-52389-7_53}
}

@article{frosch1,
  author  = {Alexander Heinlein and Axel Klawonn and Oliver Rheinbach},
  title   = {A Parallel Implementation of a Two-Level Overlapping {S}chwarz Method with Energy-Minimizing Coarse Space Based on {T}rilinos},
  journal = {SIAM Journal on Scientific Computing},
  volume  = {38},
  number  = {6},
  pages   = {C713--C747},
  year    = {2016},
  doi     = {10.1137/16M1062843}
}

@article{SousedikSistekMandel2013,
  author  = {Bed{\v{r}}ich Soused{\'\i}k and Jakub {\v{S}}{\'\i}stek and Jan Mandel},
  title   = {Adaptive-Multilevel {BDDC} and its Parallel Implementation},
  journal = {Computing},
  volume  = {95},
  number  = {12},
  pages   = {1087--1119},
  year    = {2013},
  doi     = {10.1007/s00607-013-0299-5}
}

@article{geneo,
  author  = {Nicole Spillane and Victor Dolean and Pierre-Henri Tournier
             and Fr{\'e}d{\'e}ric Nataf and Cyril Jolivet and
             Haidar Khamraev},
  title   = {Abstract Robust Coarse Spaces for Systems of {PDE}s via Generalized Eigenproblems in the Overlaps},
  journal = {Numerische Mathematik},
  volume  = {126},
  number  = {4},
  pages   = {741--770},
  year    = {2014},
  doi     = {10.1007/s00211-013-0576-y}
}

@article{shem,
  author  = {Petter E. Bj{\o}rstad and Martin J. Gander and
             Andreas Loneland and Talal Rahman},
  title   = {Analysis of a Spectral Harmonically Enriched Multiscale Coarse Space
             for Two-Level {S}chwarz Methods},
  journal = {SIAM Journal on Scientific Computing},
  volume  = {40},
  number  = {4},
  pages   = {A2360--A2385},
  year    = {2018},
  doi     = {10.1137/17M114348X}
}

@incollection{MandelSousedik2007,
  author    = {Jan Mandel and Bed{\v{r}}ich Soused{\'i}k},
  title     = {Adaptive Coarse Space Selection in the {BDDC} and the {FETI-DP} Iterative Substructuring Methods},
  booktitle = {Domain Decomposition Methods in Science and Engineering XVII},
  series    = {Lecture Notes in Computational Science and Engineering},
  volume    = {60},
  pages     = {93--101},
  publisher = {Springer},
  year      = {2008}
}

@article{PechsteinDohrmann2017,
  author  = {Clemens Pechstein and Clark R. Dohrmann},
  title   = {A Unified Framework for Adaptive {BDDC}},
  journal = {Electronic Transactions on Numerical Analysis},
  volume  = {46},
  pages   = {273--336},
  year    = {2017},
  url     = {https://etna.math.kent.edu/volumes/2011-2020/vol46/abstract.php?vol=46&pages=273-336}
}

@article{KlawonnRadtkeRheinbach2016,
  author  = {Axel Klawonn and Patrick Radtke and Oliver Rheinbach},
  title   = {A Comparison of Adaptive Coarse Spaces for Iterative Substructuring Methods in Two Dimensions},
  journal = {Electronic Transactions on Numerical Analysis},
  volume  = {45},
  pages   = {75--106},
  year    = {2016},
  url     = {https://etna.math.kent.edu/volumes/2011-2020/vol45/abstract.php?vol=45&pages=75-106}
}

@article{nl_left1,
  author  = {Xiao-Chuan Cai and David E. Keyes},
  title   = {Nonlinearly Preconditioned Inexact {N}ewton Algorithms},
  journal = {SIAM Journal on Scientific Computing},
  volume  = {24},
  number  = {1},
  pages   = {183--200},
  year    = {2002},
  doi     = {10.1137/S1064827500361235}
}

@article{nl_left2,
  author  = {Victor Dolean and Martin J. Gander and Walid Kheriji
             and Felix Kwok and Roland Masson},
  title   = {How to Use a Nonlinear {S}chwarz Method to Precondition {N}ewton's Method},
  journal = {SIAM Journal on Scientific Computing},
  volume  = {38},
  number  = {6},
  pages   = {A3357--A3380},
  year    = {2016},
  doi     = {10.1137/15M102887X}
}

@article{nl_right1,
  author  = {Axel Klawonn and Martin Lanser and Oliver Rheinbach},
  title   = {Nonlinear {FETI-DP} and {BDDC} Methods},
  journal = {SIAM Journal on Scientific Computing},
  volume  = {36},
  number  = {2},
  pages   = {A737--A765},
  year    = {2014},
  doi     = {10.1137/130920563}
  }

@article{nl_right2,
  author  = {Julien Pebrel and Christian Rey and Pierre Gosselet},
  title   = {A Nonlinear Dual Domain Decomposition Method:
             Application to Structural Problems with Damage},
  journal = {International Journal for Multiscale Computational Engineering},
  volume  = {6},
  number  = {3},
  pages   = {251--262},
  year    = {2008},
  doi     = {10.1615/IntJMultCompEng.v6.i3.40}
}

@article{nl_left3,
  author  = {Alexander Heinlein and Martin Lanser},
  title   = {Additive and Hybrid Nonlinear Two-Level {S}chwarz Methods and
             Energy Minimizing Coarse Spaces for Unstructured Grids},
  journal = {SIAM Journal on Scientific Computing},
  volume  = {42},
  number  = {4},
  pages   = {A2461--A2488},
  year    = {2020},
  doi     = {10.1137/19M1276972}
}

@article{nl_adaptive,
  author  = {Alexander Heinlein and Axel Klawonn and Martin Lanser},
  title   = {Adaptive Nonlinear Domain Decomposition Methods with an Application to the {$p$}-{L}aplacian},
  journal = {SIAM Journal on Scientific Computing},
  volume  = {45},
  number  = {3},
  pages   = {S152--S172},
  year    = {2023},
  doi     = {10.1137/21M1433605}
}

@article{HKLW:2018:ML_acc, 
    AUTHOR = {Alexander Heinlein and Axel Klawonn and Martin Lanser and Janine Weber}, 
    TITLE = {{M}achine {L}earning in {A}daptive {Do}main {D}ecomposition {M}ethods - {P}redicting the {G}eometric {L}ocation of {C}onstraints},
     JOURNAL = {SIAM J. Sci. Comput.},
  FJOURNAL = {SIAM Journal on Scientific Computing},
    VOLUME = {41},
      YEAR = {2019},
    NUMBER = {6},
     PAGES = {A3887--A3912},
}

@phdthesis{Weber:2022:Diss,
  title={Efficient and robust {FETI-DP} and {BDDC} methods -- {A}pproximate coarse spaces and deep learning-based adaptive coarse space},
  author={Weber, Janine},
  year={2022},
  school={Universit{\"a}t zu K{\"o}ln},
  note = {Online available at \url{http://kups.ub.uni-koeln.de/id/eprint/55179}}
}

@article{KLW:2024:JCP,
title = {Learning adaptive coarse basis functions of {FETI-DP}},
journal = {Journal of Computational Physics},
volume = {496},
pages = {112587},
year = {2024},
issn = {0021-9991},
author = {Axel Klawonn and Martin Lanser and Janine Weber},
}

@incollection{HKLW:2023:DD_schwarz,
  title={Predicting the geometric location of critical edges in adaptive {GDSW} overlapping domain decomposition methods using deep learning},
  author={Heinlein, Alexander and Klawonn, Axel and Lanser, Martin and Weber, Janine},
  booktitle={Domain Decomposition Methods in Science and Engineering XXVI},
  pages={307--315},
  year={2023},
  publisher={Springer}
}

@inproceedings{KLW:2023:DD_nonlin,
  title={{L}earning {A}daptive {C}onstraints in {N}onlinear {FETI-DP} {M}ethods},
  author={Klawonn, Axel and Lanser, Martin and Weber, Janine},
  booktitle={European Conference on Numerical Mathematics and Advanced Applications},
  pages={45--53},
  year={2023},
  organization={Springer}
}

@InProceedings{Klawonn:2002:DPF_face_constraints,
author={Klawonn, Axel and Widlund, Olof B. and Dryja, Maksymilian},
editor={Pavairno, Luca F. and Toselli, Andrea},
title={Dual-{P}rimal {FETI} Methods with Face Constraints},
booktitle={Recent Developments in Domain Decomposition Methods},
year={2002},
publisher={Springer Berlin Heidelberg},
address={Berlin, Heidelberg},
pages={27--40},
}

@article{USDpt:SciML,
title = {{B}rochure on {B}asic {R}esearch {N}eeds for {S}cientific {M}achine {L}earning: {C}ore {T}echnologies for {A}rtificial {I}ntelligence. {USDOE} {O}ffice of {S}cience {(SC)} ({U}nited {S}tates)},
author = {Baker, Nathan and Alexander, Frank and Bremer, Timo and Hagberg, Aric and Kevrekidis, Yannis and Najm, Habib and Parashar, Manish and Patra, Abani and Sethian, James and Wild, Stefan and Willcox, Karen},
doi = {10.2172/1484362},
url = {https://www.osti.gov/biblio/1484362}, 
publisher = {USDOE Office of Science (SC)},
place = {United States},
year = {2018}
}

@misc{tensorflow2015-whitepaper,
title={ {TensorFlow}: Large-Scale Machine Learning on Heterogeneous Systems},
url={https://www.tensorflow.org/},
note={Software available from \url{https://www.tensorflow.org/}},
author={
    Mart\'{\i}n~Abadi and
    Ashish~Agarwal and
    Paul~Barham and
    Eugene~Brevdo and
    Zhifeng~Chen and
    Craig~Citro and
    Greg~S.~Corrado and
    Andy~Davis and
    Jeffrey~Dean and
    Matthieu~Devin and
    Sanjay~Ghemawat and
    Ian~Goodfellow and
    Andrew~Harp and
    Geoffrey~Irving and
    Michael~Isard and
    Yangqing Jia and
    Rafal~Jozefowicz and
    Lukasz~Kaiser and
    Manjunath~Kudlur and
    Josh~Levenberg and
    Dandelion~Man\'{e} and
    Rajat~Monga and
    Sherry~Moore and
    Derek~Murray and
    Chris~Olah and
    Mike~Schuster and
    Jonathon~Shlens and
    Benoit~Steiner and
    Ilya~Sutskever and
    Kunal~Talwar and
    Paul~Tucker and
    Vincent~Vanhoucke and
    Vijay~Vasudevan and
    Fernanda~Vi\'{e}gas and
    Oriol~Vinyals and
    Pete~Warden and
    Martin~Wattenberg and
    Martin~Wicke and
    Yuan~Yu and
    Xiaoqiang~Zheng},
  year={2015},
}

@article{klawonn2024survey,
  title={Machine learning and domain decomposition methods-a survey},
  author={Klawonn, Axel and Lanser, Martin and Weber, Janine},
  journal={Computational Science and Engineering},
  volume={1},
  number={1},
  pages={2},
  year={2024},
  publisher={Springer}
}

@article{chung2021learning,
  title={Learning adaptive coarse spaces of {BDDC} algorithms for stochastic elliptic problems with oscillatory and high contrast coefficients},
  author={Chung, Eric and Kim, Hyea-Hyun and Lam, Ming-Fai and Zhao, Lina},
  journal={Mathematical and Computational Applications},
  volume={26},
  number={2},
  pages={44},
  year={2021},
  publisher={MDPI}
}

@article{Dolean:2026:NO,
  title={When can a neural operator replace a coarse solve? {A}rchitectural principles for two-level preconditioning},
  author={Melchers, Hugo and Abdelmalik, Michael and Dolean, Victorita},
  journal={arXiv preprint arXiv:2605.19867},
  year={2026}
}

\end{document}